\documentclass{amsart}

\usepackage{amsmath}
\usepackage{amssymb}
\usepackage{mathtools}
\usepackage{url}
\usepackage{hyperref}
\usepackage{autonum}
\usepackage{xspace}
\usepackage{xcolor}
\usepackage{hyperref}

\usepackage{graphicx}
\usepackage{tikz}
\usepackage{pgfplots} 

\usepackage{caption}
\usepackage{subcaption}
\usepackage{comment}
\usepackage{booktabs}
\usepackage{float}

\numberwithin{equation}{section}

\hypersetup{
		colorlinks,
		linkcolor={red!50!black},
		citecolor={blue!50!black},
		urlcolor={blue!80!black}
	}

\definecolor{color0}{rgb}{0.12156862745098,0.466666666666667,0.705882352941177}
\definecolor{color1}{rgb}{1,0.498039215686275,0.0549019607843137}
\definecolor{color2}{rgb}{0.172549019607843,0.627450980392157,0.172549019607843}
\definecolor{color3}{rgb}{0.83921568627451,0.152941176470588,0.156862745098039}
\definecolor{color4}{rgb}{0.580392156862745,0.403921568627451,0.741176470588235}
\definecolor{color5}{rgb}{0.549019607843137,0.337254901960784,0.294117647058824}
\definecolor{color6}{rgb}{0.890196078431372,0.466666666666667,0.76078431372549}

\usepackage{enumitem}
\setenumerate{label=\upshape(\alph*),itemsep=3pt,topsep=3pt}

\newtheorem{theorem}{Theorem}[section]

\newtheorem{myalgorithm}[theorem]{Algorithm}
\newtheorem{proposition}[theorem]{Proposition}

\newcommand{\Nedelec}{N\'ed\'elec\xspace}

\newcommand{\ci}{\mathrm{i}}

\newcommand{\kommentare}[1]{}

\newcommand*{\dint}[1]{\mathop{}\!\mathrm{d}#1}

\newcommand{\h}{h}
\newcommand{\hlarge}{h}

\newcommand{\VhLag}{V_h}
\newcommand{\VhLagDivC}{\widetilde{V}_h}

\newcommand{\VhNed}{\mathcal{R}_\h}

\renewcommand{\div}{\operatorname{div}}
\newcommand{\curl}{\operatorname{curl}}

\newcommand{\divh}{\operatorname{div_\hlarge}}

\newcommand{\sol}{u}

\newcommand{\MagF}{\mathbf{h}_{\mathrm{ext}}}
\newcommand{\MagFscalar}{\textup{h}_{\mathrm{ext}}}

\newcommand{\MagPot}{\mathbf{A}}

\newcommand{\R}{\mathbb{R}}

\newcommand{\C}{\mathbb{C}}
\newcommand{\N}{\mathbb{N}}

\newcommand{\VS}{H^1} 
\newcommand{\VSh}{V_{\h}}

\newcommand{\bfD}{\mathbf{D}}
\newcommand{\bfB}{\mathbf{B}}
\newcommand{\bfC}{\mathbf{C}}

\renewcommand{\L}{\mathbf{L}}
\renewcommand{\H}{\mathbf{H}}

\newcommand{\Hdiv}{\H(\div)}
\newcommand{\Hcurl}{\H(\curl)}
\newcommand{\Hcurlzero}{\H_0(\curl)}
\newcommand{\Hdivzero}{\H_0(\div)}
\newcommand{\Hcurldiv}{\H(\curl,\div)}
\newcommand{\Hcurldivzero}{\H_0(\curl,\div)}
\newcommand{\Hsol}{\H_0(\curl,\div^0)}
\newcommand{\Hsolh}{\VhNed^0}

\renewcommand{\Re}{\operatorname{Re}}

\usepackage{geometry}
\begin{document}

\title[GLENN: NN-enhanced computation of GL energy minimizers]
{GLENN: Neural network-enhanced computation of Ginzburg--Landau energy minimizers}

	\date{\today}

\author[M.~Crocoll]{Michael Crocoll}
\address{Institute for Applied and Numerical Mathematics, 
	Karlsruhe Institute of Technology, 76149 Karlsruhe, Germany}
\email{\{michael.crocoll,benjamin.doerich,roland.maier\}@kit.edu}
\author[C.~D\"oding]{Christian D\"oding}
\address{Institute for Numerical Simulation, University of Bonn, 53115 Bonn, Germany}
\email{doeding@ins.uni-bonn.de}
\author[B. D{\"o}rich]{Benjamin D{\"o}rich}

\author[R.~Maier]{Roland Maier}

\date{} 

\begin{abstract}
    In this work, we propose a neural network-enhanced finite element
    strategy to compute the minimizer of the Ginzburg--Landau energy based on an unsupervised deep Ritz-type strategy. We treat the parameter $\kappa$ as a variable input parameter to obtain possible minimizers for a large range of $\kappa$-values. This allows for two possible strategies: 1) The neural network may be extensively trained to work as a stand-alone solver. 2) Neural network results are used as starting values for a subsequent classical iterative minimization procedure. The latter strategy particularly circumvents the missing reliability of the neural network-based approach. Numerical examples are presented that show the potential of the proposed strategy.
\end{abstract}

\maketitle

\section{Introduction}

In 1911, the Dutch physicist Onnes \cite{Onn11} realized that when mercury is cooled below some critical temperature, the electrical resistance vanishes, i.e., the mercury is in a superconducting state. At the same time, one can observe that an applied magnetic field below some critical strength is expelled from the superconducting material. On the other hand, above a critical strength, the magnetic field penetrates and destroys the superconducting property. 
However, for so-called type-II superconductors, there is an intermediate regime in which the magnetic field only penetrates at isolated points while leaving the main part of the material in the superconducting state. Such a property is of particular interest in applications where 
the superconductor then becomes robust to weak perturbative magnetic fields.

The mathematical model to describe the phenomena of superconductivity on a macroscopic level is the Ginzburg--Landau (GL) model \cite[Sec.~3]{DuGP92}. Let $\Omega \subset \mathbb{R}^d$ with $d = 2$ or $d = 3$ denote the region occupied by the superconducting material. The superconducting properties are described by a complex-valued function $u : \Omega \to \mathbb{C}$, called the \textit{order parameter}. Although $u$ itself is not directly observable, similar to a quantum mechanical wave function, its squared modulus $|u|^2$ can be measured in physical experiments and interpreted as the density of superconducting electrons (\emph{Cooper pairs}). The value of $|u|^2$ ranges between $0$ and $1$ where $|u(x)|^2 = 1$ belongs to regions of superconductivity and $|u(x)|^2 = 0$ to regions where the material is in its normal (non-superconducting) state. Superconductivity is rare in nature but has a wide rage of possible applications.

Beyond the fully superconducting state ($|u|^2 \equiv 1$) and the normal state ($|u|^2 \equiv 0$), mixed states with $0 \le |u(x)|^2 \le 1$ arise, in particular, when an external magnetic field $\MagF$ is applied. The mixed normal-superconducting state can include a co-existence of superconducting and normal states -- a so-called \emph{Abrikosov vortex lattice}~\cite{Abr04}, where the density is zero in vortex centers. Such configurations are considered in this work, which, however, can only occur for special superconductors with appropriate magnetic fields. In particular, the internal magnetic field, represented by $\curl \MagPot$, tends to align with $\MagF$, where $\MagPot : \Omega \to \mathbb{R}^d$ is the  (typically unknown) magnetic vector potential within the superconductor. From a mathematical perspective, for the state of the superconductor the GL model prescribes that $(u, \MagPot)$ minimizes the \textit{GL free energy} (cf.~\cite[Sec.~3]{DuGP92}) given by
\begin{equation}\label{eq:energy_functional}
    E(u,\MagPot) \coloneqq \frac12 \int_{\Omega} 
	|\frac{\ci}{\kappa} \nabla u + \MagPot u |^2
	+
	\frac12 
	(1- |u|^2)^2 
	+
	|\curl \MagPot- \MagF|^2 \dint{x},
\end{equation}
where $\kappa > 0$ is a material parameter, also called the \emph{GL parameter}.

As indicated above, one is interested in vortex states, where $\kappa$ is linked to the structure of the vortex lattice, see~\cite{SaS07,SaS12,Ser99,SeS10}. There are no vortices for small values of $\kappa$, while larger values lead to more and more (very localized) vortices~\cite{Aftalion99,SaS07}. High values of $\kappa$ are particularly interesting and challenging, considering especially the numerical approximation of corresponding minimizers. 

The numerical approximation of minimizers of~\eqref{eq:energy_functional} has first been studied in~\cite{DuGP92,DuGP93}. Therein, error estimates in $H^1$ have been derived for both the order parameter $u$ and the magnetic vector potential $\MagPot$ when using a finite element discretization. The convergence results, however, do not take into account the role of $\kappa$, which can generally lead to constraints on the mesh size as investigated in~\cite{DoeH24} for a reduced setting with a fixed potential $\MagPot$. In that setting, the energy in~\eqref{eq:energy_functional} is only minimized with respect to $\sol$. 
The derived constraints on the choice of the mesh size may be overcome by the use of appropriate multiscale constructions as realized in~\cite{BDH24}. 
The full problem is considered in~\cite{DoeDH25*}, where different meshes and spaces for $\sol$ and $\MagPot$ are used as an adaptation to the different scaling behavior in $\kappa$. 
More recently,
in~\cite{ChaFH26} 
the 
precise coupling of not only the spatial mesh width~$h$ and the parameter~$\kappa$, but also the polynomial degree~$p$ is studied. The authors are able to precisely quantify the pre-asymptotic regime and the corresponding pollution effect,
and further show how higher-order polynomials reduce these effects.
An extension of this technique to multiscale approximation spaces is presented in~\cite{DoeH25*}.
A generalization to non-convex domains is provided in \cite{Doe25*}. Here, \Nedelec elements  and a discrete divergence-free condition have to be used to anticipate the low regularity of the vector potential. This is also the framework used in this work.
Besides, we mention here the works~\cite{DuNicolaidesWu98} and~\cite{QuJu05}, which conduct an error analysis for a covolume and a finite volume method, respectively (without explicit rates in the mesh size or $\kappa$).  

All of these works focus on the approximability of minimizers. For the practical computation in the numerical examples, different choices of minimization procedures are employed: 
in~\cite{BDH24} and~\cite{DoeH24}, a semi-implicit $L^2$-gradient flow is used, 
whereas in \cite{DoeDH25*} a Sobolev and in \cite{ChaFH26,DoeH25*} a conjugate Sobolev gradient flow approach are employed.
Even though it becomes clear that it appears to be always preferable to apply a conjugate Sobolev gradient flow for the GL problem, the results are not satisfying in the sense of reliably capturing the best minimizers without trial and error to find good starting values. 
An alternative approach is to use physics-based neural networks to approximate minimizers. In particular, the GL energy could be naturally used as a loss function in the spirit of a deep Ritz method~\cite{WeiY17} that is to be minimized by the training of a neural network. Nonetheless, a pure machine learning-based strategy can also be very sensitive to small changes in the network architecture.

In this work, we propose a neural network-enhanced strategy for the approximation of minimizers (of the discrete version) of~
\eqref{eq:energy_functional}. An essential issue of classical iterative algorithms is the choice of a starting value, which can heavily influence the number of required iterations and even the actual limit. 
To overcome this, we propose to train a neural network based on a certain parameter range of $\kappa$ to obtain a reasonable initial guess that leads to better GL minimizers. In essence, this also circumvents the problematic properties of neural network-based schemes, which to date are not fully reliable in the sense of a provable convergence behavior. More precisely, even if the suggested initial guess is not optimal, the subsequent classical iterative process still follows deterministic and theoretically grounded principles. Nonetheless, we also consider the setting where a network is trained to directly output a solution pair $(u,\MagPot)$ for a chosen value of $\kappa$. This, however, requires a certain amount of pre-computations to set up a valid model. The latter approach is particularly meant to show that alternatives to classical iterative algorithms based on finite element spaces have the potential to lead to even better minimizers. 

There exist similar approaches in the literature, where classical algorithms are enhanced by neural networks. Recent examples are in the context of nonlinear Schrödinger eigenvalue problems~\cite{PetPPR26} or to speed up iterative algorithms for convex optimization~\cite{VenA21} or inverse problems \cite{AdlO17}. 
However, to the best of our knowledge they all use supervised learning (i.e., training with data), while our approach is unsupervised, i.e., it does not need data and solely relies on the energy-based loss function for the training.

The rest of the paper is structured as follows.
In Section~\ref{sec:GLE_intro}, we recall the most important results on the minimizer of the GL energy 
and explain the established methods and how to derive our neural network-enhanced method from it. 
We present the neural networks approach in 
Section~\ref{sec:NN} 
and explain the finite element solver
in Section~\ref{sec:iterative_FEM}.
Finally, we provide numerical experiments in
Section~\ref{sec:num_exp}.

\section{Ginzburg-Landau problem and discrete minimizers}
\label{sec:GLE_intro}

\subsection{Preliminaries}

For the mathematical framework, we require the following notation: For $2 \le p \le \infty$ we denote by $L^p = L^p(\Omega,\C)$ the standard Lebesgue space of scalar complex-valued function on $\Omega$ equipped with the usual norm $\| \cdot \|_{L^p}$. For $p =2$ the Hilbert space $L^2$ is considered as a \textit{real} Hilbert space with the inner product
$   (v,w) = \Re \int_{\Omega} v \, \overline{w} \dint{x}$,
where $\overline{w}$ denotes the complex conjugate of $w$. Similarly, we denote by $\L^p = L^p(\Omega,\R^d)$ the space of $p$-Lebesgue-integrable vector fields, with norm $\| \cdot \|_{\L^p}$, and the inner product on $\L^2$ is given by
$    (\bfB,\bfC) = \int_{\Omega} \bfB \cdot \bfC \dint{x}$.
Note, that we use the same notation $(\cdot,\cdot)$ for both inner products on $L^2$ and $\L^2$ as its meaning is clear from the context. Accordingly, we denote by $H^1 = H^1(\Omega,\C)$ with norm $\| \cdot \|_{H^1}$ the Sobolev space of complex-valued functions whose weak derivative exists in $L^2$. The same notation applies to vector fields and we write $\H^1 = H^1(\Omega,\R^d)$. For $d = 3$, we denote by $\Hcurl$ the space of vector fields with weak rotation in $\L^2$ and by $\Hdiv$ the space of vector fields with weak divergence in $L^2$ . Both spaces can be equipped with natural zero boundary conditions leading to the subspaces
\begin{align}
    & \Hcurlzero = \{ \bfB \in \Hcurl \, | \, \bfB \times \nu = 0 \text{ on } \partial \Omega \}, \\
    & \Hdivzero = \{ \bfB \in \Hdiv \, | \, \bfB \cdot \nu = 0 \text{ on } \partial \Omega \},
\end{align}
where $\nu$ denotes the outer unit normal on $\partial \Omega$. Finally, we set
\begin{align}
    \Hcurldiv = \Hcurl \cap \Hdiv, \quad \Hcurldivzero = \Hcurl \cap \Hdivzero.
\end{align}

For $d = 2$, there are two $\curl$ operators that act either on scalar fields or on vector fields. For a scalar field $\MagF = \MagFscalar \, e_z$, it maps to the vector field $\curl \MagF = (-\partial_y \MagFscalar, \partial_x \MagFscalar)^\top$ whereas for a vector field $\MagPot$ the $\curl$ maps to the scalar field $\curl \MagPot = \partial_x \MagPot_2 - \partial_y \MagPot_1$.
For simplicity, we will use the same symbol, $\curl$, for both operators, as the context will make it clear which is meant. For any vector $a,b \in \R^2$ we formally set $a \times b = a_1 b_2 - a_2 b_1$. These definitions extend the spaces $\Hcurl$ and $\Hdiv$ and their corresponding subspaces of zero traces to two dimensions in a straightforward way. 

As mentioned above, the state of a superconductor is described by a minimizing pair $(u,\MagPot) \in H^1 \times \Hcurl$ of the GL free energy, i.e.,
\begin{align}
    E(u, \MagPot) = \min_{(v,\bfB) \in H^1 \times \Hcurl} E(v,\bfB).
\end{align}

The first Fr\'echet derivative of the energy functional
\begin{align}
    E'(u,\MagPot) = \big( \partial_u E(u,\MagPot), \partial_{\MagPot}E(u,\MagPot) \big): H^1 \times \Hcurl \rightarrow \R
\end{align}
 is given by
 \begin{subequations} \label{eq:E_prime}
\begin{align} 
    \langle \partial_u E(u,\MagPot), v \rangle & = \Big( \frac{\ci}{\kappa} \nabla u + \MagPot u, \frac{\ci}{\kappa} \nabla v + \MagPot v \Big) + \big( (|u|^2 - 1) u, v \big)
    \label{eq:E_prime_u}
    \\
    \langle \partial_{\MagPot} E(u,\MagPot), \bfB \rangle & = \Big(|u|^2 \MagPot + \frac{1}{\kappa} \Re (\ci \overline{u} \nabla u),\bfB \Big) + (\curl \MagPot - \MagF, \curl \bfB)
    \label{eq:E_prime_A}
\end{align}
\end{subequations}
for $(v,\bfB) \in H^1 \times \Hcurl$. By the first-order optimality condition each minimizer $(u,\MagPot)$, if it exists, needs to be a critical point of $E$, i.e., $\langle E'(u,\MagPot), (v,\bfB) \rangle = 0$ for all $(v,\bfB) \in H^1 \times \Hcurl$. This leads to the famous GL equations
\begin{align}
    \Big( \frac{\ci}{\kappa} \nabla u + \MagPot u, \frac{\ci}{\kappa} \nabla v + \MagPot v \Big) + \big( (1 - |u|^2) u, v \big) & = 0, \\
    \Big(|u|^2 \MagPot + \frac{1}{\kappa} \Re (\ci \overline{u} \nabla u),\bfB \Big) + (\curl \MagPot - \MagF, \curl \bfB) & = 0,
\end{align}
for all $(v,\bfB) \in H^1 \times \Hcurl$. One proves that indeed there exists at least one minimizing pair. Nevertheless, any such minimizer cannot be unique due to so-called \textit{gauge-transformations}.  Let $(u,\MagPot) \in H^1 \times \Hcurl$ be arbitrary. For any $\varphi \in H^2$ we define the gauge transformation
\begin{align}
    G_\varphi: H^1 \times \Hcurl \rightarrow H^1 \times \Hcurl, \quad G_\varphi(u,\MagPot) = (e^{\ci \kappa \varphi} u, \MagPot + \nabla \varphi).
\end{align}
It now holds that
$
E(u,\MagPot) = E \big( G_\varphi(u,\MagPot) \big),
$
i.e., the GL energy is invariant with respect to the gauge transformation, and hence a minimizer in $H^1 \times \Hcurl$ cannot be unique. Although we have to deal with the symmetry under complex phase shifts the order parameter $u$, we can remove the gauge invariance with respect to the vector potential $\MagPot$ by restricting the minimization to divergence-free vector potentials. For this so-called \emph{Coulomb gauge}, we employ the space
\begin{align}
    \Hsol =\{ \bfB \in \Hcurldivzero \, | \, \div \bfB = 0 \}.
\end{align}
Next, we recall a result,
which shows that at least one minimizer of the GL energy \eqref{eq:energy_functional} exists. Moreover, the minimizer can be chosen from the space $\Hsol$, making use of the gauge transformation, see 
\cite[Lemma~2.6]{Doe25*}
and \cite[Theorem~3.5]{DuGP92}.

\begin{proposition} \label{prop:existence_and_gauge} Let $\Omega \subset \R^d$, $d = 2,3$ be a bounded, simply connected Lipschitz domain and $\MagF \in L^2(\Omega,\R^{\ell_d})$. Then, the following statements hold.
\begin{enumerate}
    \item[(1)] \label{item3:existence} Let $(u,\MagPot) \in H^1 \times \Hcurl$ be arbitrary. Then, there exists a function $\varphi \in H^1$ such that $(\tilde u, \tilde\MagPot) \coloneqq G_\varphi(u,\MagPot) \in H^1 \times \Hsol$ and $E(\tilde{u}, \tilde{\MagPot}) = E(u,\MagPot)$.
    \item[(2)] \label{item2:existence} There exists at least one minimizer $(u,\MagPot) \in H^1 \times \Hsol$ of $E$.
    
\end{enumerate}
\end{proposition}

According to Proposition~\ref{prop:existence_and_gauge}, it suffices to seek for a minimizing pair $(u,\MagPot)$ with a divergence-free vector potential. We enforce the vanishing divergence of the vector potential by a weak side constraint leading to the following optimization problem: Find $(u,\MagPot) \in H^1 \times \Hcurl$ such that
\begin{equation} \label{eq:minimizer_constrained}
    E(u,\MagPot) = \min_{(v,\bfB) \in H^1 \times \Hcurl} E(v,\bfB) \quad \text{subject to } (\MagPot, \nabla v)_{\L^2} = 0 \text{ for all } v \in H^1.
\end{equation}
Therefore, every such minimizer automatically belongs to $H^1 \times \Hsol$, i.e., $\MagPot$ is divergence-free and satisfies the boundary condition $\MagPot \cdot \nu = 0$ on $\partial \Omega$.

Due to the constraint of vanishing divergence in $\MagPot$ the first-order conditions read as a saddle point problem involving a Lagrange multiplier to enforce the side constraint. A minimizer $(u,\MagPot) \in H^1 \times \Hcurl$ of \eqref{eq:minimizer_constrained} is also characterized by the solution triple~$(u,\MagPot,\lambda) \in H^1 \times \Hcurl \times H^1$ that solves the saddle point problem
\begin{equation} \label{eq:saddle_point}
\begin{split}
    \langle E'(u,\MagPot), (v,\bfB) \rangle + ( \nabla \lambda, \bfB )_{\L^2} & = 0 \quad \text{for all } (v,\bfB) \in H^1 \times \Hcurl, \\
     (\MagPot, \nabla w)_{\L^2} & = 0 \quad \text{for all } w \in H^1.
\end{split}
\end{equation}
Note that $\lambda \in H^1$ denotes the associated Lagrange multiplier. Conversely, every critical point solving~\eqref{eq:saddle_point}
corresponds to either a local minimum or a saddle point of the energy functional $E$. Loosely speaking, such a critical point with a low energy is a good candidate for the global minimizer if we verify the second-order optimality condition for a local minimum. That is, we have to exclude the possibility of a saddle point.

\subsection{A reduced Ginzburg--Landau model with given $\MagPot$}

To obtain first impressions on the applicability and potential of our approach, in this work we will also consider a reduced model of~\eqref{eq:energy_functional} where we assume $\MagPot$ to be given by the relation $\curl \MagPot = \MagF$. The corresponding modified energy then reads
\begin{equation}\label{eq:energy_functional_reduced}
    E(u) \coloneqq \frac12 \int_{\Omega} 
	|\frac{\ci}{\kappa} \nabla u + \MagPot u |^2
	+
	\frac12 
	(1- |u|^2)^2 
     \dint{x},
\end{equation}
and is only minimized for $\sol$. This has the advantage of fewer degrees of freedom while it still provides deep insight into many challenges.

\subsection{Numerical computation of Ginzburg--Landau minimizers}
\label{sec:num_comp_gle}

In the following section, we briefly discuss two approaches which are purely neural network-driven or finite element-driven, and focus mainly on their advantages and disadvantages. A thorough explanation is postponed to  Section~\ref{sec:NN} and to Section~\ref{sec:iterative_FEM} for the neural network and the 
finite element solver, respectively.
The main contribution of our work is to study suitable combinations of both approaches, which we will develop at the end of this section. The corresponding numerical experiments are presented in Section~\ref{sec:num_exp}. 

\subsubsection*{\textbf{Neural network approximation}}

A modern approach is the use of \emph{neural networks} (NNs). They are able to exploit modern architectures within \emph{graphics processing units} (GPUs) and their compute-heavy training is easily parallelizable. Several approaches for NNs with supervised and unsupervised learning have been developed over recent years ranging from data-driven models to physics-based NNs.
For most of these approaches, the strategy to train the NN is to minimize a prescribed loss functional. In a supervised setting, this could be a certain distance between the NN output and the labels of the training data.
However, our setting fits much better to an unsupervised learning strategy, since the GL energy \eqref{eq:energy_functional} gives rise to a natural loss function. This could be seen as a variant of the deep Ritz method. 
Another huge advantage of the usage of NNs is that there are many powerful optimization algorithms freely available in standard Python libraries such as \texttt{PyTorch}, which allow for the training of very large models.

However, finding a suitable architecture for a network as well as achieving a successful training process is rather challenging. More precisely, the choice of the precise network architecture (including many tunable parameters) itself is a high-dimensional optimization problem. Further, the resulting minimization problem in the training is generally non-convex. Thus, it is very difficult to reliably find a good network that is not just an unsatisfactory local minimum. 
Due to the large number of degrees of freedom, randomness at different stages of the training process is essential in order to obtain reasonable gradient descent methods. Examples are the initialization of the weights and the use of stochastic gradient descent methods. 
While randomness is computationally advantageous, it also means that from a theoretical perspective deterministic and reliable error estimates are typically not available. 

\subsubsection*{\textbf{Finite element approximations}}

A more classical approach to compute minimizers such as in~\eqref{eq:energy_functional} consists in the use of finite element (FE) methods. Here, one replaces the infinite-dimensional spaces by finite-dimensional approximation spaces, e.g., spaces of piecewise polynomials. Then one ends up with a finite-dimensional optimization problem, for which well-known methods can be used. One advantage is that the resulting matrices are in general sparse and thus pave the way for efficient computations in the minimization procedure.
In addition, it is possible to show rigorous approximation rates and resolution conditions 
in terms of the mesh size $h$ and the parameter $\kappa$ as mentioned in the introduction. 
Further, computing a-posteriori the smallest eigenvalues of $E''(u)$ can be used to avoid ending up in a saddle point.
In summary, this provides more reliability in the computed approximations compared to the pure NN approach. 

However, there are also some disadvantages of this approach.
Compared to the NN implementation, 
leveraging the computational potential of GPUs in a FE solver is not trivial,
and thus much harder to compete with the resulting computational speed-up,
even using an MPI (Message Passing Interface) implementation as done in our code.
In our experiments, we also observed a strong sensitivity to the chosen initial values 
leading to completely different discrete minimizers on different energy levels. 
Unfortunately, the high cost for a single minimization run prevents us from exploring the initial value landscape in a systematic manner.

\subsubsection*{\textbf{A neural network enhanced two-step approach (hybrid-approach)}}

Our idea is to combine the advantages of both above-mentioned methods in a so-called \emph{hybrid approach}. 
More precisely, we use the NN to compute a first approximation of the minimizer,
where we can exploit the great parallelization on GPUs.
We then interpolate the minimizer into the FE space and use it as an initial value for the iterative FE solver. Thus, one could interpret the NN outcome as a 
prediction
for the exact FE minimizer.
We thus have the advantage of replacing the heuristic choices of initial values with the more sophisticated surrogates of the NN while
computing the final minimizer with the reliability of the FE solver.
In addition, we observe later on that it is beneficial for the NN to compute the minimizer on a whole interval of admissible $\kappa$ values. 
The aim is then that a good choice of a small number of initial values excludes several local minimizers. 
Overall, we pursue the following two approaches:

\begin{itemize}[leftmargin=5mm]

    \item We perform a training of the NN with a problem-specific objective to compute approximations of GL energy minimizers. We train the network on a range of $\kappa$ and interpolate it for a fixed $\kappa$ into the FE space to compute energies that are comparable to the FE solver. Since the resulting network can be used on the full range of $\kappa$, one can then decide a-priori for how many values of $\kappa$ the network shall be used in order to budget the time for the training accordingly. The evaluation of the NN after the training to generate solutions for different $\kappa$ values is negligible regarding the runtime when compared to a FE solver.

    \item We use the same procedure as in the previous point, but use this interpolation as an initial value for the FE solver, to possibly achieve states of even lower energy. This enhances the FE solver with the NN and thus also indirectly allows the FE solver to profit from GPU computations without any adaptation on the FE side of the code.
    
\end{itemize}

\noindent
We elaborate on the precise strategies in our numerical examples in Section~\ref{sec:num_exp}.

\section{Neural networks}
\label{sec:NN}

In this section, we introduce the  NN approach that we use to solve the GL energy minimization problem~\eqref{eq:minimizer_constrained}. As mentioned before, this approach can be used as a stand-alone method to compute approximations of GL energy minimizers or as the first part of the hybrid-approach to compute an initial guess for the second level conjugate Sobolev gradient iteration. The approach is applied to the full problem \eqref{eq:energy_functional}, as well as the
reduced problem \eqref{eq:energy_functional_reduced}.\\
Before introducing the architecture and learning procedure of the NN for the GL problem, we first recall basic machine learning preliminaries for the sake of completeness. For further reading on the topic, we recommend~\cite{GooBC16}.

\subsection{Machine learning basics}
\subsubsection{MLPs, training, optimizers, and more}
One of the simplest NN architectures is realized by \emph{multi-layer perceptron (MLP)}. MLPs are defined as a function 
$   F=f^L \circ g^{L-1} \circ \ldots \circ g^0$,
where $g^{j} = \phi^j \circ f^j$. Here, $f^j$ are affine-linear maps and $\phi^j$ are non-linear functions that are applied entry-wise. The $\phi^j$ are called \emph{activation functions}. They raise the NN from being just an affine-linear map to a more expressive function. %
Note that NNs are not limited to said choice of the functions $g^j$. In fact, we use a different scheme in this work, see later. The individual maps $g^j$ are also referred to as the \emph{layers} of the NN, and the output dimension is typically called the \emph{width} of the layer. Often, several layers with a specific structure are grouped together to form \emph{blocks}. One very popular example of this are \emph{residual layers/blocks}. A residual block can, for example, consist of a few MLP layers, the output of which is called the \emph{residual branch} and is added to the input, hence the name residual. More formally, this means that instead of realizing 
$x \mapsto f_\mathrm{NN}(x),$
where $f_\mathrm{NN}$ is the output of the MLP layers, the NN block realizes
$\mathrm{resBlock}(x) = (I + f_{\mathrm{NN}})(x)$,
where $I$ is the identity map.

An MLP is parametrized by the matrices $W^j$ (\emph{weights}) and vectors $b^j$ (\emph{biases}) that define the affine-linear maps $f^j = W^jx + b^j$. A simple but typical machine learning task could now be to interpolate a set of points (\emph{samples}) and values (\emph{labels/targets}) with the NN. This set is then referred to as the \emph{training dataset}. In order to interpolate the training samples, a loss functional, e.g., a norm or metric in the space of the labels, is evaluated. This loss is then minimized over the space of network parameters in an iterative process called \emph{training}. To do this, one uses gradient descent like methods. The step size in these methods is referred to as the \emph{learning rate}. The learning rate is often not constant over the gradient descent steps, but is scheduled. One way this could look like is a gradual increase of the learning rate in the beginning (\emph{warm-up})~\cite{GoyDG18}, followed by a decay of the learning rate afterwards. Many of the \emph{optimizers}, i.e., the algorithms used for the minimization task in machine learning, are so-called \emph{momentum methods}. While a simple gradient descent gets stuck in a local minimum, momentum methods can overcome local minima. One very popular optimizer is \emph{Adam}~\cite{KingB14}, which stands for ``adaptive moment estimation". A variant of Adam is \emph{AdamW}. Here, the \emph{weight decay}, which is a penalty term that penalizes large model weights, is decoupled from the rest of the optimization. AdamW is a default choice in machine learning, which is why we utilize it in this work. Another optimizer is \emph{Muon}~\cite{KelJB24}, which has shown good results in the training of large language models~\cite{LiuSYJ25}. Due to said success, we use Muon in this work to evaluate whether it is a viable option for our setting as well. Muon is short for ``momentum orthogonalized by Newton-Schulz'' and essentially orthogonalizes the updates to the weight matrices, which can be beneficial. However, Muon can only be used on 2D parameters, which in our case only includes weight matrices. Other parameters still have to be trained with AdamW. For some of the networks trained, we used only AdamW, others were trained with Muon (with the adaptation from~\cite{LiuSYJ25}) and AdamW in combination. When using both optimizers, AdamW was used for the scaling parameter of the residual branch (see below for further details). Overall, we did not find a clear advantage of one optimizer over the other. Nonetheless, we include both in our experiments for the sake of completeness.

\paragraph{\textbf{Batches and epochs}}
Typically, the training dataset is processed in \emph{batches}, which are elements in a random partitioning of the dataset. They are fed into the NN one after another and on each of the batches a training step is performed. After all batches have been processed, i.e., the entire dataset has been seen, one says that an \emph{epoch} is completed. The batches are then re-drawn from the dataset and a new epoch begins. During training, one often computes \emph{validation} metrics to check how good the NN performs outside of the training dataset.

The set of parameters defining the architecture and training strategy of a NN but are not trainable are called \emph{hyperparameters} of the NN. These include, for example, the depth, width, and learning rate.

\paragraph{\textbf{Unsupervised learning}}
Previously, we mentioned a training dataset consisting of samples and labels or targets. This style of training an NN is called \emph{supervised learning}. But there is also \emph{unsupervised learning}, which refers to the situation where no labels or targets are known. Instead, there are only samples without corresponding labels. The quality of the approximation of the target values is simply encoded in the loss functional. One example for such unsupervised learning tasks are \emph{physics-informed neural networks (PINNs)}~\cite{RaiPK19}, where an NN is trained to solve a differential equation, by approximating the solution. The loss functional for a PINN is the mean squared error of the point-wise residual of the differential equation that, to some degree, measures how close the NN is to the solution without requiring the knowledge of said solution. The derivatives of the NN with respect to the input, that are required to compute the residual, can be computed using automatic differentiation. Another machine learning method to solve differential equations is the \emph{deep Ritz method}~\cite{WeiY17}. Instead of the strong formulation, a weak formulation is used here. This weak formulation is stated as a minimization problem, which then is a natural problem formulation for the NN training. For the application of machine learning to the GLE, we follow an unsupervised method similar to the deep Ritz method. 

\subsection{Machine learning for the Ginzburg--Landau problem}
\paragraph{\textbf{Loss functions}}
As indicated above, we use unsupervised learning for the GL problem. A major reason for this is that we do not have a dataset with target values. This means that we have to encode the metric of how good the NN performs in the loss function. For the GL problem, a canonical idea of a loss functional is the GL free energy. However, we also want to include the material parameter~$\kappa$ in the samples such that the trained model can be used for many runs of the solver. In order to do this, we integrate over a range of $\kappa$-values and also scale with the corresponding $\kappa$. In the discrete setting, this weighted integral reduces to a weighted sum.

The main motivation of this scaling is the asymptotic energy estimate in \cite[Theorem 8.1]{SanS07}, which (for constant $\MagF$) yields
$E(u,\MagPot,\kappa)\sim \frac{\log \kappa}{\kappa}$
as $\kappa\to\infty$. However, since we do not treat arbitrarily large~$\kappa$ and we deviate from the constant magnetic field, we simply neglect the logarithmic term, but keep the linear term to compensate vanishing energies for large $\kappa$ in the training. An additional benefit of this is that this gives a larger weight to samples corresponding to larger value of $\kappa$, for which the problems tends to be harder to solve.
The loss functional for the GL problem reads
\begin{equation}
{\mathcal{E}}(u,\MagPot)  \coloneqq \int_{\kappa_{\mathrm{min}}}^{\kappa_{\mathrm{max}}} \kappa \, \hat{E}(u,\MagPot,\kappa) \dint{\kappa},
\end{equation}
where
\begin{equation}
    	\hat{E}(u,\MagPot,\kappa)
        \coloneqq \frac12 \int_{\Omega} 
	|\frac{\ci}{\kappa} \nabla u + \MagPot u |^2
	+
	\frac12 
	(1- |u|^2)^2 
	+
	|\curl \MagPot- \MagF|^2 
    +|\mathrm{div} \MagPot |^2 \dint{x}.
\end{equation}

For the reduced GL problem \eqref{eq:energy_functional_reduced}, we simply plug in the considered $\MagPot$.
We also experimented with models trained for only a single $\kappa$ value, where we used $\hat{E}(u,\MagPot,\kappa)$ as the training loss. However, one of the findings of the preliminary experiments is that -- especially in the case of solving only for the order parameter -- training for a range of admissible $\kappa$ seems to yield better results for higher values of $\kappa$ than just training for the specific value. Particularly, increasing $\kappa_\textrm{max}$ beyond the actual target range to, e.g., $\kappa_\textrm{max}=130$, improved the results for $\kappa = 100$.

The procedure of learning the dependence of the solution on $\kappa$ can be interpreted for the GL problem as learning the parameter-to-solution map
\begin{equation}
    S\colon [\kappa_{\mathrm{min}},\kappa_{\mathrm{max}}] \rightarrow H^1(\Omega , \C \times \R^d), 
     \quad
    \kappa \mapsto (u_{\kappa}, \MagPot_{\kappa}) = \mathrm{argmin}(E(u,\MagPot,\kappa)).  
\end{equation}
A similar perspective can be chosen for the reduced problem.

\paragraph{\textbf{Neural network architecture}}
In contrast to the previously discussed simple architectures of NNs via MLPs, we use more sophisticated modern architectures, which turns out to be necessary and that we describe in the following. 
The network consists of a variable number of blocks with a variable width. Each block has a gated structure~\cite{DauFAG16}. Gating in this context means entry-wise multiplication of two branches, one of which is activated and called \emph{gate}. In particular, we make use of \emph{SwiGLU} layers~\cite{Sha20}, which have proven to be useful in the \emph{feed-forward} blocks of \emph{transformers}, which are network architectures that are used, e.g., in large language models. Feed-forward layers in this context can, for example, be MLP layers or, as just mentioned, SwiGLU layers. The main idea in the SwiGLU scheme is, that the input of the layer is split up into two linearly mapped branches. Then, one of the branches is activated and the result is multiplied entry-wise with the second branch. Afterwards, the result is once again mapped via a linear map, see \eqref{eq:swiglu}.
 We also propose a slight variation with a double activation scheme (which we refer to as \emph{DAGLU}). Here, an additional linear layer followed by an activation function is placed at the beginning of the block. In exchange, the final linear layer on the output is removed. This scheme can be combined with GELU or SiLU activations. Such structures have turned out to be advantageous in preliminary tests, which is why we included them in our constructions. 
The SwiGLU and DAGLU block formulas read
\begin{equation}
    \mathrm{SwiGLU}(x,W_1, W_2,V) = W_1(\mathrm{SiLU}(W_2x) \otimes Vx),
    \label{eq:swiglu}
\end{equation}
and
\begin{equation}
    \mathrm{DAGLU}(x,W_1,W_2,V,\alpha) = \alpha \big(W_2h\big) \otimes Vh, 
    \quad
    h  = \alpha(W_1x),
    \quad
    \alpha  \in \{ \mathrm{SiLU, GELU} \}.
    \label{eq:daglu}
\end{equation}
where $\otimes$ is the entry-wise multiplication of vectors and $W_1,W_2$ and $V$ are appropriately sized matrices.

The formulas for activation functions used above read 
\begin{equation}
    \mathrm{SiLU}(x) = \frac{x}{1 + \exp(-x)}\quad\text{and}\quad
    \mathrm{GELU}(x) = x\Phi(x),
\end{equation}
where $\Phi$ is the cumulative Gaussian distribution function.

The models we use have a residual structure in conjunction with the SwiGLU and DAGLU schemes to form blocks. The residual branch is scaled by a factor that can also be learned or fixed. 
We mostly use a learned entry-wise scaling vector~\cite{TouCSS21} with values in $(0,1]$, initialized by $1/(2L)$, where~$L$ is the number of blocks in the network. The small and layer-dependent initialization helps with stability in deeper networks without using normalization layers and is similar in concept to~\cite{ZhaDM19}. Using normalization layers in the NN proved detrimental for convergence in our preliminary testing.

The resulting formula for the SwiGLU block including the residual structure is
\begin{equation}
    \mathrm{RSwiGLU}(x,W_1, W_2,V,\gamma) = x + \gamma \otimes \mathrm{SwiGLU}(x,W_1, W_2,V),
    \label{eq:swiglu_block}
\end{equation}
 and $\gamma$ is a vector of scaling factors and, as before, $\otimes$ refers to the entry-wise multiplication of vectors. The DAGLU block including the residual structure reads
\begin{equation}     
\begin{aligned}
    \mathrm{RDAGLU}(x,W_1,W_2,V,\alpha,\gamma) &= x + \gamma \otimes \mathrm{DAGLU}(x,W_1,W_2,V,\mathrm{act}),
    \label{eq:daglu_block}
\end{aligned}
\end{equation}
where again, $\gamma$ is a vector of scaling factors and $\otimes$ is the entry-wise multiplication of vectors. 
For a visual representation of the SwiGLU and DAGLU block scheme, see Figure~\ref{fig:placeholder}.

\usetikzlibrary{positioning, shapes.geometric, arrows}

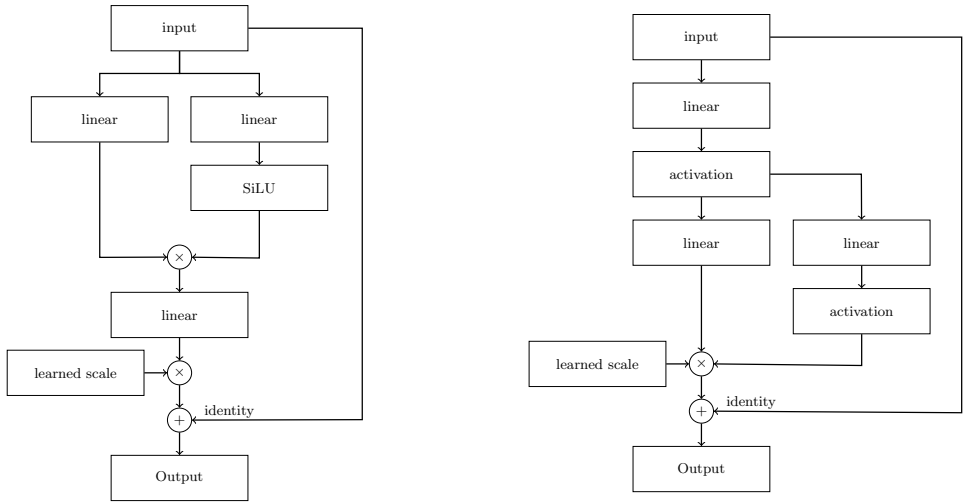
\begin{figure}
    \centering
    \begin{minipage}{0.45\textwidth}
        \scalebox{0.6}{\begingroup

\tikzset{
    block/.style={rectangle, draw, minimum width=3cm, minimum height=1cm, align=center},
    sum/.style={circle, draw, inner sep=2pt},
    arrow/.style={->, thick},
    every node/.style={font=\small}
}

\begin{tikzpicture}[node distance=1.5cm]

\node[block] (input) {input};
\node[block, below=1cm of input, xshift=-1.75cm] (Win1) {linear};
\node[block, below=1cm of input, xshift=1.75cm] (Win2) {linear};

\node[block, below=0.5cm of Win2](actgate) {SiLU};
\node[sum, below=0.75cm of actgate, xshift=-1.75cm] (mult) {$\times$};

\node[block, below=0.5cm of mult] (Wout) {linear};

\node[sum, below=0.5cm of Wout, xshift=-0cm] (mult2) {$\times$};

\node[block, left=0.5cm of mult2](scaling) {learned scale};

\node[sum, below=0.5cm of mult2] (sum) {$+$};

\node[block, below=0.5cm of sum] (output) {Output};

\draw[arrow] (input) -- ++(-0cm,-1.0cm) -- ++(-1.75cm,0) --(Win1);
\draw[arrow] (input) -- ++(-0cm,-1.0cm) -- ++(1.75cm,0) --(Win2);

\draw[arrow] (Win2) -- (actgate);
\draw[arrow] (Win1) -- ++(0,-3.05cm)--(mult);
\draw[arrow] (actgate) -- ++(0,-1.55cm)--(mult);
\draw[arrow] (scaling) -- (mult2);
\draw[arrow] (mult) -- (Wout);
\draw[arrow] (Wout) -- (mult2);
\draw[arrow] (mult2) -- (sum);

\draw[arrow] (sum) -- (output);
\draw[arrow] (input.east) -- ++(2.5cm,0) -- ++(0,-8.65cm)  -- (sum.east);

\node[right=0.15cm of sum, yshift=0.2cm] {identity};

\end{tikzpicture}

\endgroup}
    \end{minipage}
    \begin{minipage}{0.45\textwidth}
        \scalebox{0.6}{\begingroup

\tikzset{
    block/.style={rectangle, draw, minimum width=3cm, minimum height=1cm, align=center},
    sum/.style={circle, draw, inner sep=2pt},
    arrow/.style={->, thick},
    every node/.style={font=\small}
}

\begin{tikzpicture}[node distance=1.5cm]

\node[block] (input) {input};

\node[block, below=0.5cm of input] (mlp1) {linear};

\node[block, below=0.5cm of mlp1] (act) {activation};

\node[block, below=0.5cm of act] (mlp2) {linear};
\node[sum, below=1.9cm of mlp2] (attMult) {$\times$};
\node[block, left=0.5cm of attMult](scaling) {learned scale};
\node[sum, below=0.5cm of attMult] (sum) {$+$};
\node[block, right=0.5cm of mlp2](att) {linear};
\node[block, below=0.5cm of att](sigm) {activation};

\node[block, below=0.5cm of sum] (output) {Output};

\draw[arrow] (input) -- (mlp1);
\draw[arrow] (mlp1) -- (act);
\draw[arrow] (act) -- (mlp2);
\draw[arrow] (mlp2) -- (attMult);
\draw[arrow] (scaling) -- (attMult);
\draw[arrow] (attMult) -- (sum);
\draw[arrow] (sum) -- (output);
\draw[arrow] (act.east) -- ++(2cm,0) -- (att);
\draw[arrow] (att) -- (sigm);
\draw[arrow] (sigm) -- ++(0,-1.20cm) -- (attMult.east);
\draw[arrow] (input.east) -- ++(4.2cm,0) -- ++(0,-8.3cm)  -- (sum.east);

\node[right=0.15cm of sum, yshift=0.2cm] {identity};

\end{tikzpicture}

\endgroup}
    \end{minipage}
    \caption{Block schematic of a SwiGLU block (left) and a DAGLU block (right).}
    \label{fig:placeholder}
\end{figure}

\paragraph{\textbf{Hardware}}
The training of the NNs used in the experiments is done on the HoreKa Tier 2 High Performance Computing system at KIT. The training uses four nodes with four NVIDIA A100 GPUs per node for a total of 16 GPUs. For the tuning and preliminary experimentation, we use the HAICORE@KIT partition of HoreKa.
For the interpolation of the NNs, a GPU should be used. We leverage a workstation with a single NVIDIA A100 GPU for this purpose, which proved to be able to compute the forward passes reasonably fast.

\paragraph{\textbf{Learning rate}}
The scheduling of the learning rate during the training is divided into several phases. First, we do a short gradual warm-up~\cite{GoyDG18} of the learning rate to prevent instabilities. Here, the learning rate starts at a smaller value and is then increased over a set number of training steps. Second, a gradual decay of the learning rate is done using either an exponential decay or a cosine decay. The latter is a decay which follows a cosine in its shape and decays the learning rate down to a set minimum value. In the case of the cosine decay, we do an additional linear decay of the learning rate for a number of steps after the cosine scheduling ends. For the exponential scheduling, we skip this part in our experiments, i.e., the training was stopped before the linear decay phase. For the exact hyperparamters used for each of the models in this work, we refer to Table \ref{tab:parameters} and the configuration files in the source code that is available upon request.
There is no particular reason as to why which strategy was used for which model. There is no clear preference of one of the scheduling strategies over the other in our results and we include both strategies for completeness.

\paragraph{\textbf{Datasets}}
For the training of the networks including the material parameter $\kappa$, the dataset consists of random samples in 
$\Omega \times [\kappa_{\mathrm{min}},\kappa_{\mathrm{max}}] $. 
The samples are drawn in bulk every epoch. These samples are then processed in batches until each sample has been used once, then the samples are re-drawn with a new seed.

\section{Finite element approximation and gradient descent method}
\label{sec:iterative_FEM}

In order to obtain a FE formulation of the GL problem, we simply restrict the minimization in \eqref{eq:minimizer_constrained} to some finite-dimensional subspace. 
For our purposes, we choose continuous finite elements of order $p \in \N$ for the real and imaginary part of the order parameter~$u$,
denoted by $\VhLag$,
and \emph{\Nedelec elements} of first kind and order $m$, denoted by $\VhNed$.
Nedelec elements are the natural choice for our approach, since they only enforce continuity in the tangential direction on each element. Further, since they are not a subspace of $\H^1$, they can still capture the correct solution behavior even in  the case of non-convex domains. 
In order to enforce the divergence constraint, we have to define the weak divergence-free condition. To this end let 
$\VhLagDivC$ denote the space of continuous elements of order $m$, and define the discrete divergence as
\begin{equation} \label{eq:discrete_div}
    \divh \bfB_h = 0 \iff (\bfB_h , \nabla \varphi_h) = 0 
    \text{ for all } \varphi_h \in \VhLagDivC.
\end{equation}
Let us note that this does not only enforce that the weak divergence vanishes, but in addition also that the normal trace is zero in a weak sense.
Then, the discrete minimization problem reads: seek $(u_h,\MagPot_h) \in \VhLag \times \VhNed$ such that 
\begin{equation} \label{eq:minimizer_constrained_h}
    E(u_h,\MagPot_h) = \min_{(v_h,\bfB_h) \in \VhLag \times \VhNed} E(v_h,\bfB_h) \quad \text{subject to } \divh \MagPot_h = 0.
\end{equation}
By changing the spaces, the minimizer of 
\eqref{eq:minimizer_constrained_h}
satisfies the discrete analogue of \eqref{eq:saddle_point}.
Further, we define 
$
\Hsolh = \{ \bfB_h \in \VhNed \mid  \divh \bfB_h = 0  \}
$
as
the discrete counterpart of $\Hsol$
such that \eqref{eq:minimizer_constrained_h} can be formulated as
\begin{equation}
    E(u_h,\MagPot_h) = \min_{(v_h,\bfB_h) \in \VhLag \times \Hsolh} E(v_h,\bfB_h).
\end{equation}
Let us emphasize that this FE approach is non conforming since a discrete divergence free function in \eqref{eq:discrete_div} is in general not divergence free,
i.e.,
$\Hsolh \not\subseteq \Hsol$,
but it only holds the weaker inclusion 
$\Hsolh \subset \Hcurl$.

Next, we propose a special iterative solver to compute a reliable approximation of a discrete GL energy minimizer. 
Since we only apply it to the discretized energy \eqref{eq:minimizer_constrained_h}, we derive the algorithm only for the discrete case, even though this can be done in greater generality. To keep the notation simple, we also refrain from adding lower indices $h$ to every variable. We emphasize that this method was recently proposed in \cite{ChaFH26} for the case reduced case \eqref{eq:energy_functional_reduced}, and our method is a generalization of this approach to the full problem \eqref{eq:energy_functional}.

Starting from 
an initial guess $(u^0, \MagPot^0)$ in the FE space (which for example can be obtained by some NN),
the iteration takes the form of a descent method given by
\begin{align}
    (u^{k+1}, \MagPot^{k+1}) = (u^{k}, \MagPot^{k}) + \tau_k (d^k, \bfD^k), \quad k = 0,1,2,\dots,
\end{align}
where $(d^k, \bfD^k) \in \VhLag \times \Hsolh \subset \VS \times \Hcurl$
is the descent direction that we specify in the following. For this purpose, 
we equip the tensor space $X^k \coloneqq \VhLag \times \Hsolh$ with the $k$-dependent metric
\begin{align}
    \big( (v,\bfB), (w,\bfC) \big)_{X^k} := a_k(v,w) + b_k(\bfB,\bfC), \quad (v,\bfB),(w,\bfC) \in X^k,
\end{align}
where
\begin{align}
    a_k(v,w) & := ( \tfrac{\ci}{\kappa} \nabla v + \MagPot^k v,  \tfrac{\ci}{\kappa} \nabla w + \MagPot^k w)_{L^2} + \big( (\beta + |u^k|^2 + |\MagPot^k|^2 ) v,w \big)_{L^2} \\
    b_k(\bfB,\bfC) & := ( \curl \bfB, \curl \bfC )_{\L^2} + \big( (\beta + |u^k|^2)\bfB, \bfC \big)_{\L^2}
\end{align}
for some stabilization parameter $\beta > 0$. In particular, $(X^k)_{k \in \N}$ forms a sequence of Hilbert spaces defined by the same underlying vector space, but with a varying metric in each iteration which is indicated by the super index.
Note that if $\beta$ is sufficiently large (independent of $k$), then $a_k(\cdot,\cdot)$ and $b_k(\cdot,\cdot)$ are coercive and bounded on $\VhLag$ and $\Hsolh$, respectively. Therefore, $(\cdot,\cdot)_{X^k}$ is indeed is an inner product on the vector space $X^k$. However, the subsequent derivation of the conjugate Sobolev gradient method also applies for other choices as long as $(X^k, (\cdot,\cdot)_{X^k} )$ is a Hilbert space and $X^k \subset H^1 \times \Hcurl$.

For each $X^k$, we can formally define the Sobolev gradient $\nabla_{X^k} E(u^k,\MagPot^k)$ of the GL energy functional at the current iterate $(u^k,\MagPot^k) \in X^k$ by
\begin{equation} \label{eq:SobolevGradient}
    \big( \nabla_{X^k} E(u^k,\MagPot^k), (v,\bfB) \big)_{X^k} = \langle E'(u^k,\MagPot^k), (v,\bfB) \rangle \quad\text{for all } v \in \VhLag,\; \bfB \in \Hsolh,
\end{equation}
where $E'(u^k,\MagPot^k) \in (\VhLag \times \Hsolh)^*$ 
denotes the Frech\'et derivative and $\langle \cdot, \cdot \rangle$ is the dual pairing between $H^1 \times \Hcurl$ and its dual. 
Roughly speaking, the metric $(\cdot,\cdot)_{X^k}$ is chosen in a way such that the Sobolev gradient $\nabla_{X^{k}} E(u,\MagPot)$ is close to the identity map.
Recalling the Frech\'et derivatives in the directions of $u$ and $\MagPot$ given in \eqref{eq:E_prime_u} and \eqref{eq:E_prime_A}, respectively,
a straightforward calculation shows that the Sobolev gradient \eqref{eq:SobolevGradient} is given by
\begin{align}
    \nabla_{X^k} E(u^k,\MagPot^k) = (u^k, \MagPot^k) - \big(\delta^k,  \Delta^k \big)
\end{align}
where $\delta^k \in \VhLag$ is the solution of the elliptic problem
\begin{align}
    a_k(\delta^k, v) = \big( (1 + \beta + |\MagPot^k|^2) u^k, v \big)_{L^2} \quad \text{for all } v \in \VhLag
\end{align}
and $\Delta^k \in \Hsolh$ is the solution of the elliptic problem
\begin{align}
     b_k(\Delta^k, \bfB) = (\beta \MagPot^k - \tfrac{1}{\kappa} \Re (\ci \overline{u}^k \nabla u^k), \bfB)_{\L^2} + (\MagF, \curl \bfB)_{\L^2} \quad \text{for all } \bfB \in \Hsolh.
\end{align}
The descent direction is now defined in the spirit of a nonlinear conjugate Sobolev gradient descent method by
\begin{align}
    (d^k,\bfD^k) := \begin{cases} - \nabla_{X^k} E(u^k,\MagPot^k), & k =0 \\ - \nabla_{X^k} E(u^k,\MagPot^k) + \gamma_k (d^{k-1},\bfD^{k-1}), & k \ge 1
    \end{cases}
\end{align}
where $\gamma_k$ is the Polak-Ribi\`ere dissipation parameter, cf. \cite{PolR69}, defined by
     \begin{align}
         \gamma_k = \max \Big\{ 0, \frac{\big( \nabla_{X^k} E(u^k,\MagPot^k), \nabla_{X^k} E(u^k,\MagPot^k) - \nabla_{X^{k-1}} E(u^{k-1},\MagPot^{k-1}) \big)_{X^k}}{\big(   \nabla_{X^{k-1}} E(u^{k-1},\MagPot^{k-1}),\nabla_{X^{k-1}} E(u^{k-1},\MagPot^{k-1}) \big)_{X^{k-1}}} \Big\}.
     \end{align}
We further orthogonalize the search direction for $\MagPot$ in the sense of \eqref{eq:discrete_div} such that $\divh \bfD^k = 0$, which ensures the discrete divergence condition for the next approximation $\MagPot^k$.
Along the descent direction, we can then look for the optimal step size $\tau_k$ that minimizes the energy
\begin{align}
    \tau_k = \underset{\tau > 0}{\mathrm{argmin}} \, E\big( (u^{k}, \MagPot^{k}) + \tau (d^k,\bfD^k) \big).
\end{align}
Note that $E\big( (u^{k}, \MagPot^{k}) + \tau (d^k,\bfD^k) \big)$ is a fourth order polynomial in $\tau$, which we denote by $p_k$. Then, $\tau_k$ can be computed efficiently by a line search or we compute the roots of $p_k'$ directly.

We can now formulate the final algorithm for the nonlinear conjugate Sobolev gradient descent method. 

\begin{myalgorithm}[nonlinear conjugate Sobolev gradient descent method] Let initial values $u^0 \in \VSh$ and $\MagPot^0 \in \Hsolh$ with $u^0 \neq 0$ and a tolerance $\mathrm{tol} > 0$ be given. Further, set $d^0 =0$ and $\bfD^0 = 0$. Then for $k = 0,1,2,\dots$ perform the following steps:
\begin{enumerate} [leftmargin=7mm]
    \item[1.] \textbf{Solve} for $\delta^k \in \VSh$ such that
    \begin{align}
         a_k(\delta^k, v) = \big( (1 + \beta + |\MagPot^k|^2) u^k, v \big)_{L^2} \quad \text{for all } v \in \VSh
    \end{align}
    and for $\Delta^k \in \VhNed$ such that
    \begin{align}
     b_k(\Delta^k, \bfB) = (\beta \MagPot^k - \tfrac{1}{\kappa} \Re (\ci \overline{u}^k \nabla u^k), \bfB)_{\L^2} + (\MagF, \curl \bfB)_{\L^2} \quad \text{for all } \bfB \in \VhNed.
     \end{align}
     \item[2.] \textbf{Update} the descent directions
     \begin{align}
         d^k = \delta^k - u^k + \gamma_k d^{k-1}, \quad \bfD^k = \Delta^k - \MagPot^k + \gamma_k \bfD^{k-1}
     \end{align}
     with the Polak-Ribi\`ere parameter $\gamma_k = \max \{ 0, \frac{N_k}{B_k} \}$ with
     \begin{align}
    \hspace{-3mm} B_k &= \big(   (\delta^k - u^k, \Delta^k - \MagPot^k),(\delta^k - u^k, \Delta^k - \MagPot^k) \big)_{X^{k-1}}
    \\
    \hspace{-3mm} N_k &= \big( (\delta^k - u^k, \Delta^k - \MagPot^k), (\delta^k - u^k - \delta^{k-1} + u^{k-1}, \Delta^k - \MagPot^k - \Delta^{k-1} + \MagPot^{k-1}) \big)_{X^k}
   \end{align}
and orthogonalize such that $\divh \bfD^k = 0$.

     \item[3.] \textbf{Determine} the step size $\tau_k$ such that
     \begin{align}
        \tau_k = \underset{\tau > 0}{\mathrm{argmin}} \, E\big( (u^{k}, \MagPot^{k}) + \tau (d^k,\bfD^k) \big)
    \end{align}
    via a line search (e.g., golden search) or direct computation of the roots of $p_k'$.
    
    \item[4.] \textbf{Iterate} according to
         $(u^{k+1}, \MagPot^{k+1}) = (u^{k}, \MagPot^{k}) + \tau_k (d^k,\bfD^k)$.

    \item[5.] \textbf{Terminate} if
      $E(u^{k}, \bfD^k) - E(u^{k+1}, \bfD^{k+1}) < \mathrm{tol}$.

\end{enumerate}
\end{myalgorithm}

In practice, it is sufficient to set $\beta = 0$. The reason is that the iterates $(u^k, \MagPot^k)$ typically contribute enough mass through the second terms in the bilinear forms $a_k(\cdot,\cdot)$ and $b_k(\cdot,\cdot)$ to ensure coercivity on $\VhLag$ and $\Hsolh$, respectively.

\section{Numerical results}
\label{sec:num_exp}

In this section, we present the numerical results obtained using the NN approach implemented in \texttt{PyTorch} \cite{PyTorch_2}, as well as the hybrid approach that combines the NN with the FE iterative solver implemented in \texttt{DOLFINx} \cite{DOLFINx23}. We compare these results to those computed by the FE iterative solver using empirically chosen initial values to emphasize the advantages of the hybrid approach. We investigate various GL settings and training strategies for NNs. First, we consider NNs that are intensively trained for the reduced and full GL models. These networks lead to minimizers on lower energy levels than the ones obtained with the classical FE iterative approach with heuristic initial values. These models also provide a set of architectures and training parameters used in subsequent experiments involving faster-trained networks for GL models defined on different geometric domains, i.e., the unit square and the L-shaped domain. Table \ref{tab:parameters} depicts the most relevant parameters for the architecture and training strategy of the NN models used in all the experiments. All our models consist of 40 blocks of varying type with a width of 512 neurons resulting in about 31.5M trainable parameters. The main differences between the models lie in the chosen block type as well as the optimizer used for training. Regarding the learning rates, \texttt{target lr main} is the learning rate used for the weights, while \texttt{target lr scaling} is the learning rate used for the scaling vector for the residual branch. Further details on the hyperparameters can be found in the corresponding configuration files. These files are provided with the source code for the entire implementation, available at:
\begin{center}
    \url{https://github.com/BenjaminDoerich/glenn}
\end{center}

\begin{table}[htbp]
\centering
\begin{tabular}{l c c c c c}
\toprule
 & \texttt{GL-R1} & \texttt{GL-R2} & \texttt{GL-F1} & \texttt{GL-F2} & \texttt{GL-L*} \\
\midrule
\texttt{depth (blocks)}      & 40   & 40   & 40   & 40   & 40 \\
\texttt{width}               & 512  & 512  & 512  & 512  & 512 \\
\texttt{trainable param.} & 31.5M & 31.5M & 31.5M & 31.5M & 31.5M \\
\midrule
\texttt{block type}          & SwiGLU & SwiGLU & SwiGLU & DAGLU & SwiGLU \\
\texttt{activation}          & SiLU   & SiLU   & SiLU & SiLU & SiLU \\
\texttt{precision}           & 32 & 64 & 32 & 64 & 32 \\
\midrule
\texttt{batch size}          & $2^{14}$ & $2^{15}$ & $2^{15}$ & $2^{15}$ & $2^{15}$ \\
\texttt{steps per epoch}     & $2^{13}$ & $2^{13}$ & $2^{13}$ & $2^{13}$ & $2^{13}$ \\
\texttt{total epochs}        & 20 & 18  & 19 & 6 & 4 \\
\midrule
\texttt{target lr main}      & 0.001  & 0.003  & 0.003 & 0.003 & 0.01 \\
\texttt{target lr scaling}   & 0.001  & 0.0003 & 0.0003 & 0.0003 & 0.001 \\
\texttt{decay scheduler}     & cosine & exp. & exp. & exp. & exp. \\
\texttt{optimizer}           & AdamW  & AdamW  & Muon & Muon & Muon \\
\midrule
$\kappa_{\mathrm{max}}$      & 110  & 130  & 130 & 130 & 130 \\
\bottomrule
\end{tabular}
\caption{Architecture and training parameters for the NN models.}
\label{tab:parameters}
\end{table}

For the hybrid approach and for computing reference approximations, we exploit the FE iterative solver with quadratic Lagrange FE in the order parameter $u$ and \Nedelec elements of first kind and lowest order for the vector potential $\MagPot$. 
On the unit square domain, a quasi-uniform mesh is used with 608,842 real (or 304,421 complex) degrees of freedom 
for the order parameter and 228,059 degrees of freedom for the vector potential, which amounts to a total of 836,901 real degrees of freedom.
For the L-shaped domain we have 457,306 real (or 228,653 complex) degrees of freedom 
for the order parameter and 170,233 degrees of freedom for the vector potential.  This amounts to a total of 628,539 real degrees of freedom. Depending on the precise setting, one FE step takes around 5--7 seconds, leading to 80--120 minutes per 1000 iterations.
The iterative solver for the classical and hybrid approach terminates once a tolerance of $\mathrm{tol} = 10^{-12}$ is reached with respect to the energy increments.

As mentioned above, to quantify the quality of the NN and hybrid approach we compute reference approximations using the classical FE iterative solver with empirically chosen initial values. These are denoted by~$\varphi_j$, $j = 1,2,3,4,5$ and generated for $\alpha = \tfrac{1}{\sqrt{2}}(1 + \ci)$ through the functions
\begin{gather}
\psi_1(x) = \alpha e^{-|x|^2}, \qquad
\psi_2(x) = \tfrac{1}{\sqrt{\pi}} \big( \tfrac{2}{3} (x_1 + \ci x_2) + \tfrac{1}{2} \big) e^{-|x|^2},\qquad
\psi_3(x) = \alpha e^{10 \ci |x|^2},
\\
\psi_4(x) = \alpha (x_1 + \ci x_2) e^{10 \ci |x|^2}, \qquad
\psi_5(x) = \alpha. 
\end{gather}
With the transformation $\chi: (0,1)^2 \rightarrow (-1,1)^2$, $\chi(x) = (2x_1 - 1, 2x_2 - 1)$, we set the initial values for the order parameter on the unit square to interpolation into the FE space of
$
    u^0 = \varphi_j, \, \varphi_j = \psi_j \circ \chi
$
for $j = 1,2,3,4,5$. For the L-shaped domain we restrict these functions accordingly. In each GL setting, this leads to five reference approximations from which the one with the lowest energy level is the best candidate for the global minimizer obtained with the classical FE iterative solver.

\subsection{Neural networks for the reduced GL model} \label{sec:experiment_reduced_model}

For the first experimental setup, we consider the reduced model \eqref{eq:energy_functional_reduced} on the unit square $\Omega = (0,1)^2$. The magnetic vector potential $\MagPot$ is given by
\begin{equation}
    \MagPot(x) = \sqrt{2} \begin{pmatrix} 
    \sin( \pi x_1) \cos(\pi x_2) \\ - \cos(\pi x_1) \sin(\pi x_2)
    \end{pmatrix},  \quad x = (x_1,x_2) \in \Omega.
\end{equation}
Note that $\div \MagPot = 0$ and that $\curl A = \MagF$ with $\MagF(x) = 2\sqrt{2} \pi \sin(\pi x_1) \sin(\pi x_2 )$. The GL parameter is set to the values $\kappa = 10, 25, 50, 75, 100$. We emphasize again that the training of the NN does not particularly include these values, but rather provides approximations for a whole range of values of $\kappa \in [\kappa_{\mathrm{min}}, \kappa_{\mathrm{max}}]$.

\begin{table}[h!]
\centering
\small
\begin{tabular}{c c c c c c}
\toprule
$\kappa$ & 10 & 25 & 50 & 75 & 100 \\
\midrule
$\varphi_1$ & \textbf{0.10459064} & 0.08611814 & 0.06928056 & 0.03607700 & 0.04583871 \\
$\varphi_2$ & 0.11497944 & 0.07431313 & 0.06231434 & 0.04038106 & 0.04301528 \\
$\varphi_3$ & 0.10459064 & \textbf{0.06483810} & \textbf{0.04141371} & 0.03611419 & 0.03091586 \\
$\varphi_4$ & 0.10459064 & 0.06483810 & 0.04141371 & \textbf{0.03448296} & \textbf{0.03091576} \\
$\varphi_5$ & 0.10459064 & 0.08386487 & 0.06876863 & 0.04039311 & 0.04045206 \\
\midrule
\makebox[18mm][l]{\texttt{GL-R1} (NN)}
& 0.10462698 & 0.06523432 & 0.04215652 & 0.03212462 & 0.02785582 \\
\makebox[18mm][l]{\texttt{GL-R2} (NN)}
& 0.10460000 & 0.06525657 & 0.04222260 & 0.03308495 & 0.02609508 \\
\midrule
\makebox[18mm][l]{\texttt{GL-R1} (hyb.) }
& 0.10459064 & 0.06450104 & 0.04141371 & 0.03162426 & 0.02770266 \\
\makebox[18mm][l]{\texttt{GL-R2} (hyb.)}
& 0.10459064 & 0.06450091 & 0.04141371 & 0.03169532 & 0.02605176 \\
\bottomrule
\end{tabular}
\caption{Energy levels $E(u)$ of computed energy minimizers for the reduced GL model of Section \ref{sec:experiment_reduced_model}. First block: FE iterative solver with heuristic initial values~$\varphi_j$. Second block: pure NN models
\texttt{GL-R1} 
(NN) and
\texttt{GL-R2}
(NN). Third block: hybrid approach
\texttt{GL-R1}
(hyb.) and
\texttt{GL-R2}
(hyb.) with initial values from the NN. }
\label{tab:energy_only_ord}
\end{table}

\begin{figure}[h!]
    \flushleft
    \begin{minipage}{0.1\textwidth}
    \centering
    $\varphi_j$ (best)
    \end{minipage}
    \begin{minipage}{0.15\textwidth}
    \includegraphics[scale=0.07]{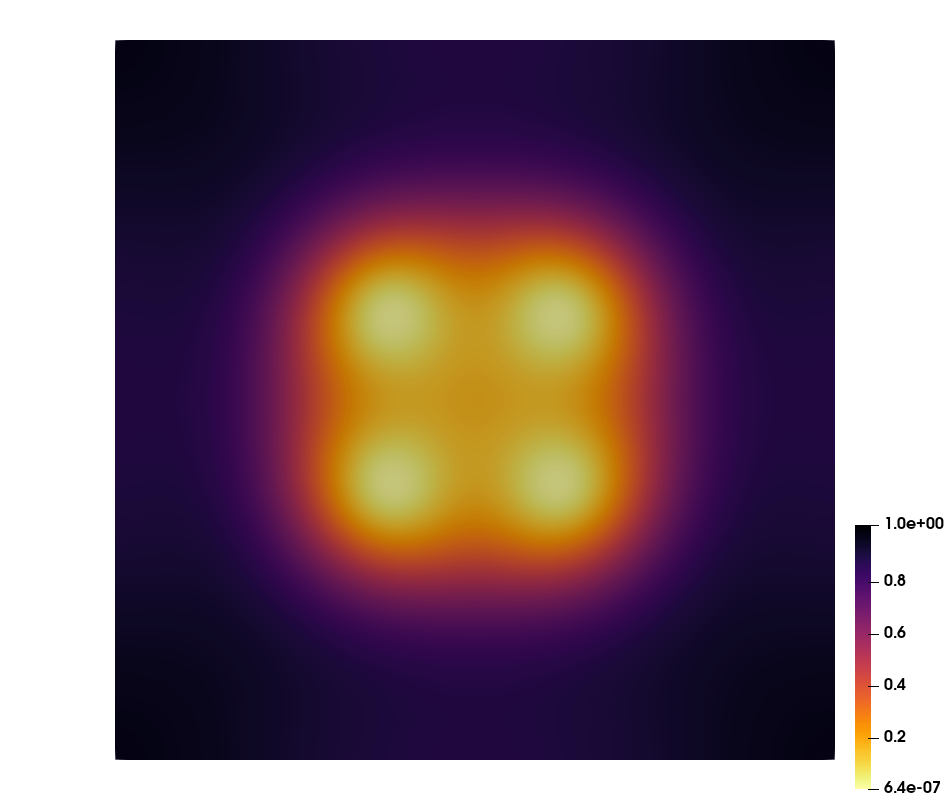}
    \end{minipage}
    \begin{minipage}{0.15\textwidth}
    \includegraphics[scale=0.07]{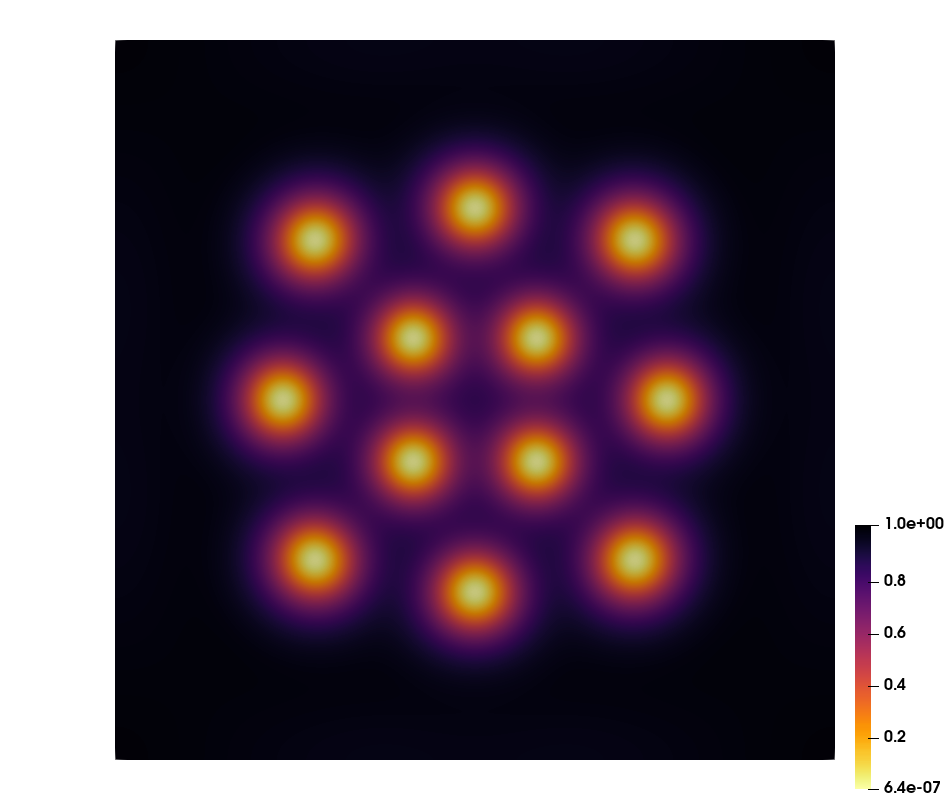}
    \end{minipage}
    \begin{minipage}{0.15\textwidth}
    \includegraphics[scale=0.07]{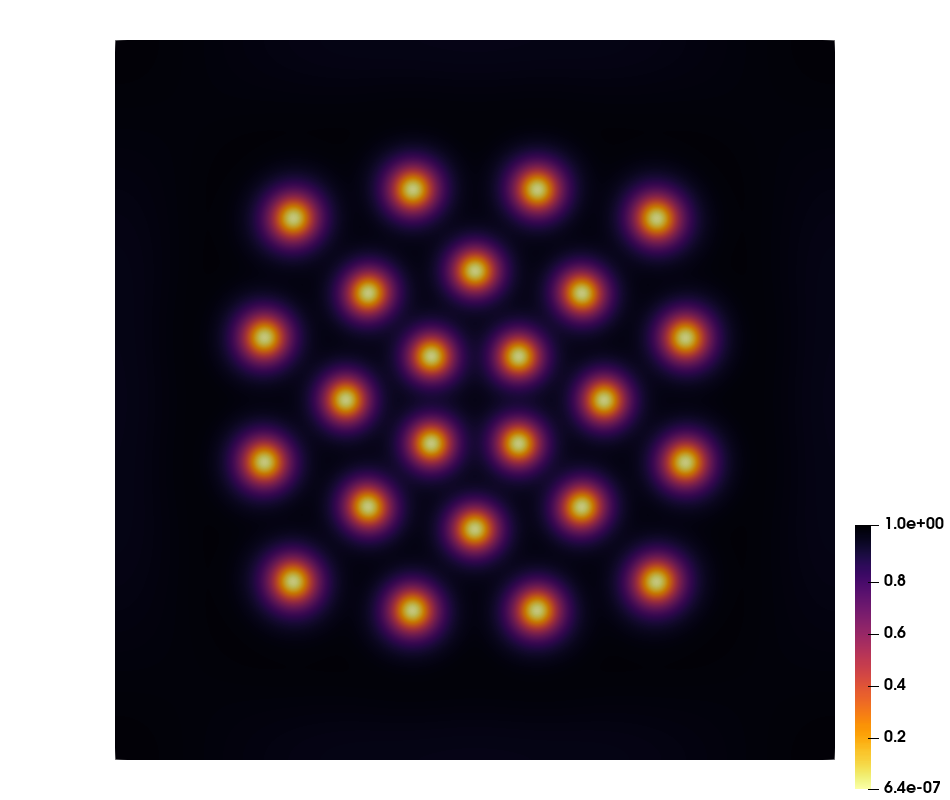}
    \end{minipage}
    \begin{minipage}{0.15\textwidth}
    \includegraphics[scale=0.07]{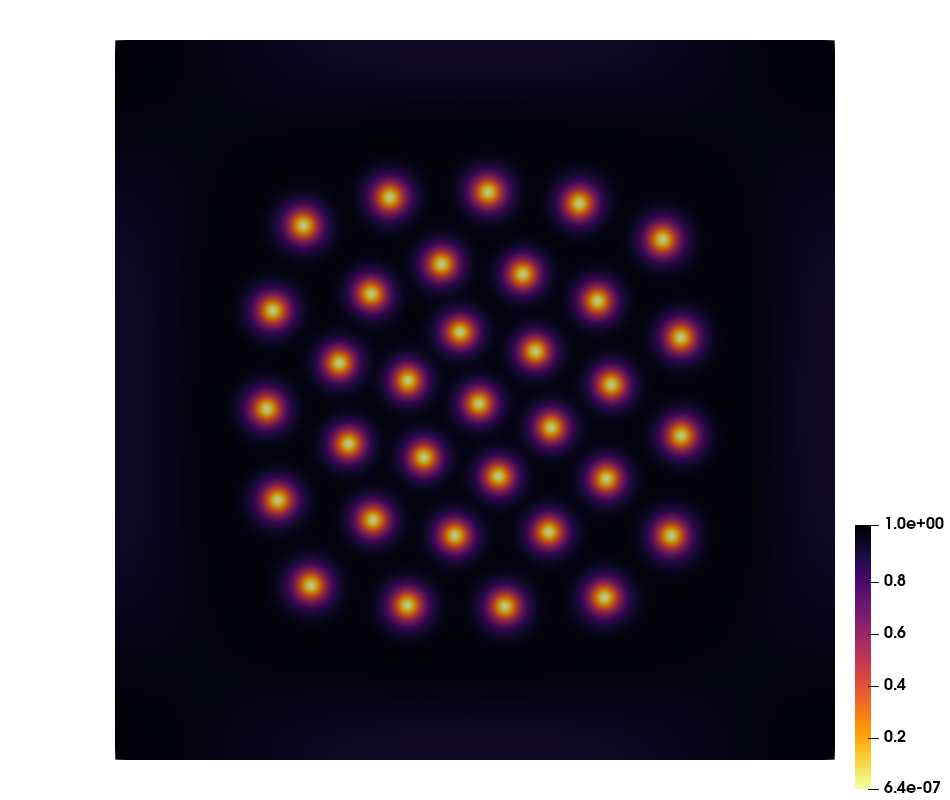}
    \end{minipage}
    \begin{minipage}{0.15\textwidth}
    \includegraphics[scale=0.07]{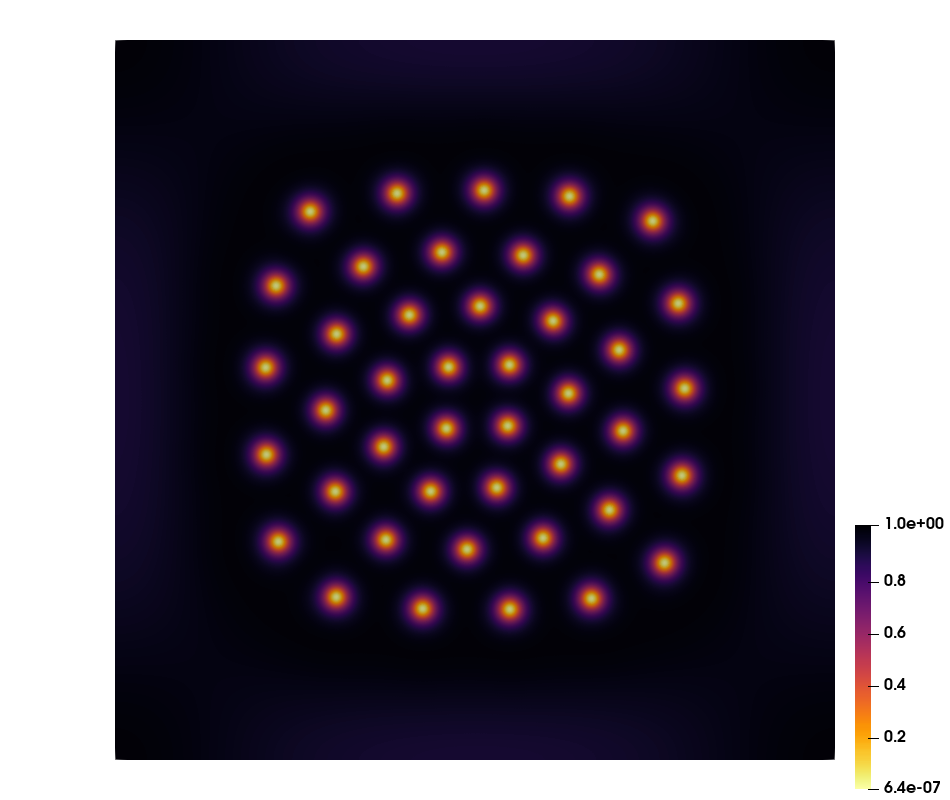}
    \end{minipage} 
    \\
    \begin{minipage}{0.1\textwidth}
    \centering
    \texttt{GL-R1} (NN)
    \end{minipage}
    \begin{minipage}{0.15\textwidth}
    \includegraphics[scale=0.07]{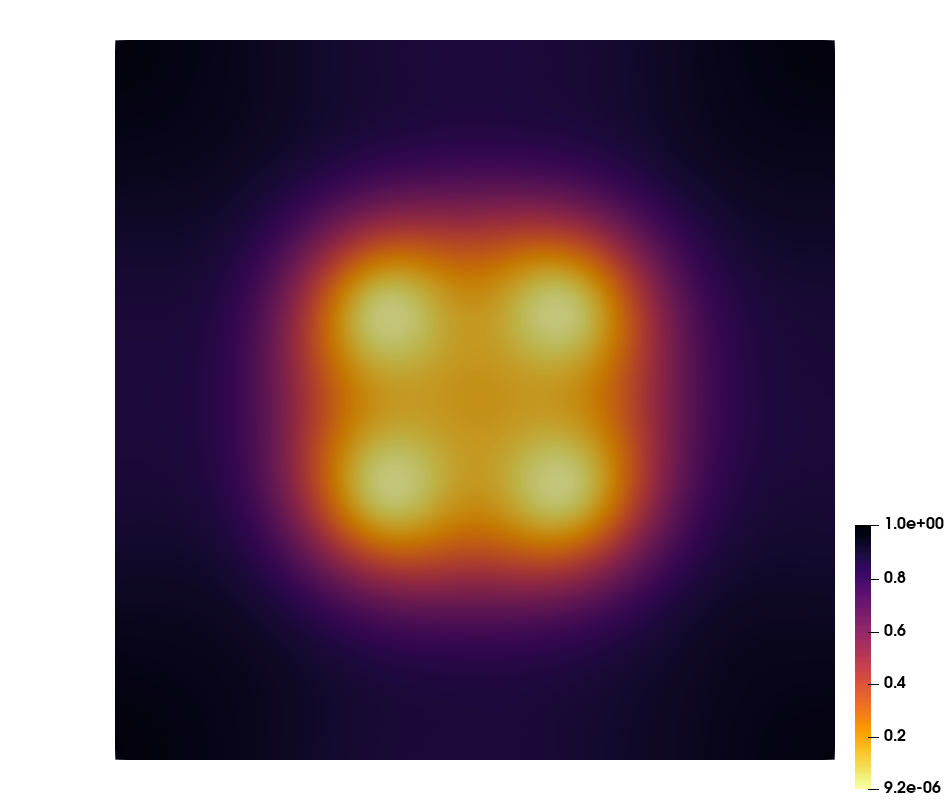}
    \end{minipage}
    \begin{minipage}{0.15\textwidth}
    \includegraphics[scale=0.07]{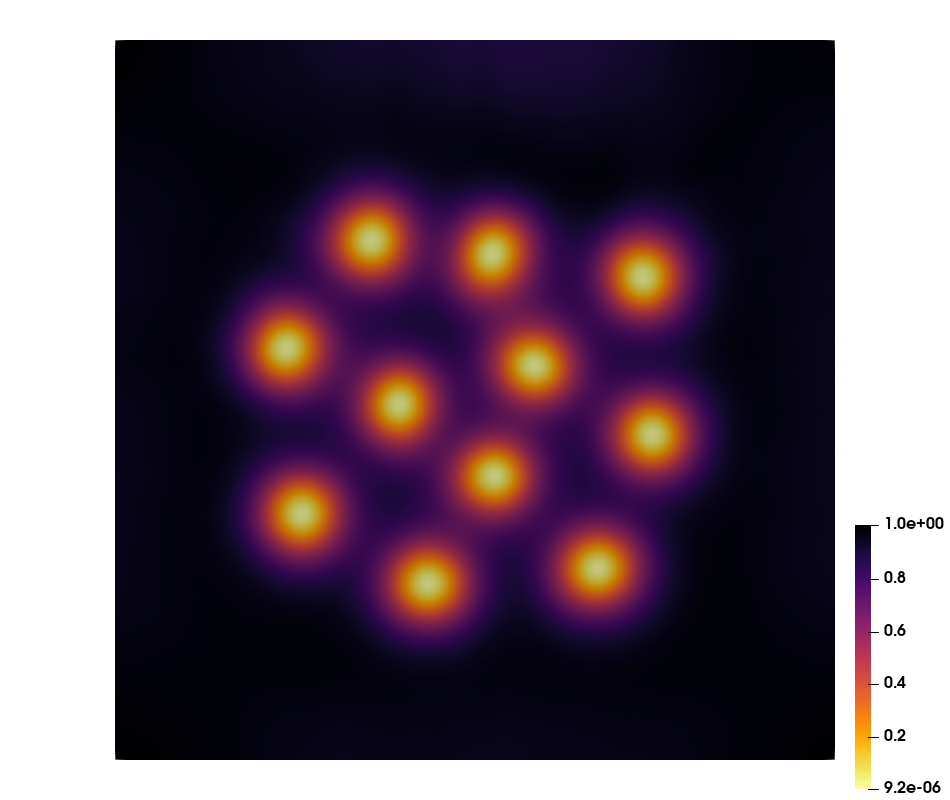}
    \end{minipage}
    \begin{minipage}{0.15\textwidth}
    \includegraphics[scale=0.07]{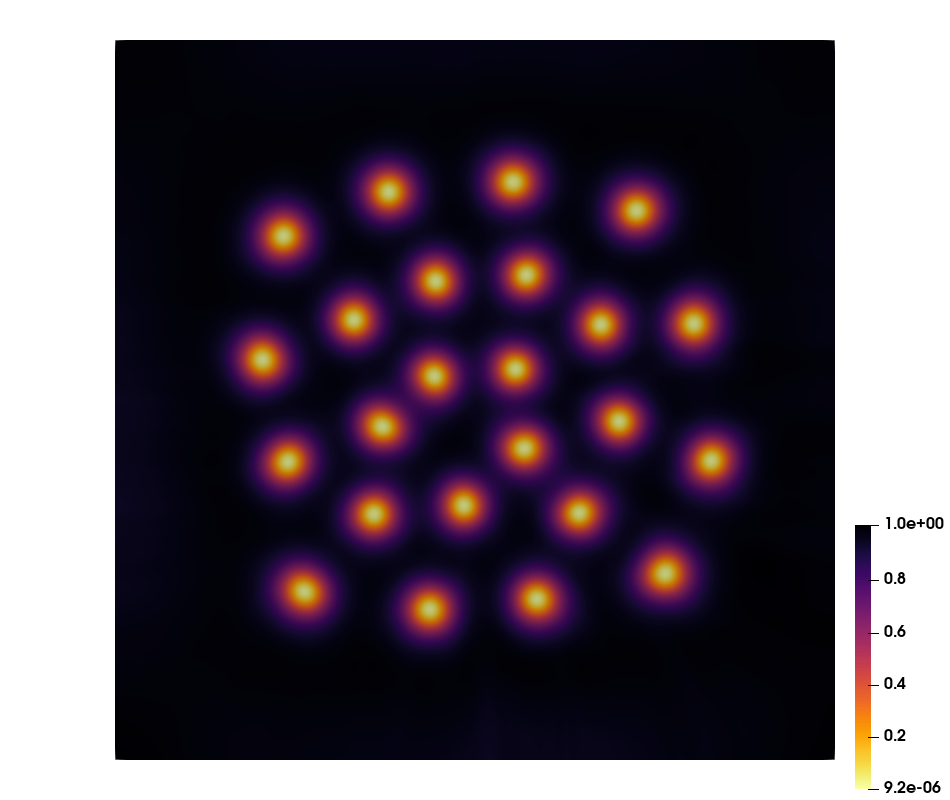}
    \end{minipage}
    \begin{minipage}{0.15\textwidth}
    \includegraphics[scale=0.07]{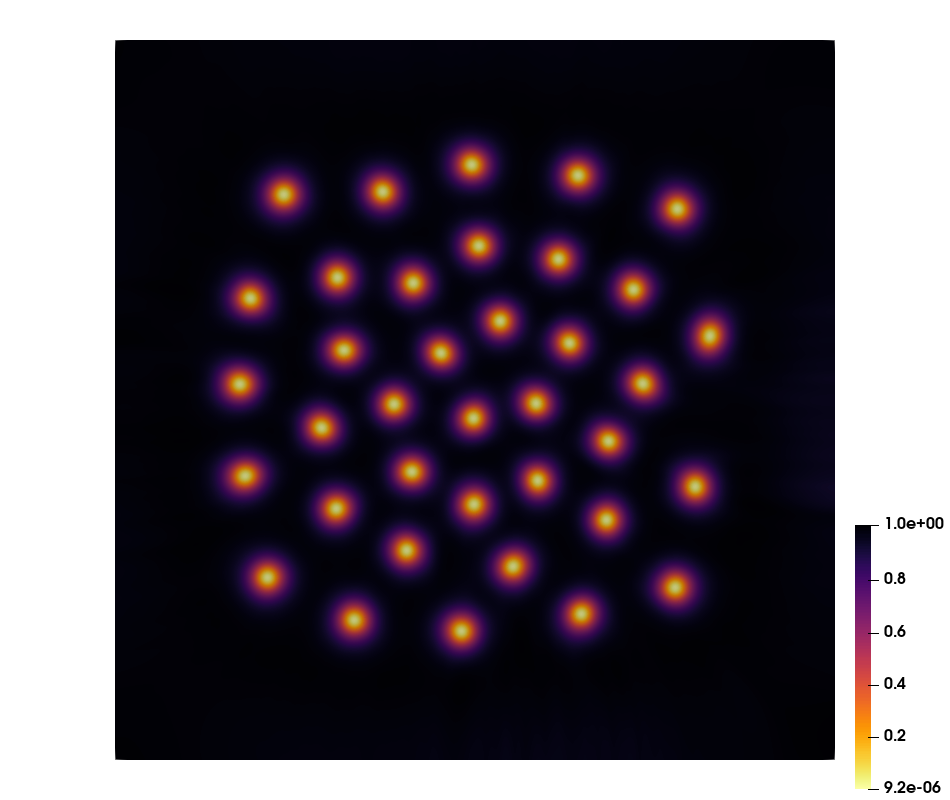}
    \end{minipage}
    \begin{minipage}{0.15\textwidth}
    \includegraphics[scale=0.07]{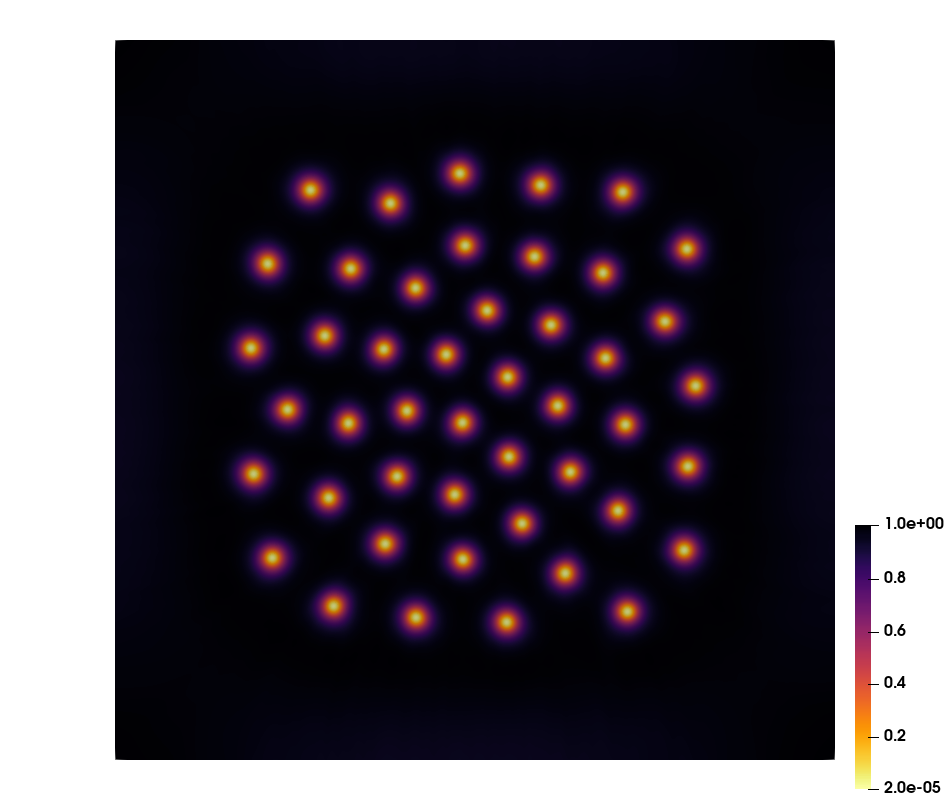}
    \end{minipage} 
    \\
    \begin{minipage}{0.1\textwidth}
    \centering
    \texttt{GL-R2} (NN)
    \end{minipage}
    \begin{minipage}{0.15\textwidth}
    \includegraphics[scale=0.07]{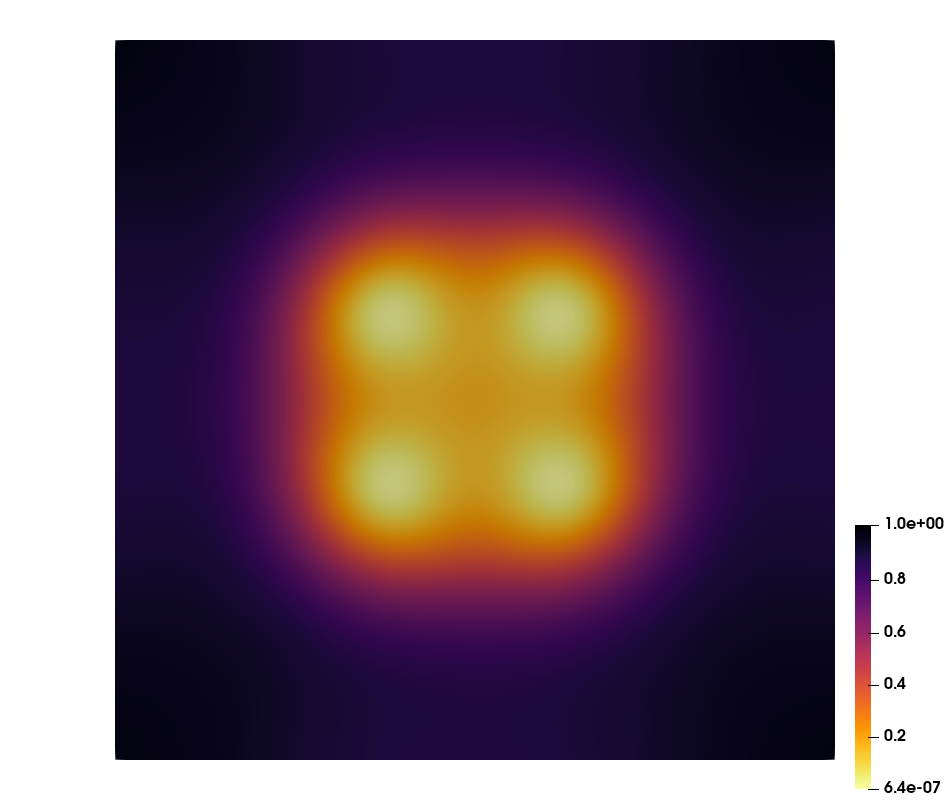}
    \end{minipage}
    \begin{minipage}{0.15\textwidth}
    \includegraphics[scale=0.07]{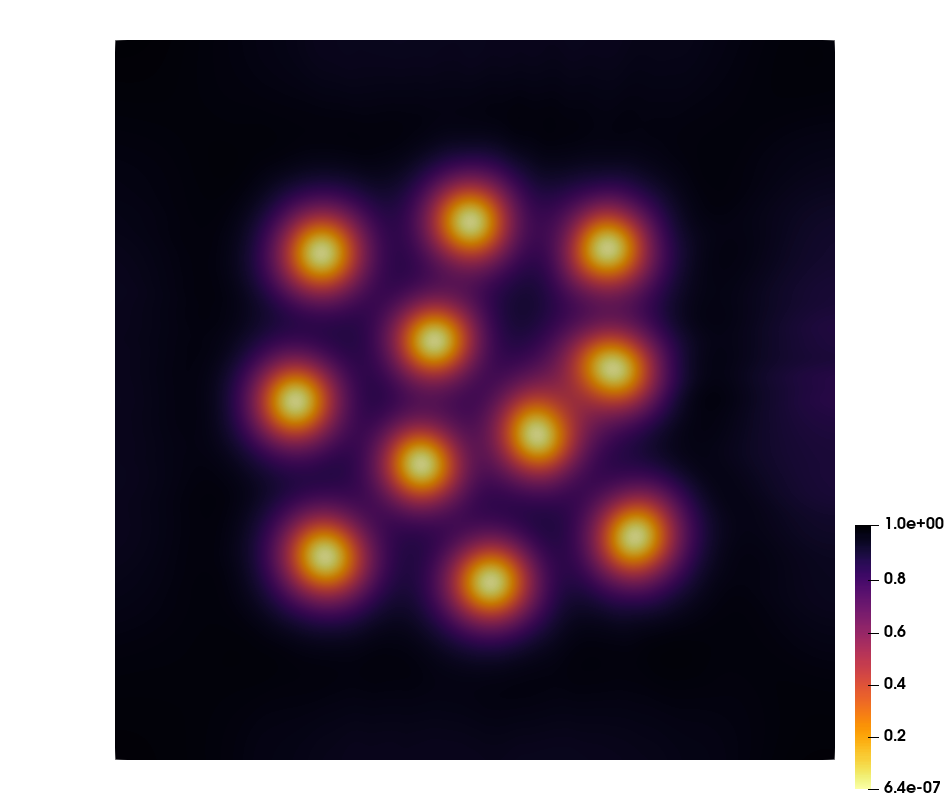}
    \end{minipage}
    \begin{minipage}{0.15\textwidth}
    \includegraphics[scale=0.07]{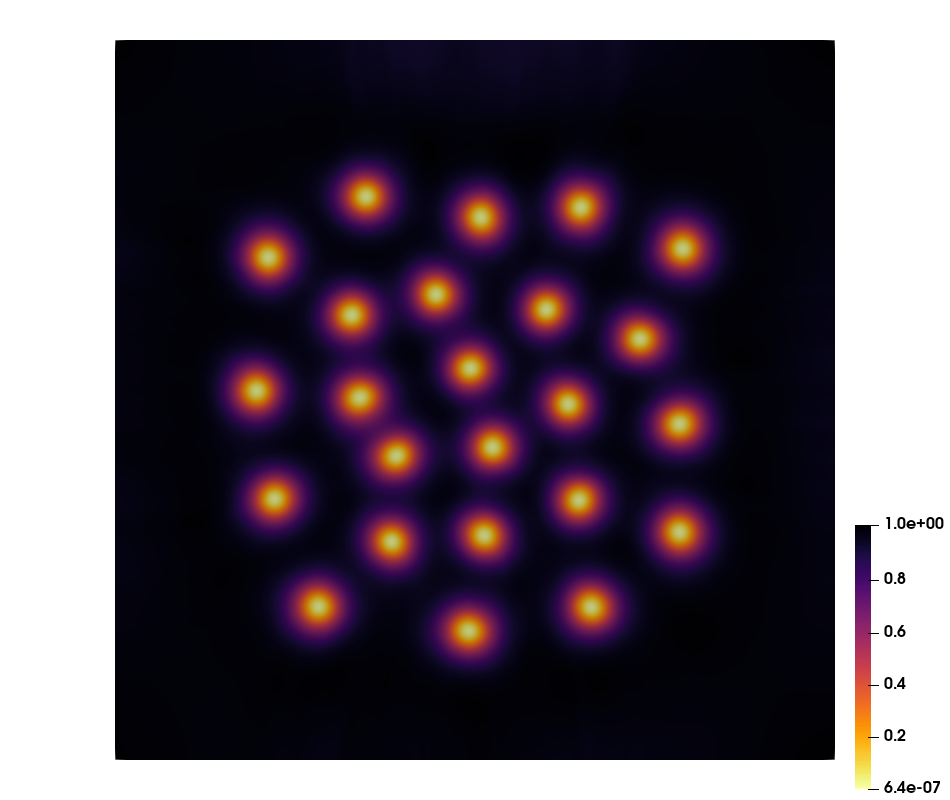}
    \end{minipage}
    \begin{minipage}{0.15\textwidth}
    \includegraphics[scale=0.07]{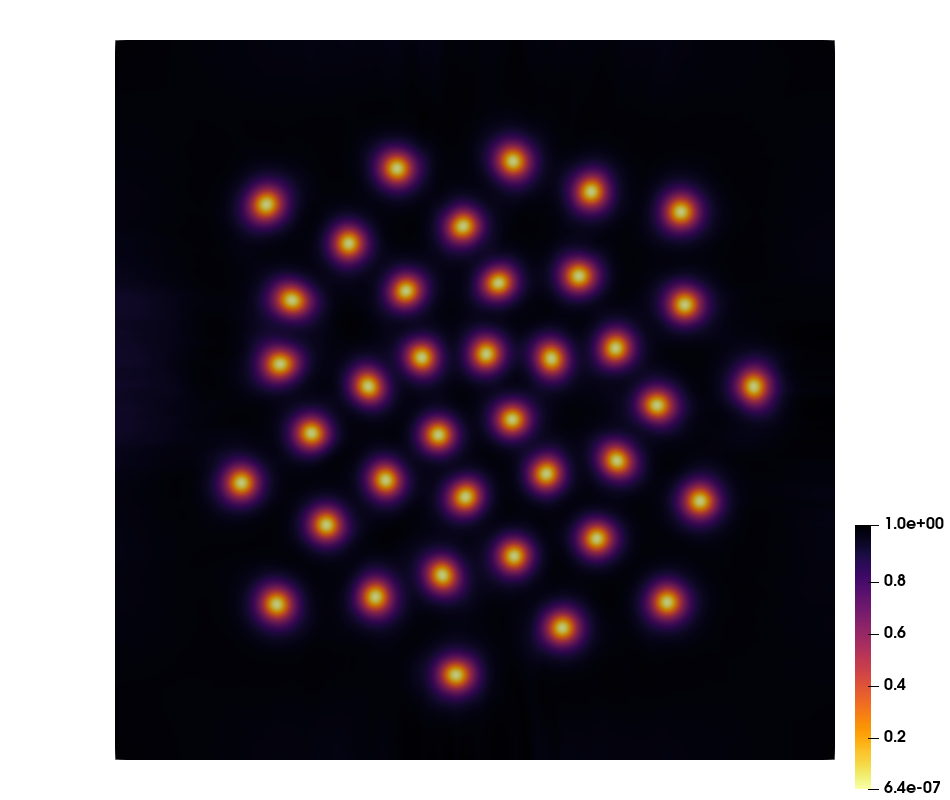}
    \end{minipage}
    \begin{minipage}{0.15\textwidth}
    \includegraphics[scale=0.07]{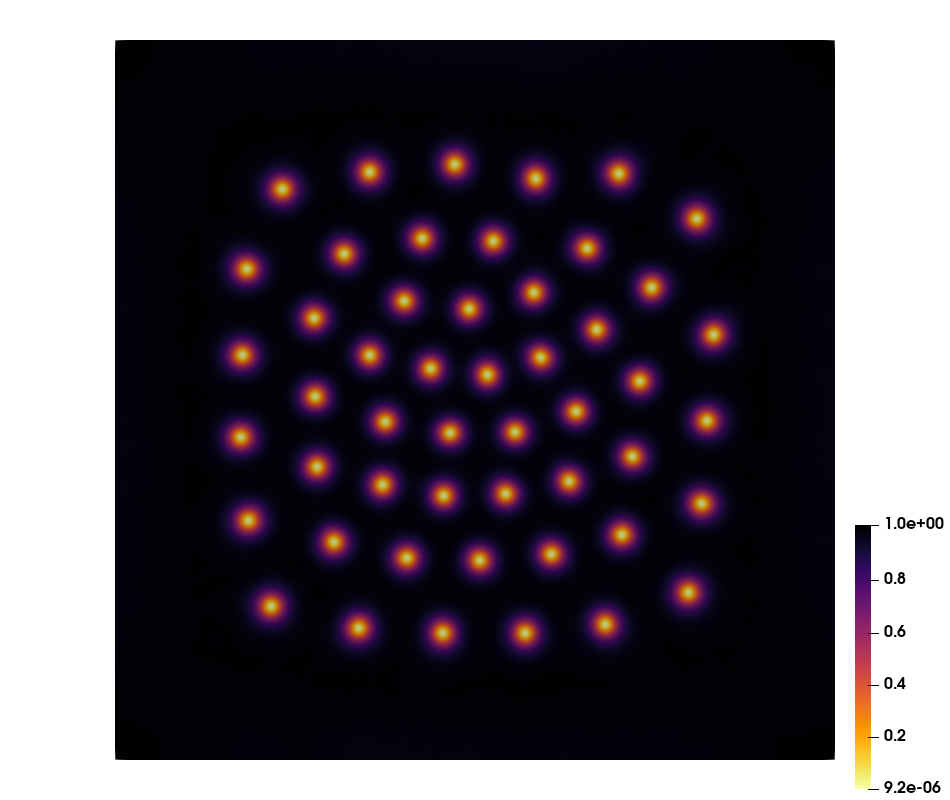}
    \end{minipage} 
    \\
    \begin{minipage}{0.1\textwidth}
    \centering
    \texttt{GL-R1} (hyb.)
    \end{minipage}
    \begin{minipage}{0.15\textwidth}
    \includegraphics[scale=0.07]{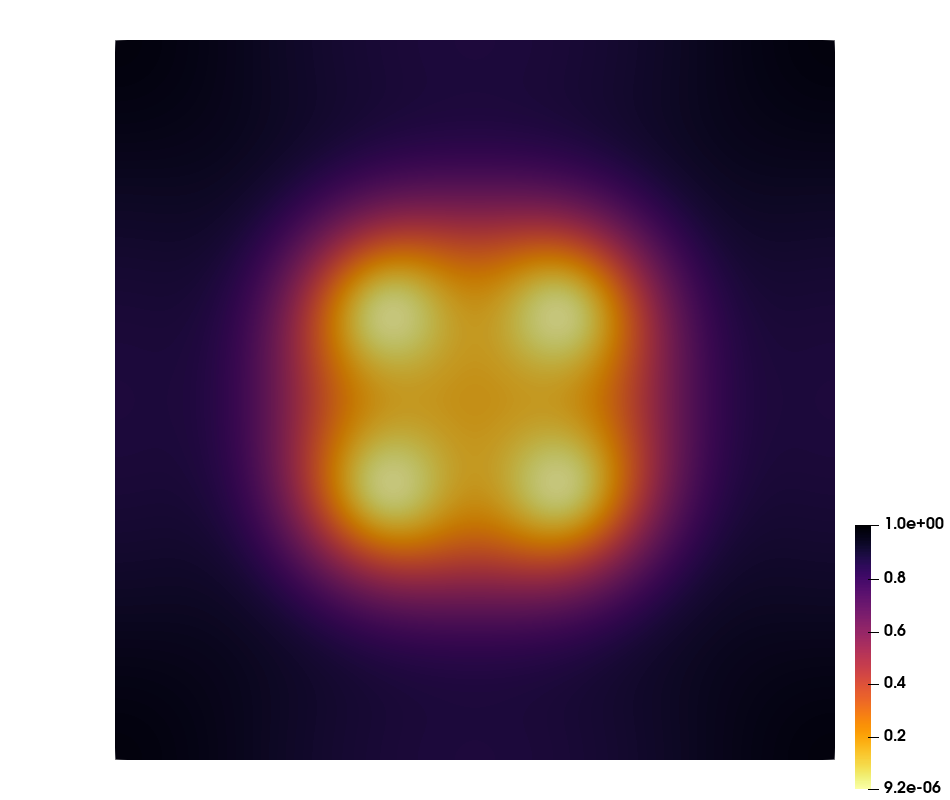}
    \end{minipage}
    \begin{minipage}{0.15\textwidth}
    \includegraphics[scale=0.07]{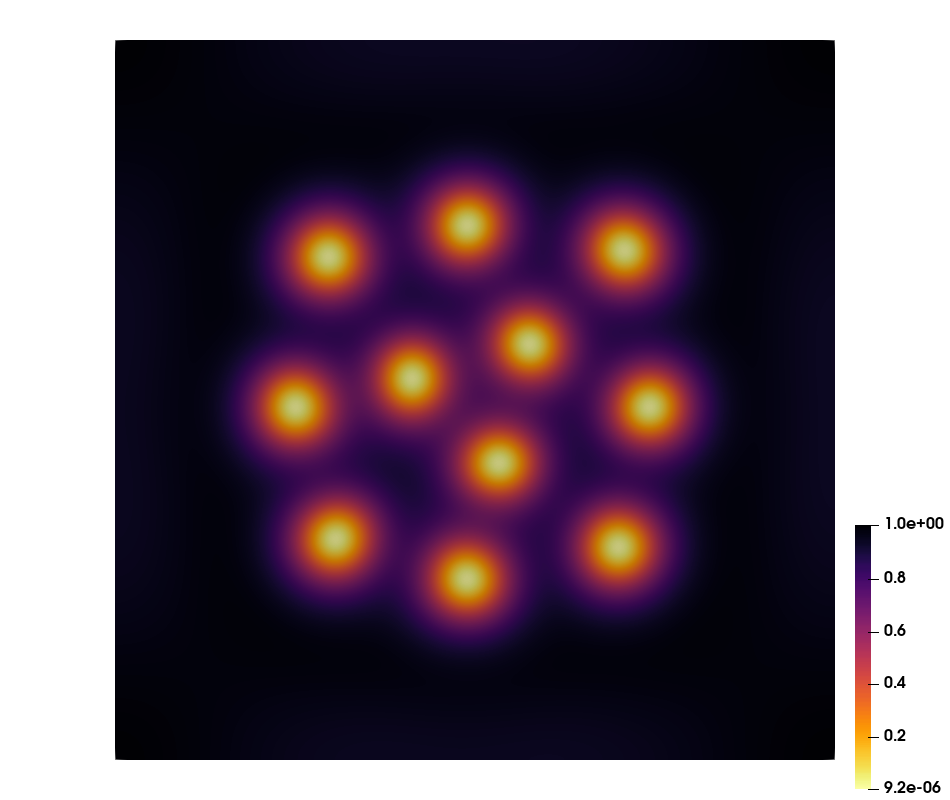}
    \end{minipage}
    \begin{minipage}{0.15\textwidth}
    \includegraphics[scale=0.07]{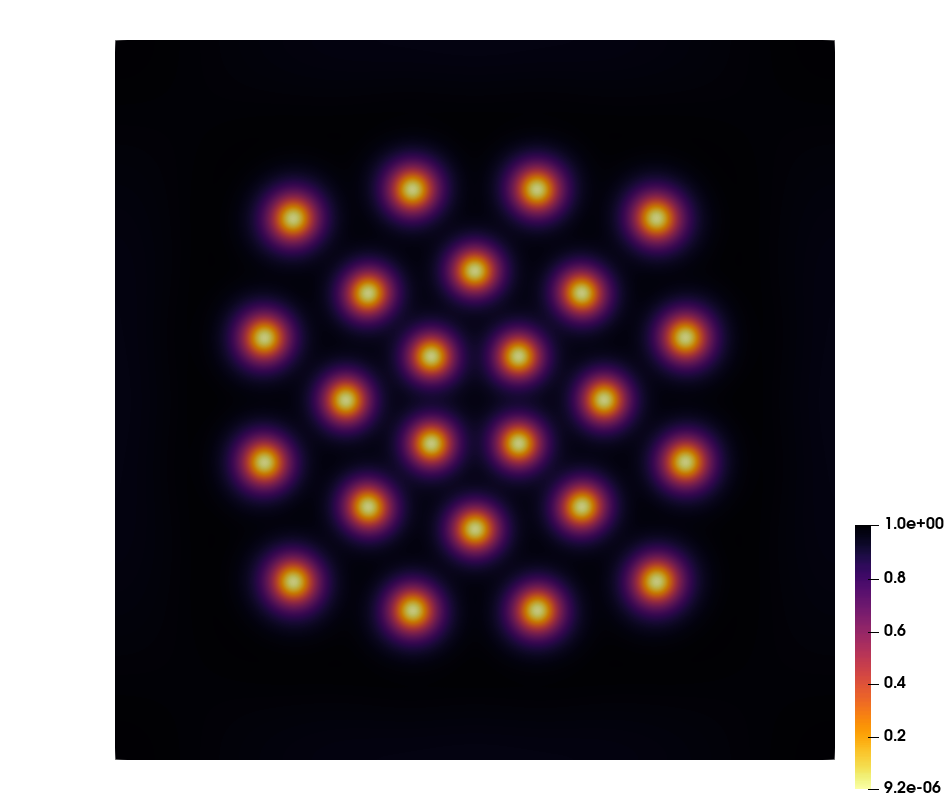}
    \end{minipage}
    \begin{minipage}{0.15\textwidth}
    \includegraphics[scale=0.07]{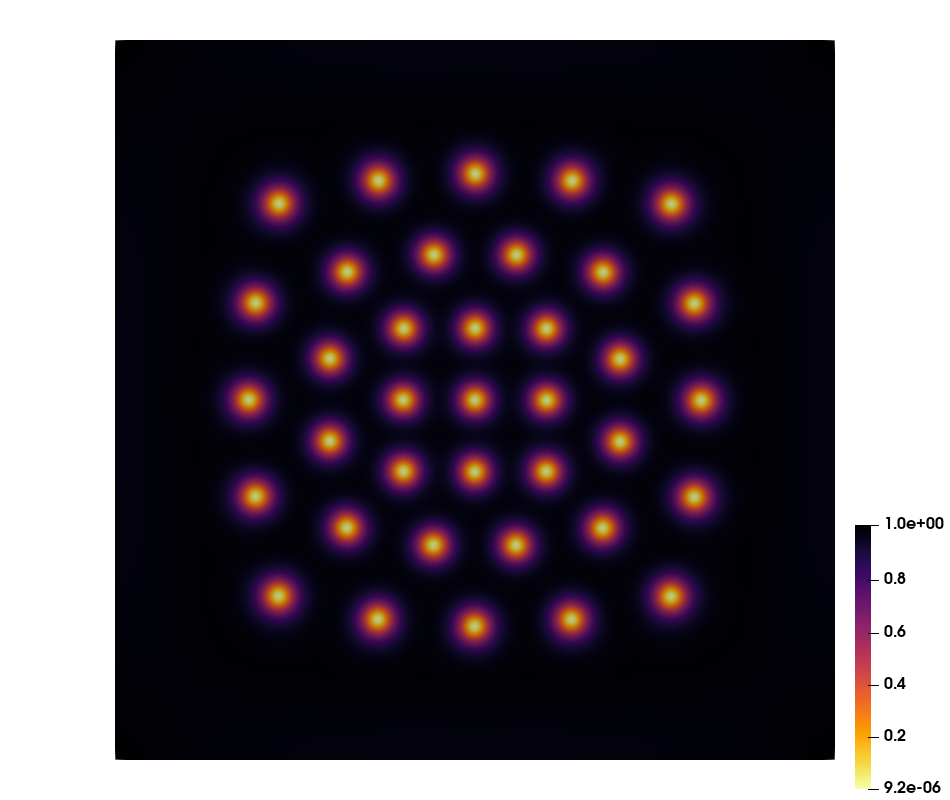}
    \end{minipage}
    \begin{minipage}{0.15\textwidth}
    \includegraphics[scale=0.07]{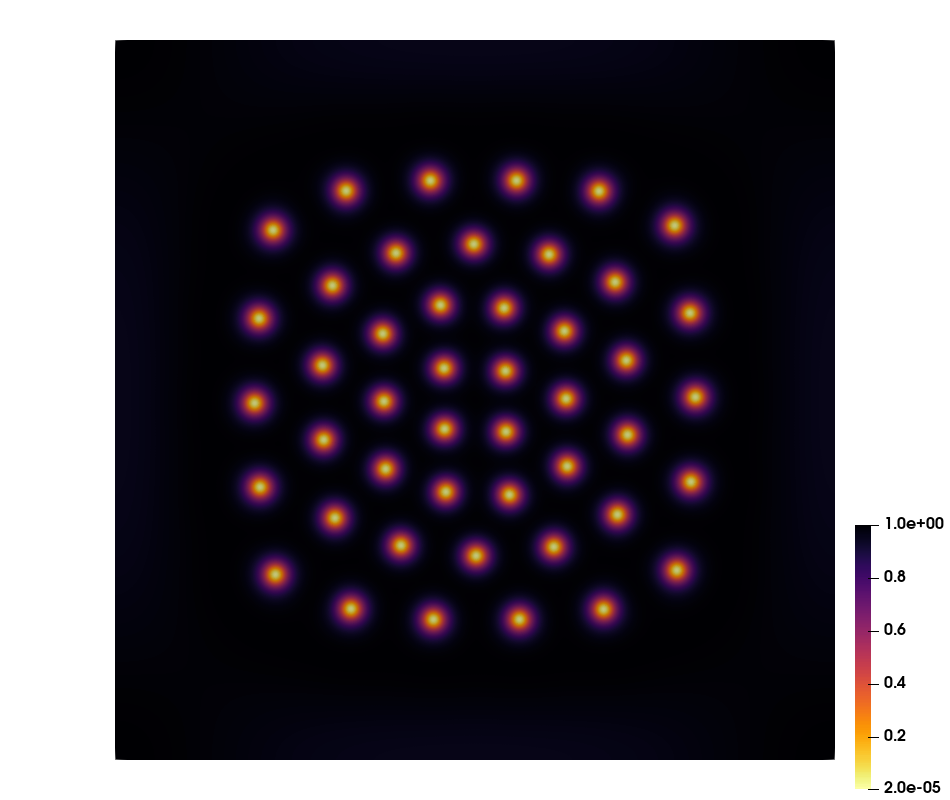}
    \end{minipage} 
    \\
    \begin{minipage}{0.1\textwidth}
    \centering
    \texttt{GL-R2} (hyb.)
    \end{minipage}
    \begin{minipage}{0.15\textwidth}
    \includegraphics[scale=0.07]{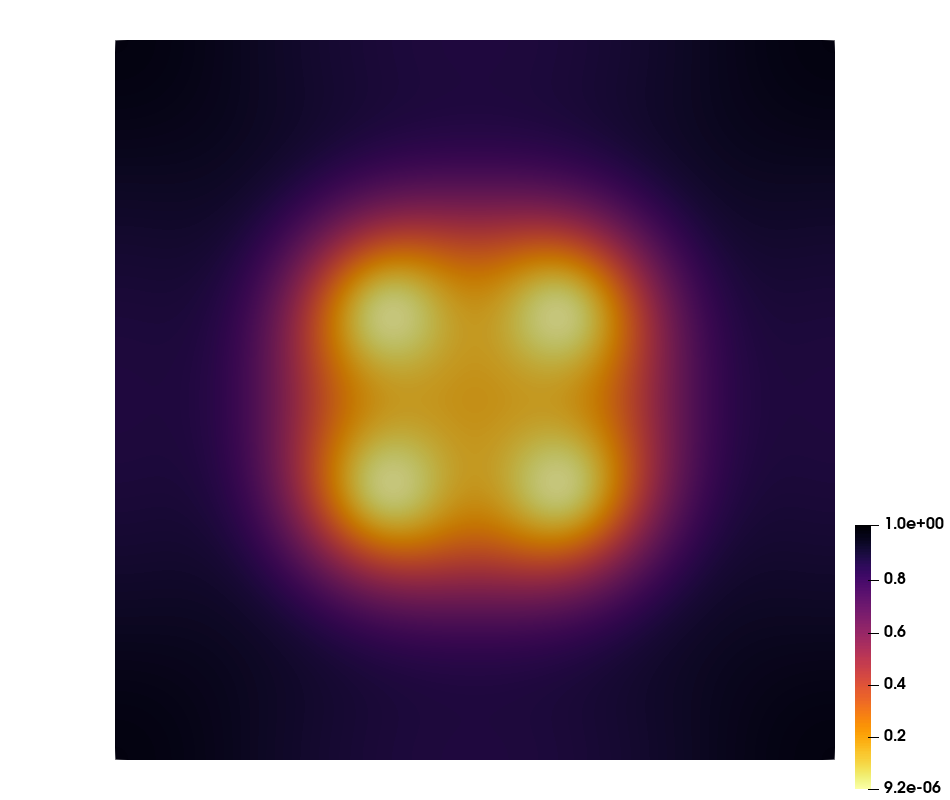}
    \end{minipage}
    \begin{minipage}{0.15\textwidth}
    \includegraphics[scale=0.07]{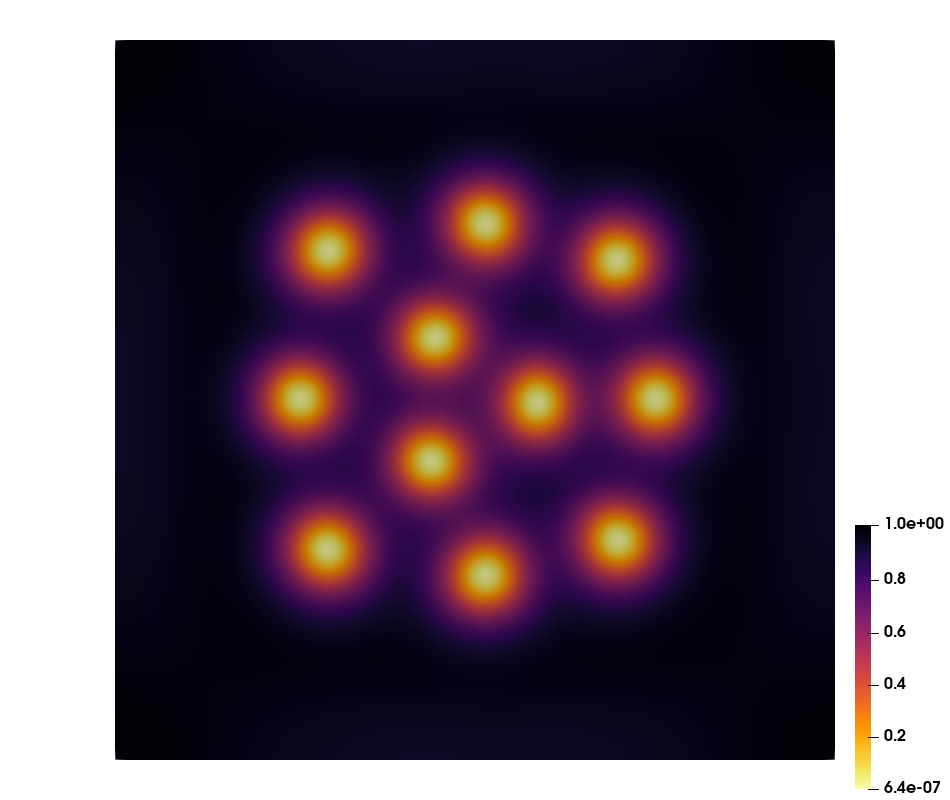}
    \end{minipage}
    \begin{minipage}{0.15\textwidth}
    \includegraphics[scale=0.07]{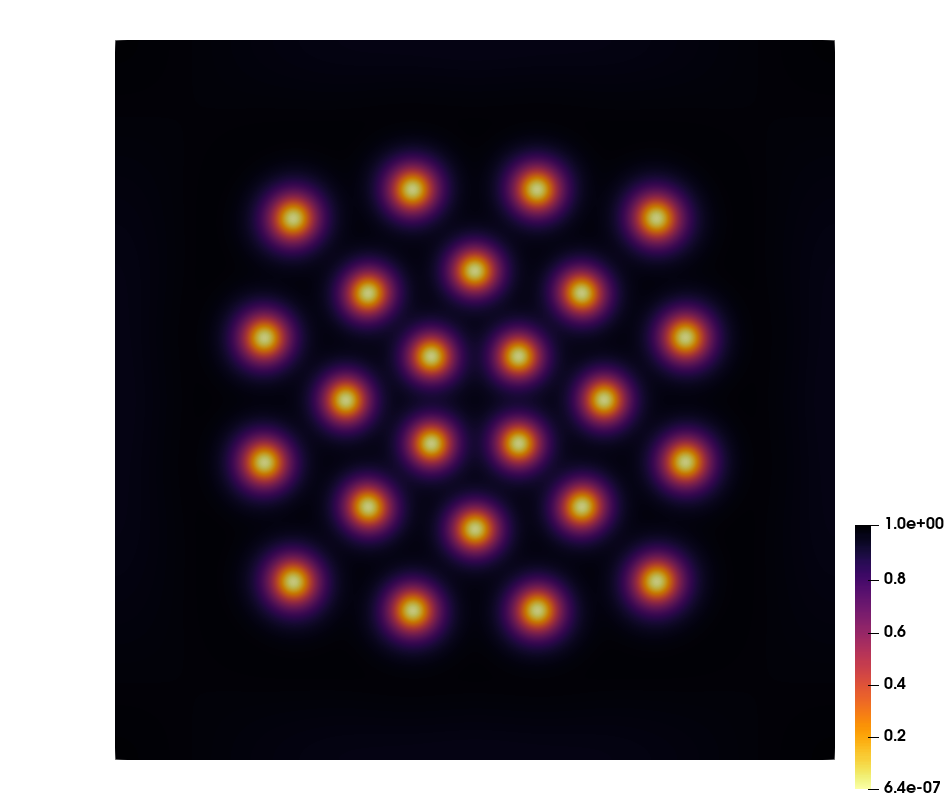}
    \end{minipage}
    \begin{minipage}{0.15\textwidth}
    \includegraphics[scale=0.07]{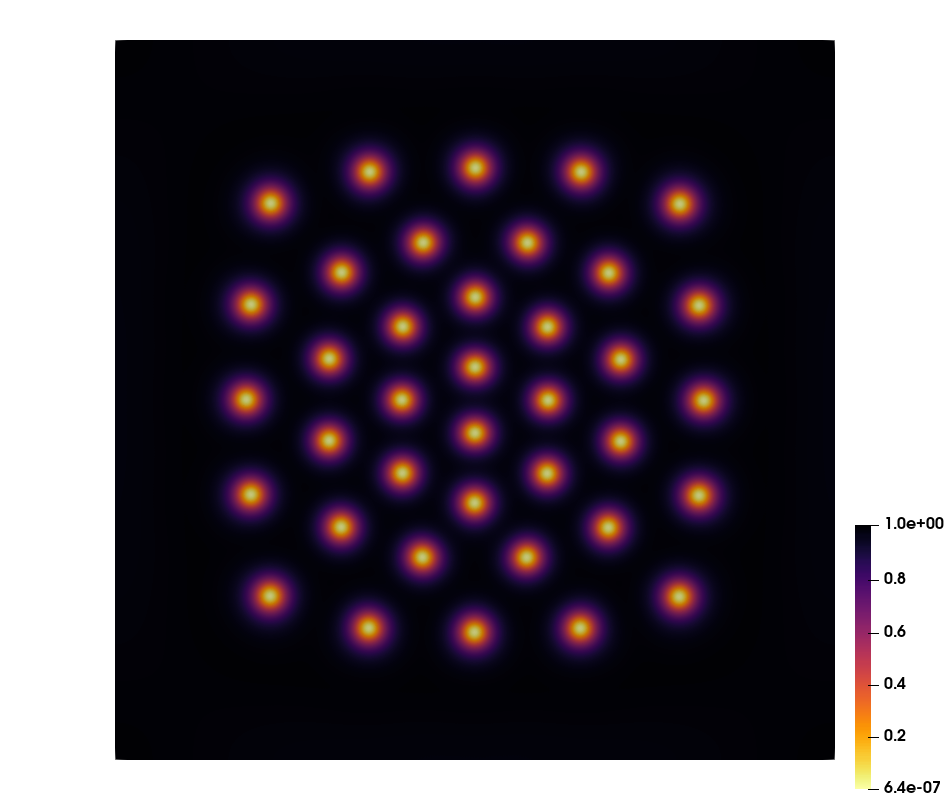}
    \end{minipage}
    \begin{minipage}{0.15\textwidth}
    \includegraphics[scale=0.07]{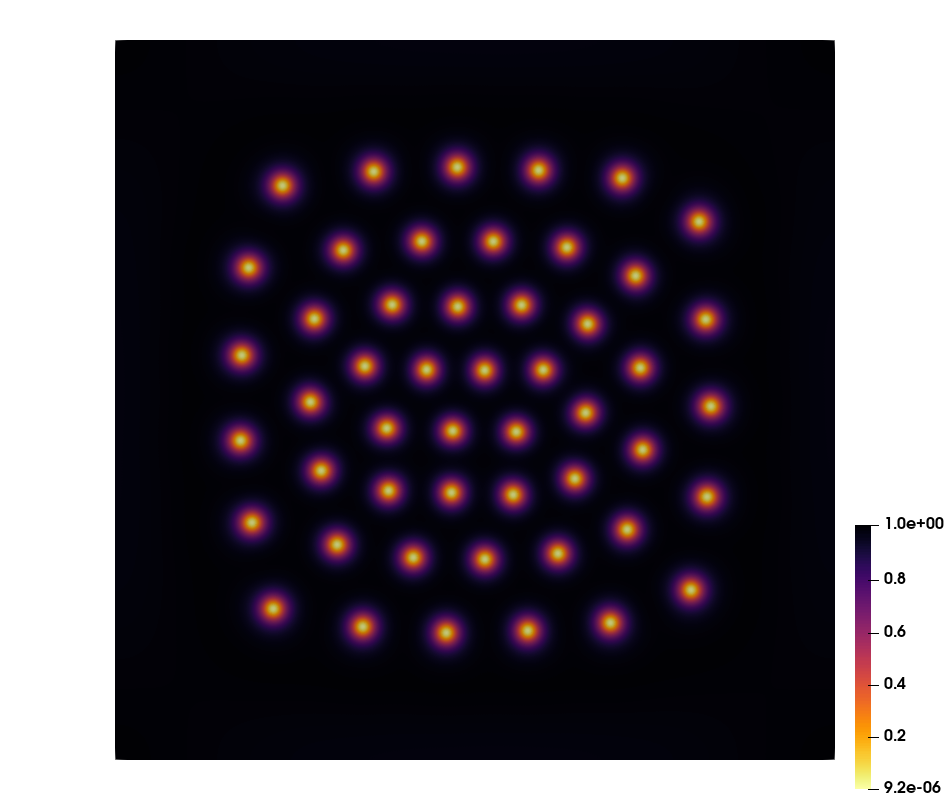}
    \end{minipage}
    \caption{Densities $|u|^2$ of computed minimizers for $\kappa = 10 ,25, 50, 75, 100$ (left to right), corresponding to Table~\ref{tab:energy_only_ord}.
    Dark indicates values close to 1; light indicates values close to 0. 
    First row: reference minimizer computed with the classical FE iterative solver and best heuristic initial values $\varphi_j$ (best). Second and third row: pure NN models
    \texttt{GL-R1} 
    (NN) and
    \texttt{GL-R2}
    (NN). 
    Fourth and fifth row: hybrid approach
    \texttt{GL-R1}
    (hyb.) and
    \texttt{GL-R2}
    (hyb.)
    using initial values from the NN.}
    \label{fig:minimizers_oo}
\end{figure}

The first block in Table \ref{tab:energy_only_ord} shows the energy levels $E(u)$ of the five computed minimizers using the classical FE iterative solver with the empirically chosen initial values $\varphi_j$, $j=1,2,3,4,5$. The bold entries indicate, for each value of $\kappa$, the lowest energy achieved by this procedure; these serve as reference approximations for comparison with the NN and the hybrid approach.

For the NN-based method, we trained two models, \texttt{GL-R1} and \texttt{GL-R2}. Their architectures and training strategies are summarized in Table~\ref{tab:parameters}. For both models, we use the SwiGLU block type and the main difference lies in the chosen learning rates and batch size. The training required approximately 3.5 hours for for \texttt{GL-R1} and 6.5 hours for \texttt{GL-R2} on the HoreKa machine at KIT. The resulting energy levels of the pure NN approach are shown in the second block of Table~\ref{tab:energy_only_ord}. For $\kappa = 10,25,50$, the NN approximations yield energy levels close to, but slightly above, the reference values. In contrast, for $\kappa = 75,100$ both models achieve lower energy levels than those obtained with the classical approach using the prescribed heuristic initial values. This suggests that the NN approach identifies better candidates for the global GL minimizer.

We then use our hybrid scheme, which takes the output of the NN as an initial guess for the classical FE iterative solver. The final energy levels after an additional iteration with the FE solver are shown in the third block of Table~\ref{tab:energy_only_ord}. Both models show comparable performance and consistently match or even outperform the reference approximations due to lower energy levels. In particular, they yield strictly lower energy states in several cases, indicating the successful identification of low energy states, which is the primary objective of our approach.

In Figure \ref{fig:minimizers_oo}, we present the densities $|u|^2$ of the computed minimizers. In all cases, we observe the formation of a vortex pattern, namely the Abrikosov vortex lattice \cite{Abr04}. As $\kappa$ is increased, the pattern becomes more intricate: the number of vortices grows while their size decreases. The minimizers computed with the classical FE approach and with the hybrid approach exhibit symmetric vortex patterns, whereas the pure NN results show a certain displacement of the vortices within the lattice. At the same time, the pure NN approach still captures the overall structure of the pattern, as well as the approximate size and number of vortices. Moreover, for the same value of $\kappa$ (e.g., $\kappa=100$), even small differences in the computed energy correspond to different vortex configurations and therefore to different local minimizers. This suggests that, with the classical approach and heuristic initial values, the global GL minimizer is not recovered, while the hybrid approach appears to provide at least a strong candidate for the global minimizer.

\begin{table}[h!]
    \centering
    \small
    \begin{tabular}{c c c c c c}
        \toprule
        $\kappa$ & 10 & 25 & 50 & 75 & 100 \\
        \midrule
        $\varphi_1$ & \textbf{40} & 60 & 172 & 10279 & 10882 \\
        $\varphi_2$ & 199 & 339 & 1137 & 4373 & 24572 \\
        $\varphi_3$ & 66 & \textbf{1166} & \textbf{4034} & 8147 & 11022 \\
        $\varphi_4$ & 91 & 1214 & 4090 & \textbf{21718} & \textbf{17831} \\
        $\varphi_5$ & 86 & 379 & 21874 & 47145 & 37079 \\
        \midrule
        \texttt{GL-R1} (hyb.) & 52 & 1314 & 7483 & 11912 & 4537 \\
        \texttt{GL-R2} (hyb.) & 47 & 5541 & 7961 & 12452 & 24568 \\
        \bottomrule
    \end{tabular}
        \caption{Number of iterations to compute the minimizers with initial value $u^0=\varphi_j$ and the hybrid-approach up to tolerance $\mathrm{tol}=10^{-12}$.}
    \label{tab:iterations_only_ord}
\end{table}

Table \ref{tab:iterations_only_ord} shows the iterations that the FE solver needed to converge up to the prescribed tolerance of $10^{-12}$ for both, the classical FE and the hybrid approach. As the table shows there is no clear indication whether either approach is favorable in that regard. Note that the numbers should only be compared for cases where the same minimizing state is achieved by both methods. While the performance is not the focus of our approach, we nonetheless implemented a refinement strategy to improve on the number of iterations needed after the NN stage. For this, we treat the trained NN as a pre-trained model that is then refined for a limited time and for the given value of $\kappa$. This can be interpreted as specializing the general model for this specific $\kappa$. In our tests, we observed clear indications that this procedure can indeed reduce the required number of iterations. However, a further detailed investigation of possibilities to reduce the number of iterations is beyond the scope of this work and left for future research.

\subsection{Neural networks for the full GL model} Next, we consider the full GL model \eqref{eq:energy_functional} with the external magnetic field
\begin{equation} \label{eq:hext}
    \MagF(x) = 2\sqrt{2} \pi \sin(\pi x_1) \sin (\pi_2), \quad x = (x_1,x_2) \in \Omega,
\end{equation}
and again choose the GL parameters as $\kappa = 10,25,50,75,100$. As before, we note that the NNs are not restricted to these values of $\kappa$, but instead are trained for a whole range.

\begin{table}[h!]
\centering
\small
\begin{tabular}{c c c c c c}
\toprule
$\kappa$ & 10 & 25 & 50 & 75 & 100 \\
\midrule
$\varphi_1$ & \textbf{0.10440628} & 0.08514562 & \textbf{0.04148069} & 0.04560323 & \textbf{0.02615202} \\
$\varphi_2$ & 0.11500423 & 0.07371087 & 0.04185499 & 0.03776825 & 0.02854666 \\
$\varphi_3$ & 0.10440628 & 0.06491247 & 0.04148069 & 0.03779534 & 0.02659625 \\
$\varphi_4$ & 0.10440628 & \textbf{0.06491247} & 0.04148069 & \textbf{0.03446474} & 0.02660840 \\
$\varphi_5$ & 0.10440628 & 0.08287191 & 0.04781070 & 0.03597329 & 0.02707063 \\
\midrule
\makebox[18mm][l]{\texttt{GL-F1} (NN)}
 & 0.10461688 & 0.06576285 & 0.04211054 & 0.03220804 & 0.02639166 \\
\makebox[18mm][l]{\texttt{GL-F2} (NN)}
& 0.10446705 & 0.06545582 & 0.04194608 & 0.03210672 & 0.02640819 \\
\midrule
\makebox[18mm][l]{\texttt{GL-F1} (hyb.)}
& 0.10440628 & 0.06445608 & 0.04148069 & 0.03171263 & 0.02615239 \\
\makebox[18mm][l]{\texttt{GL-F2} (hyb.)}
& 0.10440628 & 0.06445609 & 0.04148069 & 0.03171263 & 0.02614767 \\
\bottomrule
\end{tabular}
\caption{Energy levels $E(u,\MagPot)$ of computed energy minimizers for the full GL model of Section \ref{sec:experiment_reduced_model}. First block: classical FE iterative solver with heuristic initial values $\varphi_j$ for $u$ and $0$ for $\MagPot$. Second block: pure NN models
\texttt{GL-F1}
(NN) and
\texttt{GL-F2}
(NN). Third block: hybrid approach
\texttt{GL-F1}
(hyb.) and
\texttt{GL-F2}
(hyb.) using initial values from the NN.}
\label{tab:energy_full}
\end{table}

\begin{figure}[h!]
    \flushleft
    \begin{minipage}{0.1\textwidth}
    \centering
    $\varphi_j$ (best)
    \end{minipage}
    \begin{minipage}{0.15\textwidth}
    \includegraphics[scale=0.07]{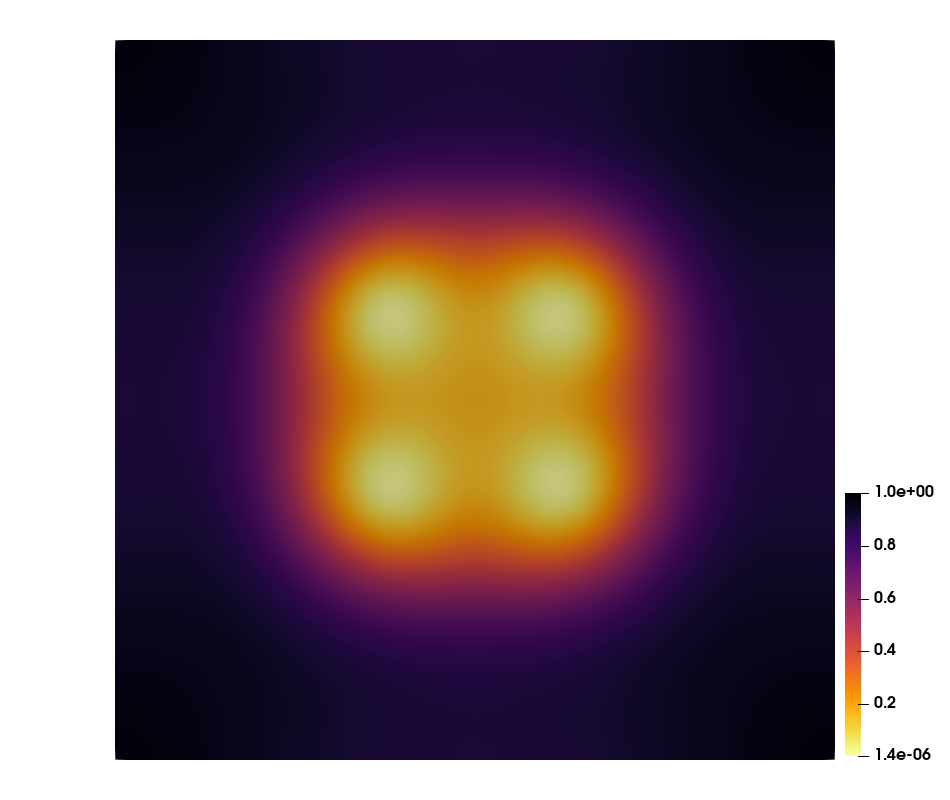}
    \end{minipage}
    \begin{minipage}{0.15\textwidth}
    \includegraphics[scale=0.07]{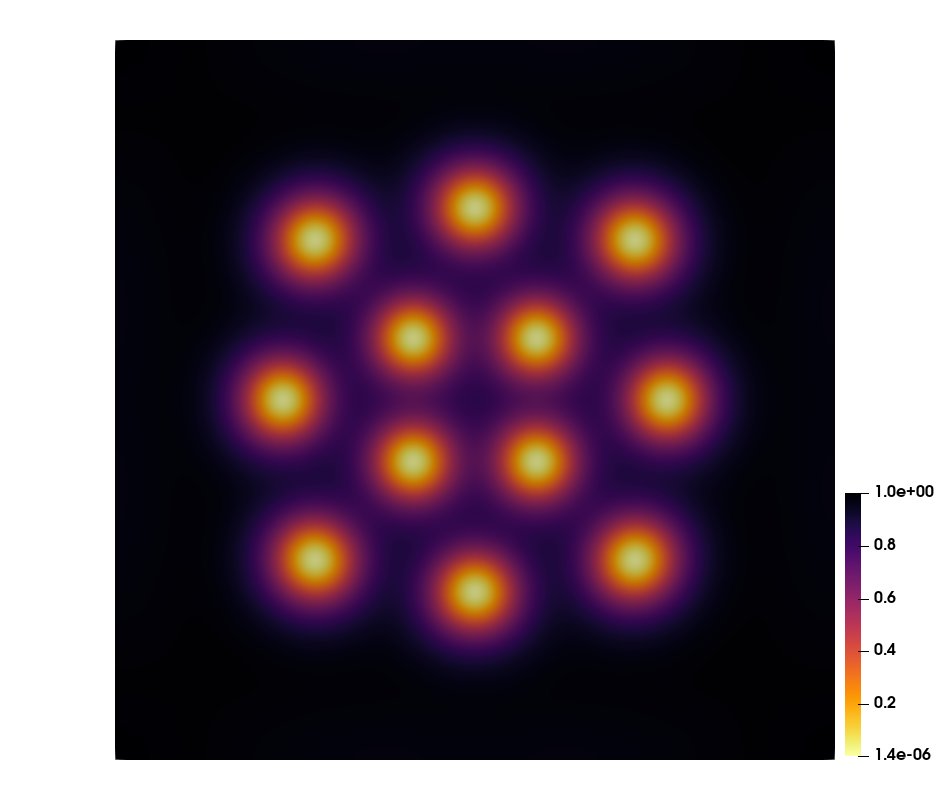}
    \end{minipage}
    \begin{minipage}{0.15\textwidth}
    \includegraphics[scale=0.07]{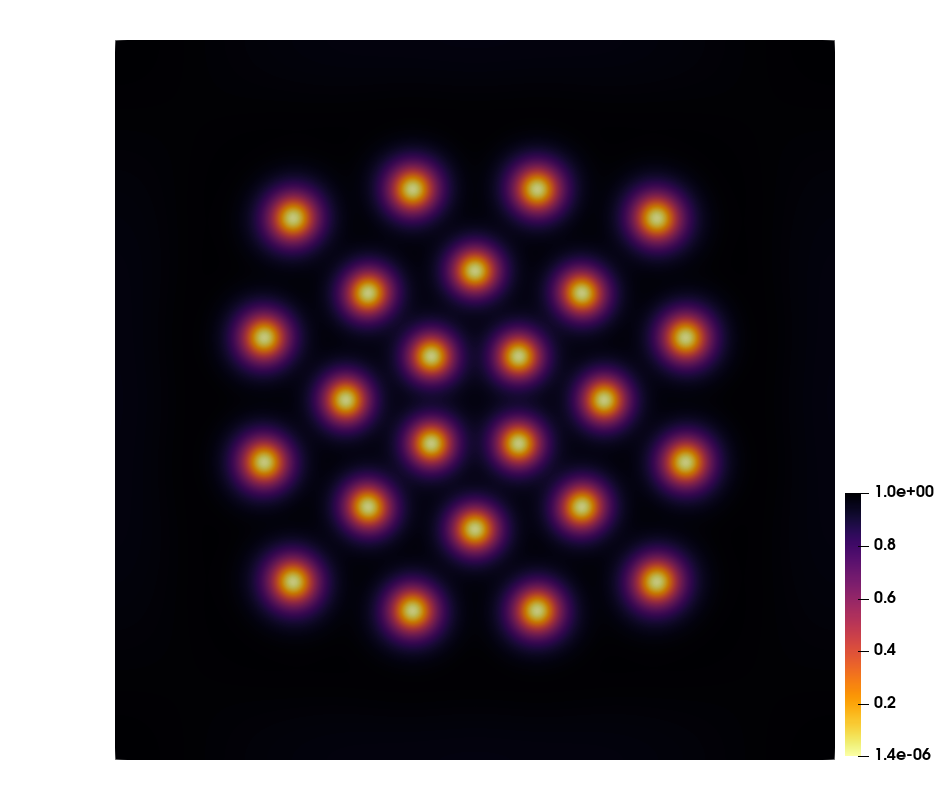}
    \end{minipage}
    \begin{minipage}{0.15\textwidth}
    \includegraphics[scale=0.07]{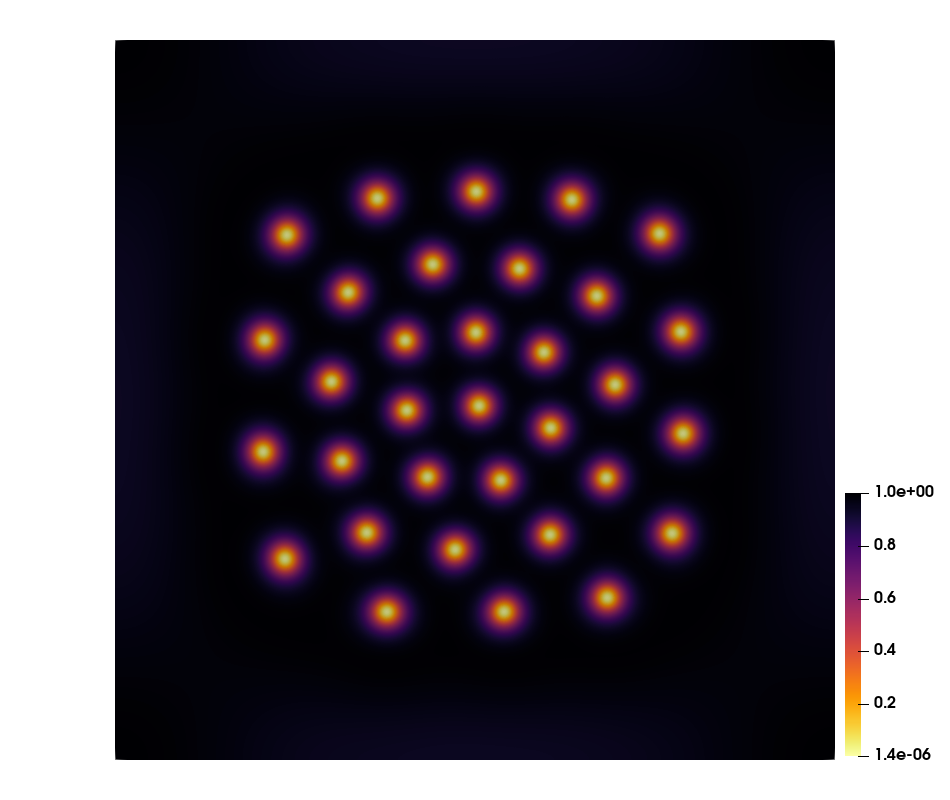}
    \end{minipage}
    \begin{minipage}{0.15\textwidth}
    \includegraphics[scale=0.07]{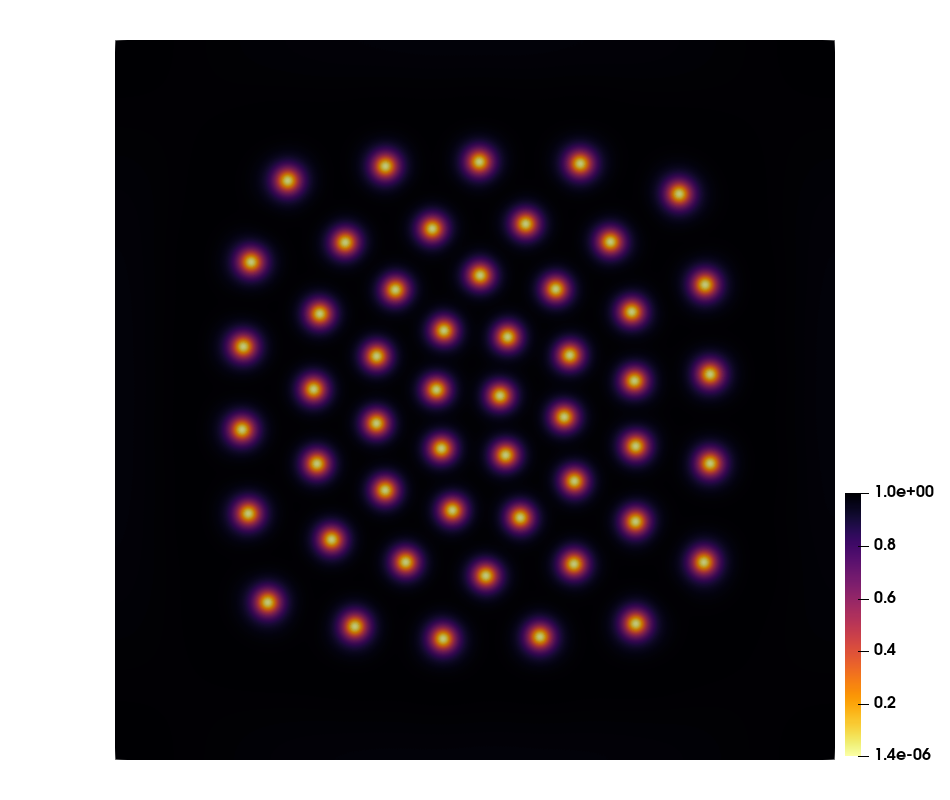}
    \end{minipage} 
    \\
    \begin{minipage}{0.1\textwidth}
    \centering
    \texttt{GL-F1} (NN)
    \end{minipage}
    \begin{minipage}{0.15\textwidth}
    \includegraphics[scale=0.07]{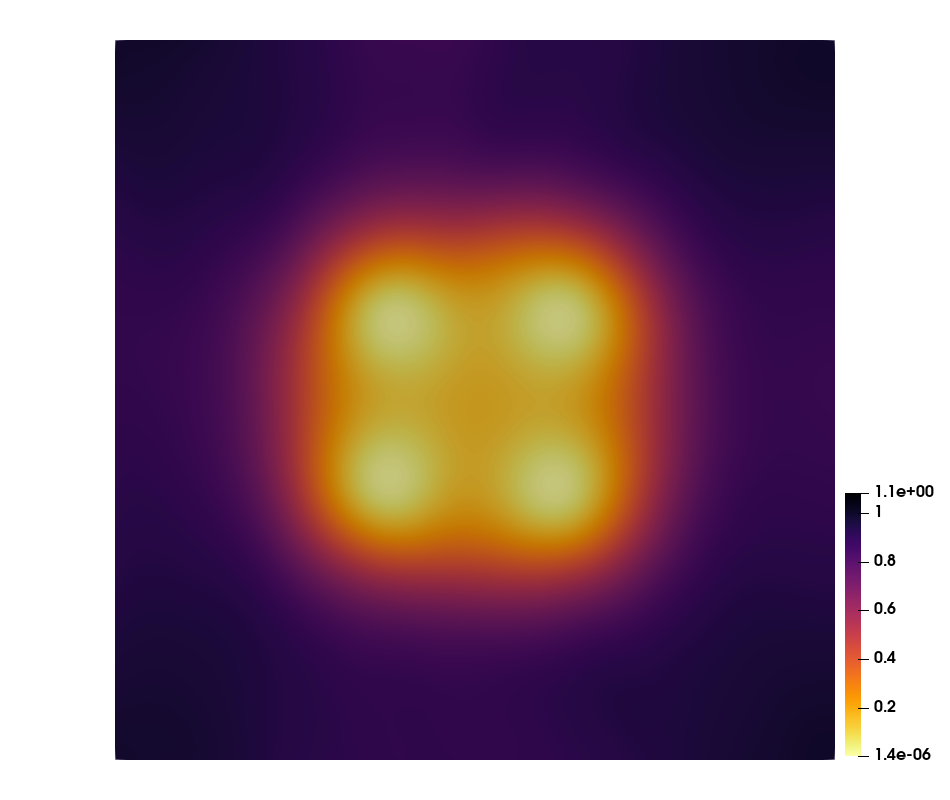}
    \end{minipage}
    \begin{minipage}{0.15\textwidth}
    \includegraphics[scale=0.07]{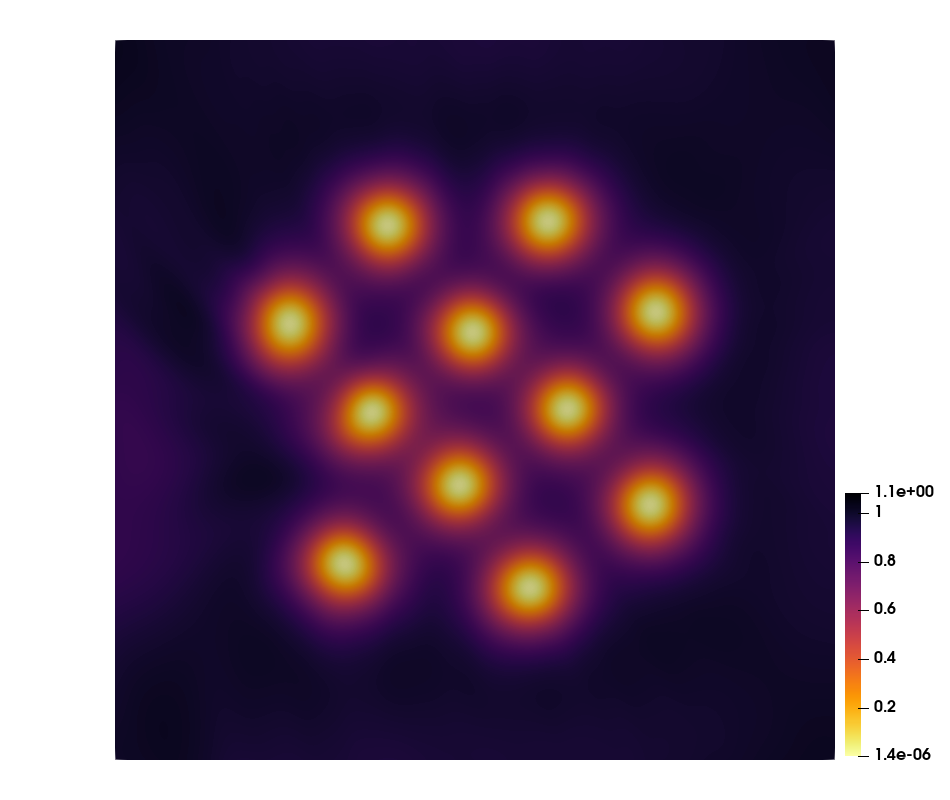}
    \end{minipage}
    \begin{minipage}{0.15\textwidth}
    \includegraphics[scale=0.07]{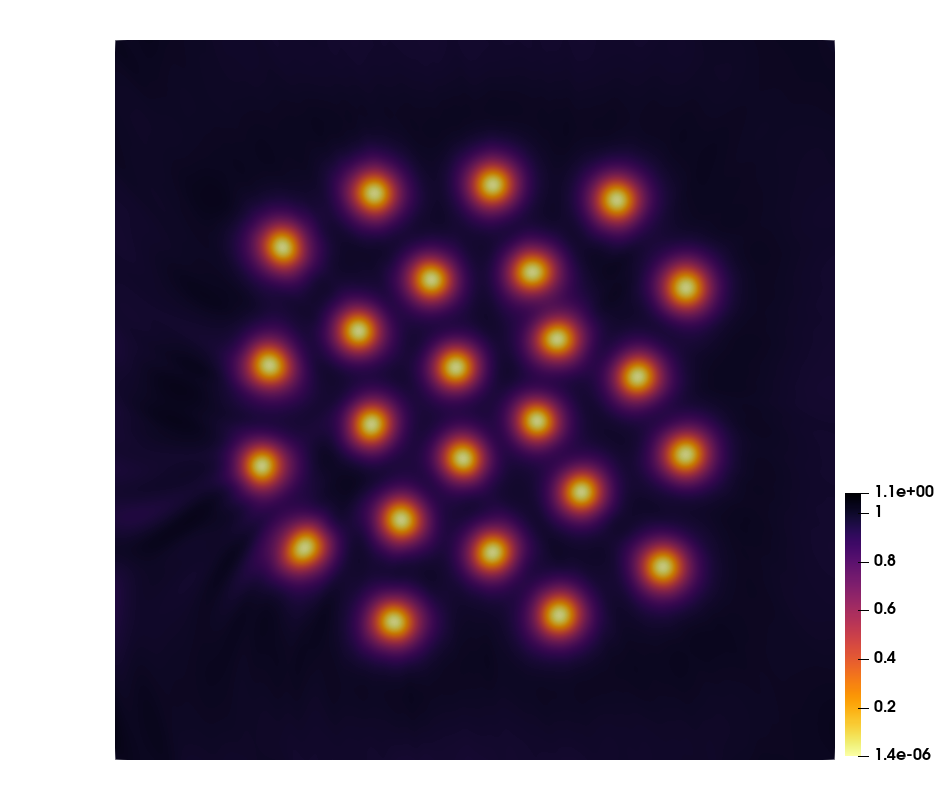}
    \end{minipage}
    \begin{minipage}{0.15\textwidth}
    \includegraphics[scale=0.07]{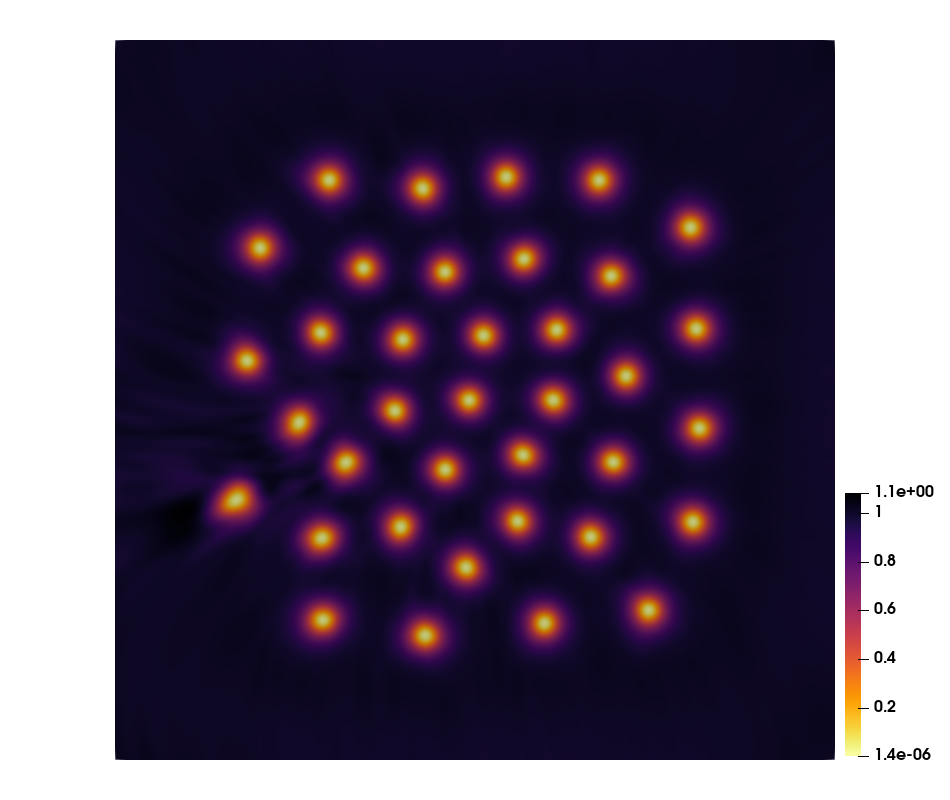}
    \end{minipage}
    \begin{minipage}{0.15\textwidth}
    \includegraphics[scale=0.07]{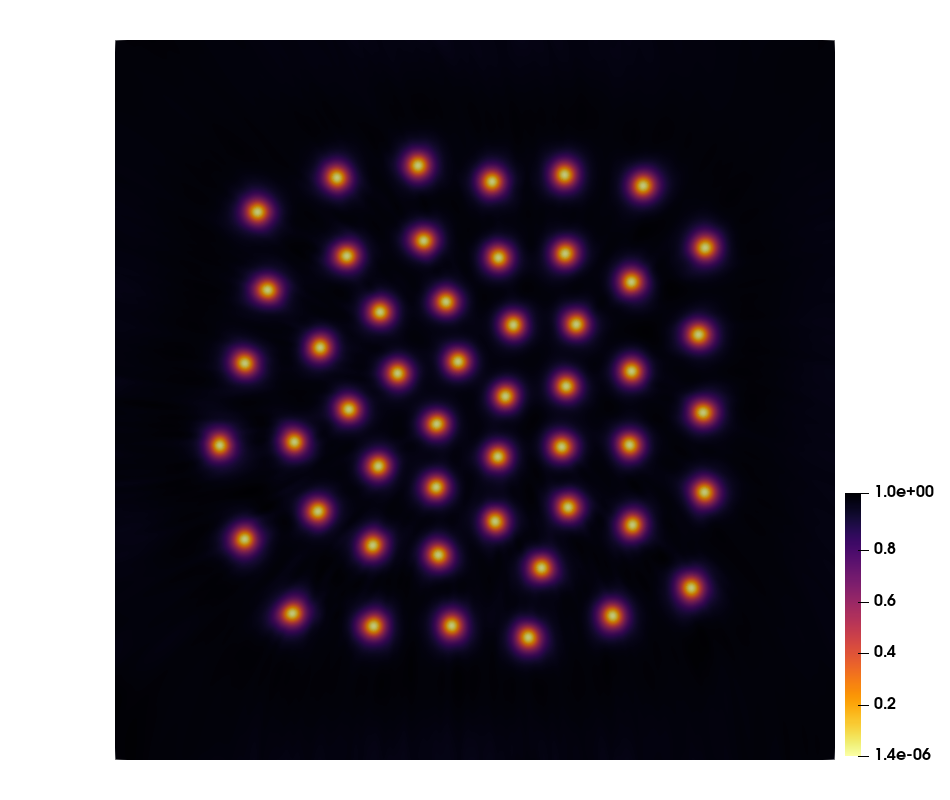}
    \end{minipage} 
    \\
    \begin{minipage}{0.1\textwidth}
    \centering
    \texttt{GL-F2} (NN)
    \end{minipage}
    \begin{minipage}{0.15\textwidth}
    \includegraphics[scale=0.07]{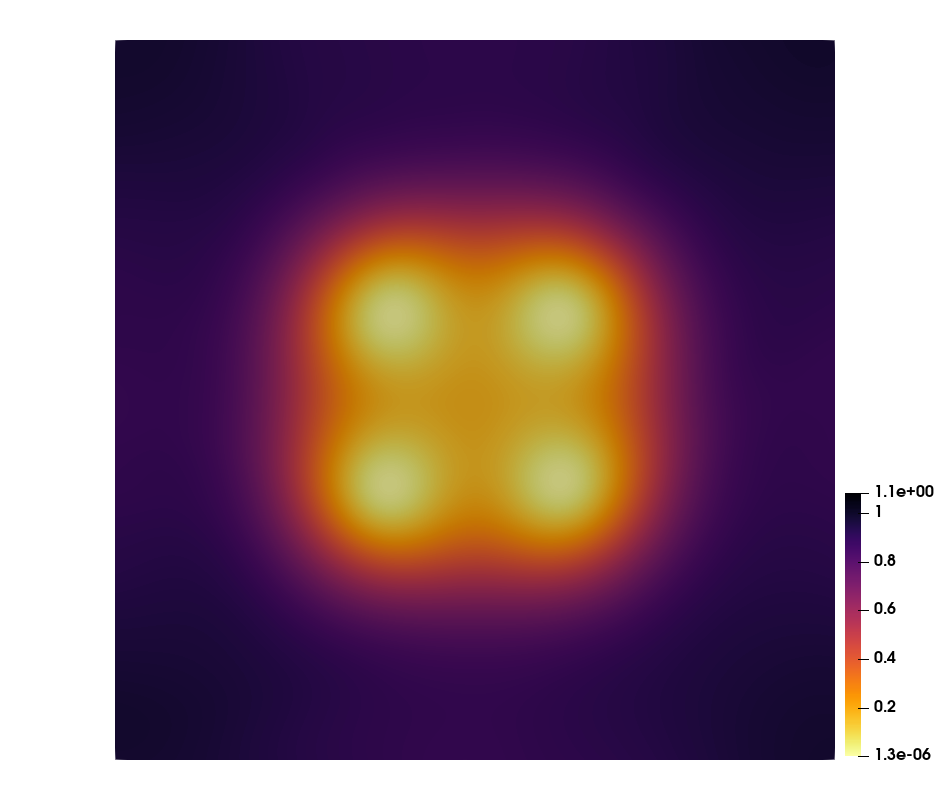}
    \end{minipage}
    \begin{minipage}{0.15\textwidth}
    \includegraphics[scale=0.07]{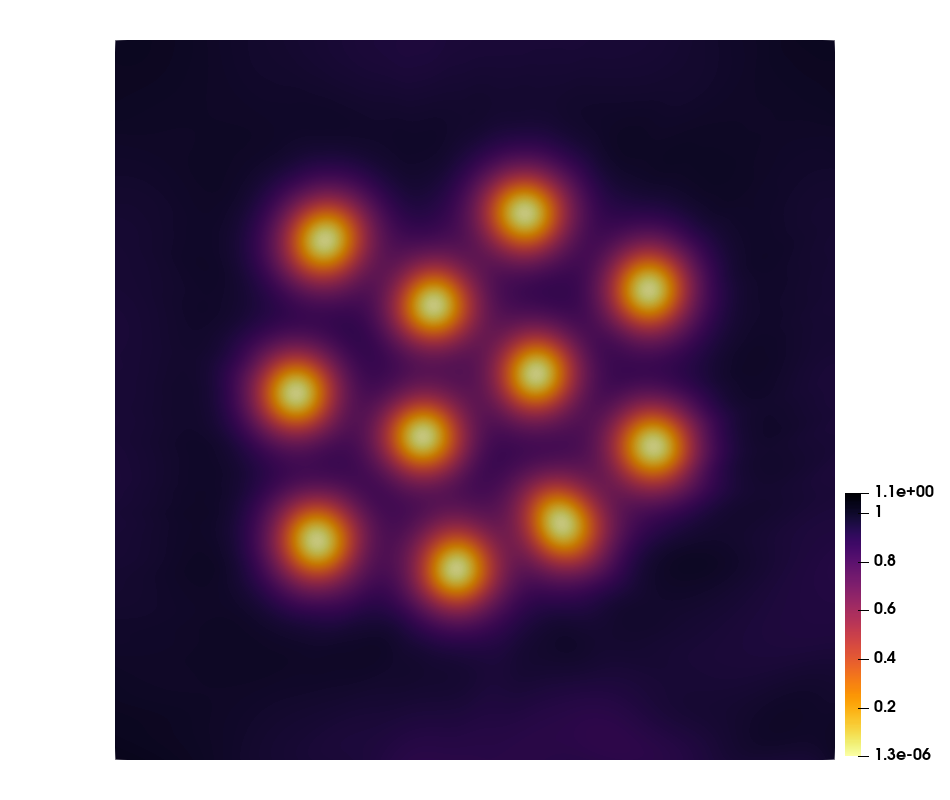}
    \end{minipage}
    \begin{minipage}{0.15\textwidth}
    \includegraphics[scale=0.07]{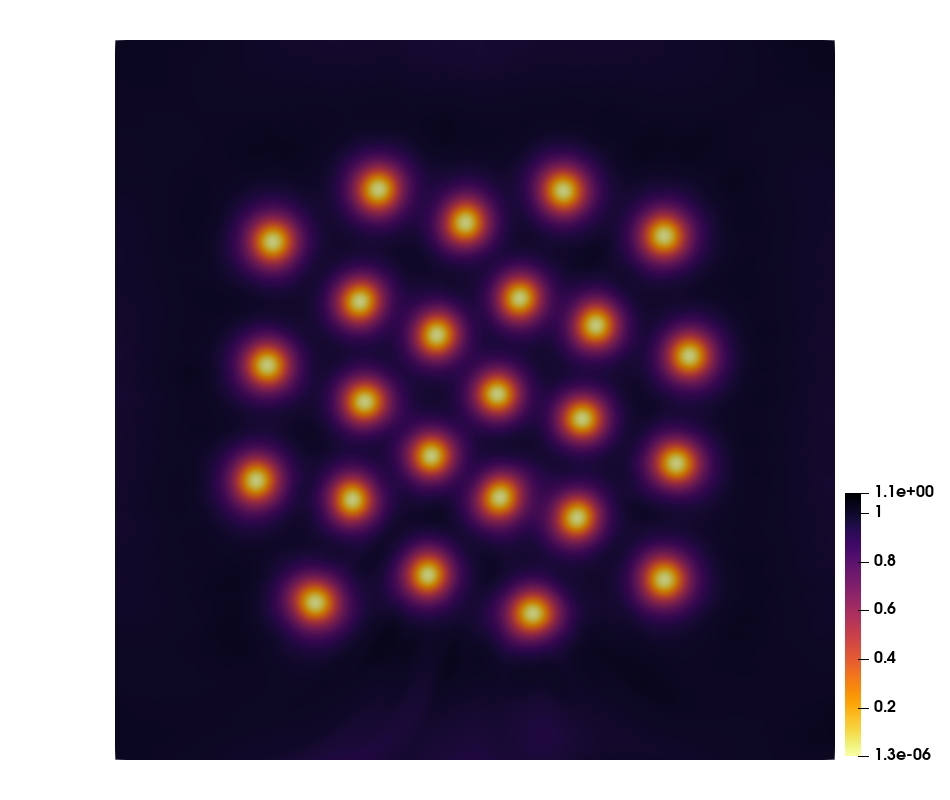}
    \end{minipage}
    \begin{minipage}{0.15\textwidth}
    \includegraphics[scale=0.07]{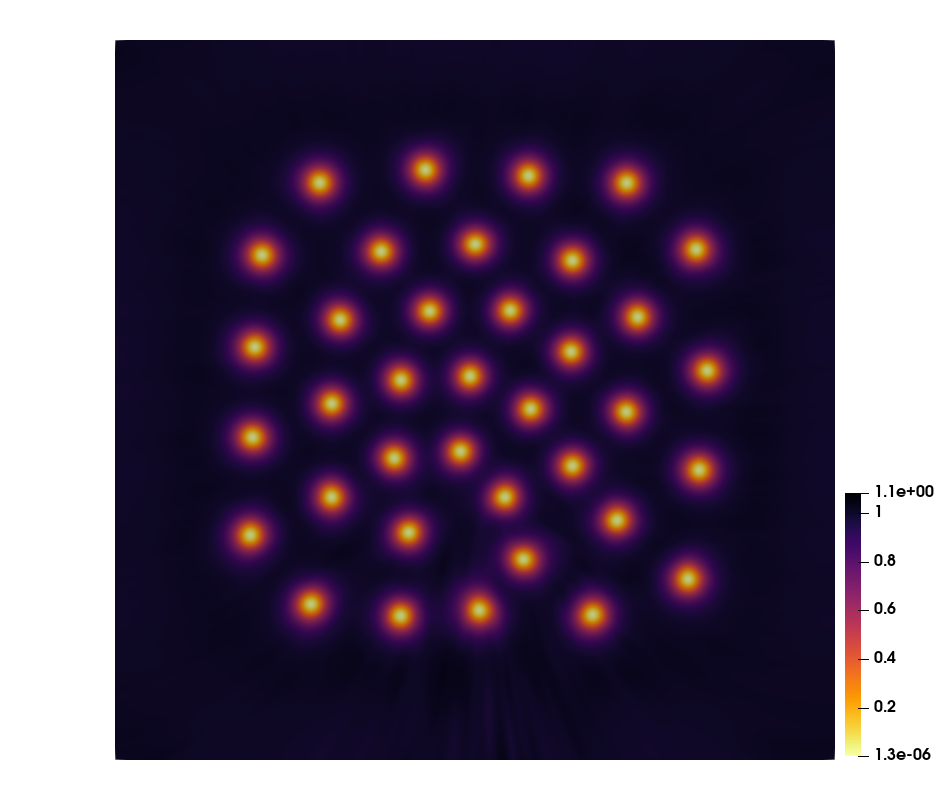}
    \end{minipage}
    \begin{minipage}{0.15\textwidth}
    \includegraphics[scale=0.07]{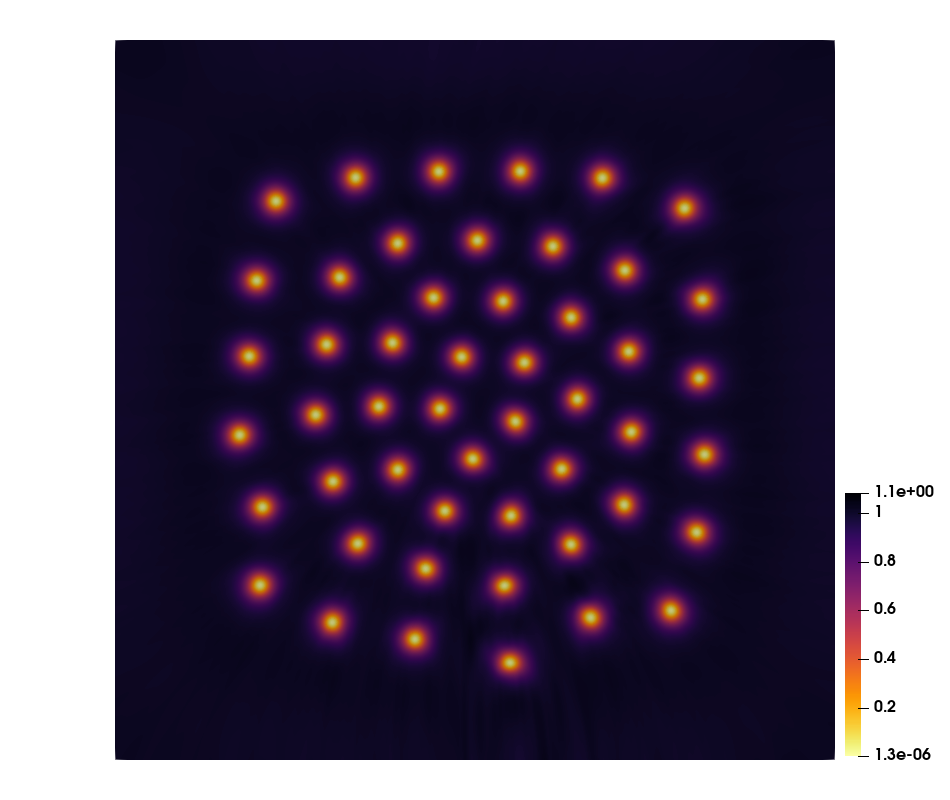}
    \end{minipage} 
    \\
    \begin{minipage}{0.1\textwidth}
    \centering
    \texttt{GL-F1} (hyb.)
    \end{minipage}
    \begin{minipage}{0.15\textwidth}
    \includegraphics[scale=0.07]{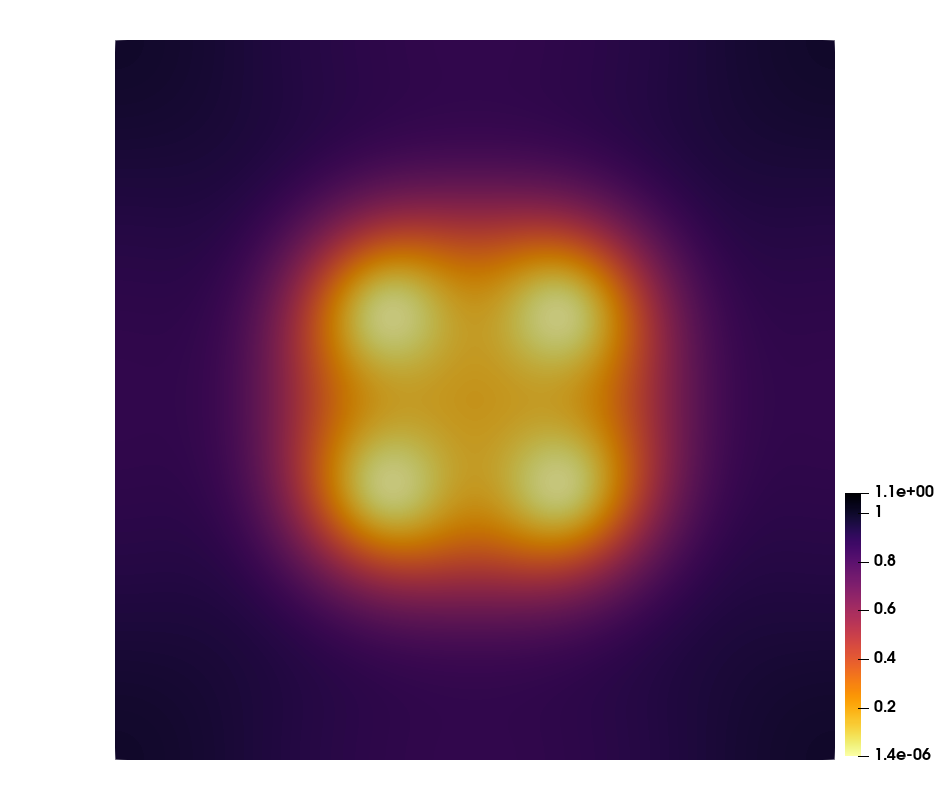}
    \end{minipage}
    \begin{minipage}{0.15\textwidth}
    \includegraphics[scale=0.07]{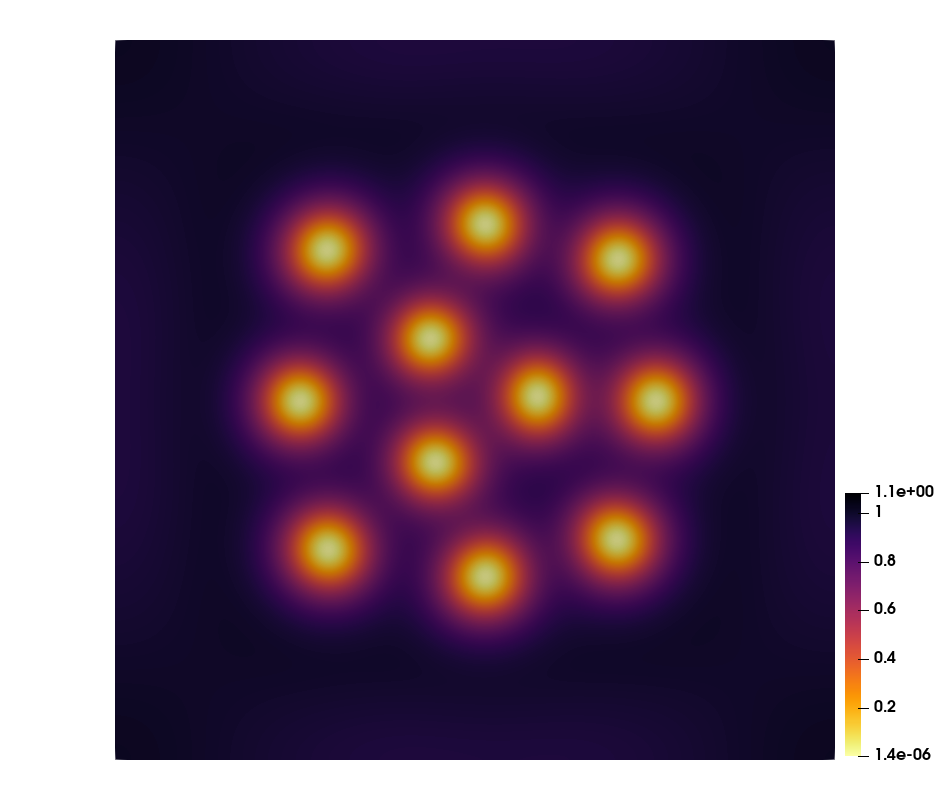}
    \end{minipage}
    \begin{minipage}{0.15\textwidth}
    \includegraphics[scale=0.07]{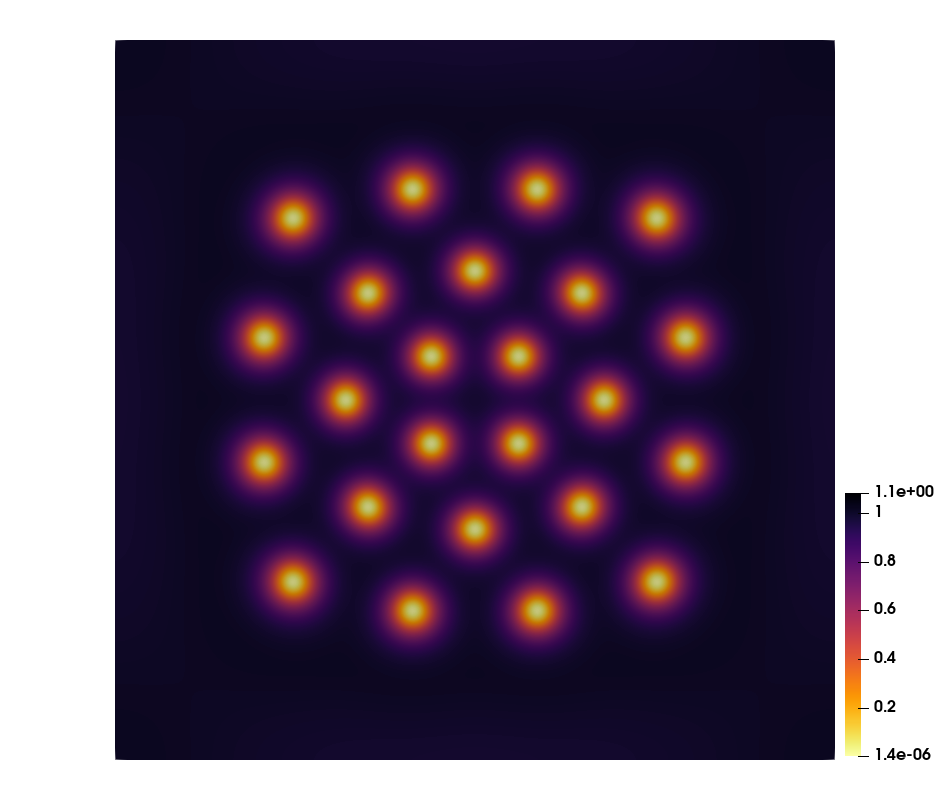}
    \end{minipage}
    \begin{minipage}{0.15\textwidth}
    \includegraphics[scale=0.07]{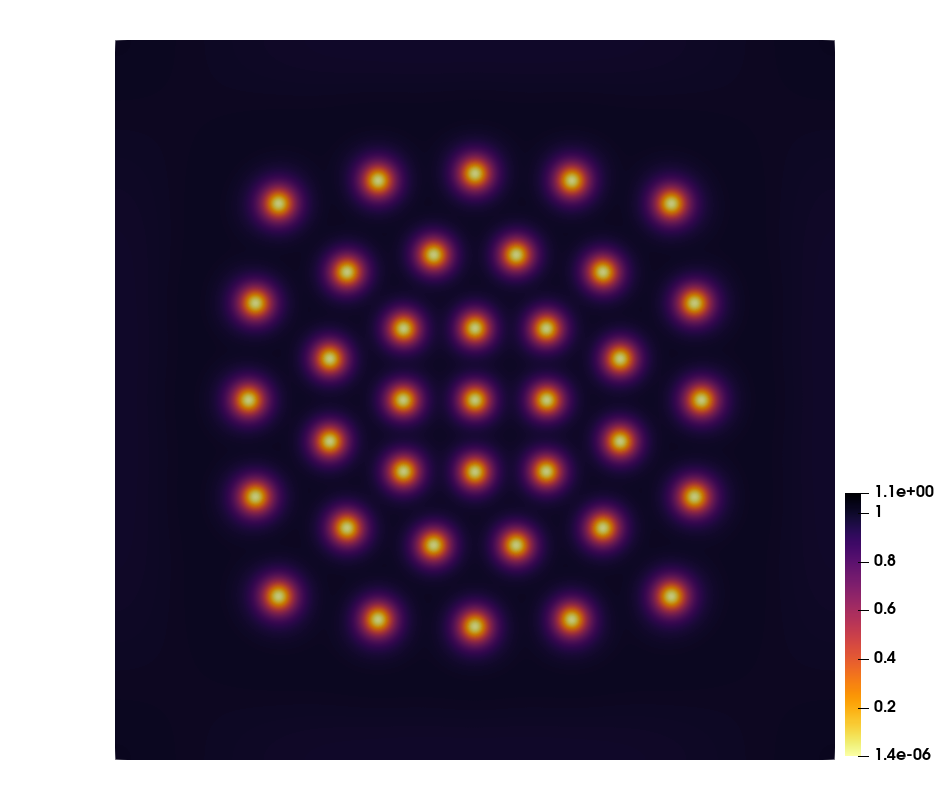}
    \end{minipage}
    \begin{minipage}{0.15\textwidth}
    \includegraphics[scale=0.07]{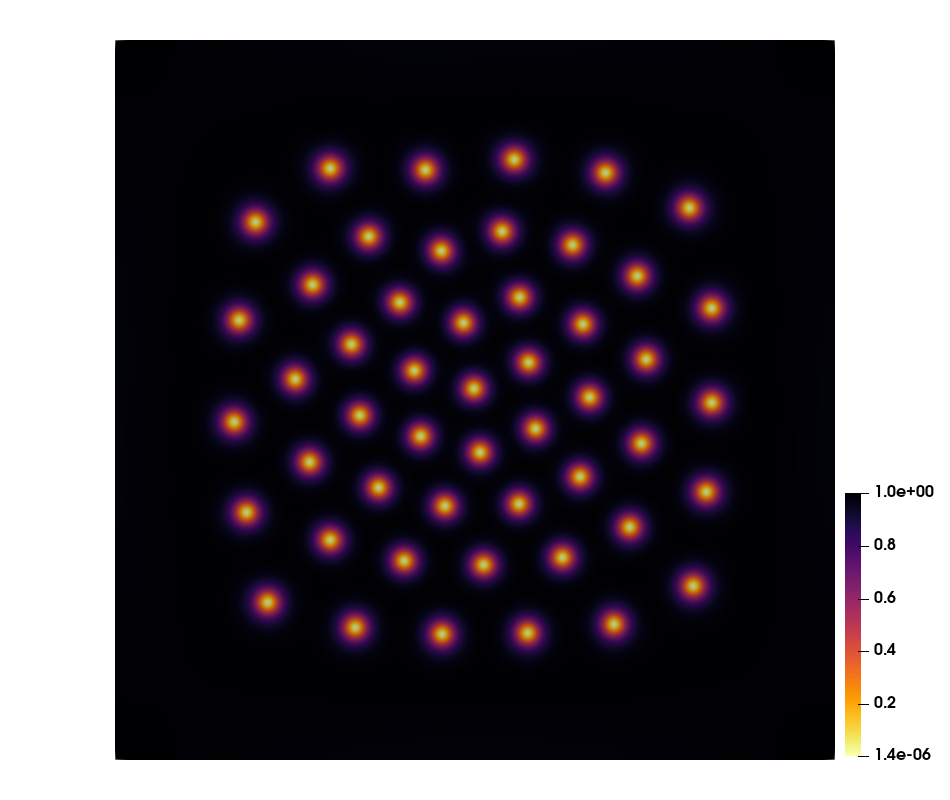}
    \end{minipage} 
    \\
    \begin{minipage}{0.1\textwidth}
    \centering
    \texttt{GL-F2} (hyb.)
    \end{minipage}
    \begin{minipage}{0.15\textwidth}
    \includegraphics[scale=0.07]{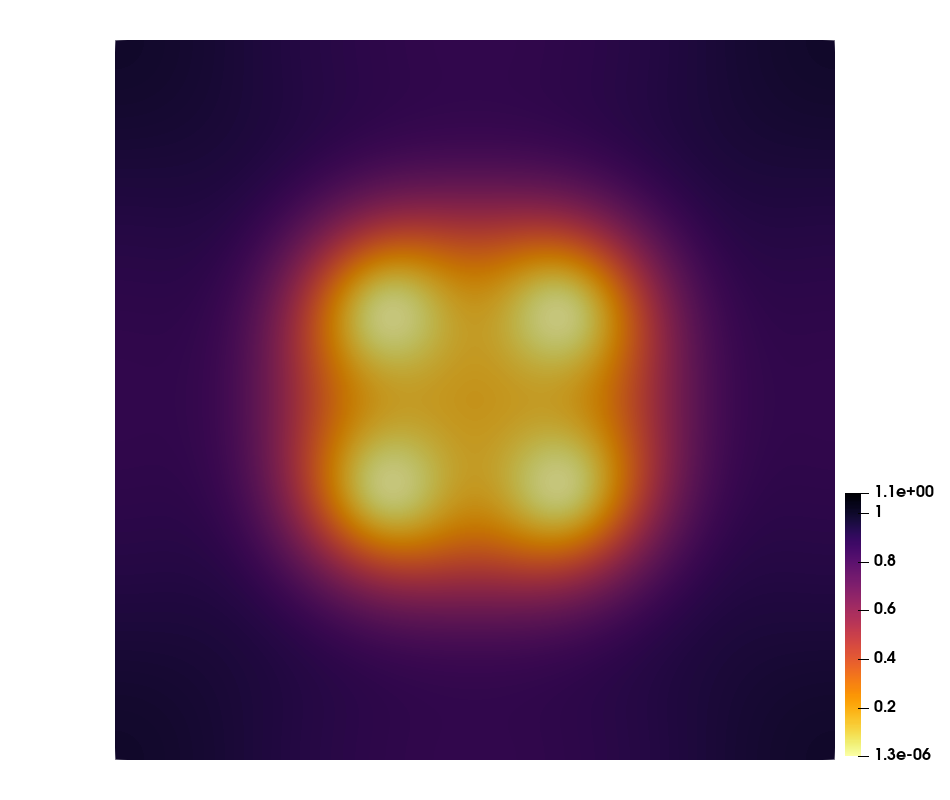}
    \end{minipage}
    \begin{minipage}{0.15\textwidth}
    \includegraphics[scale=0.07]{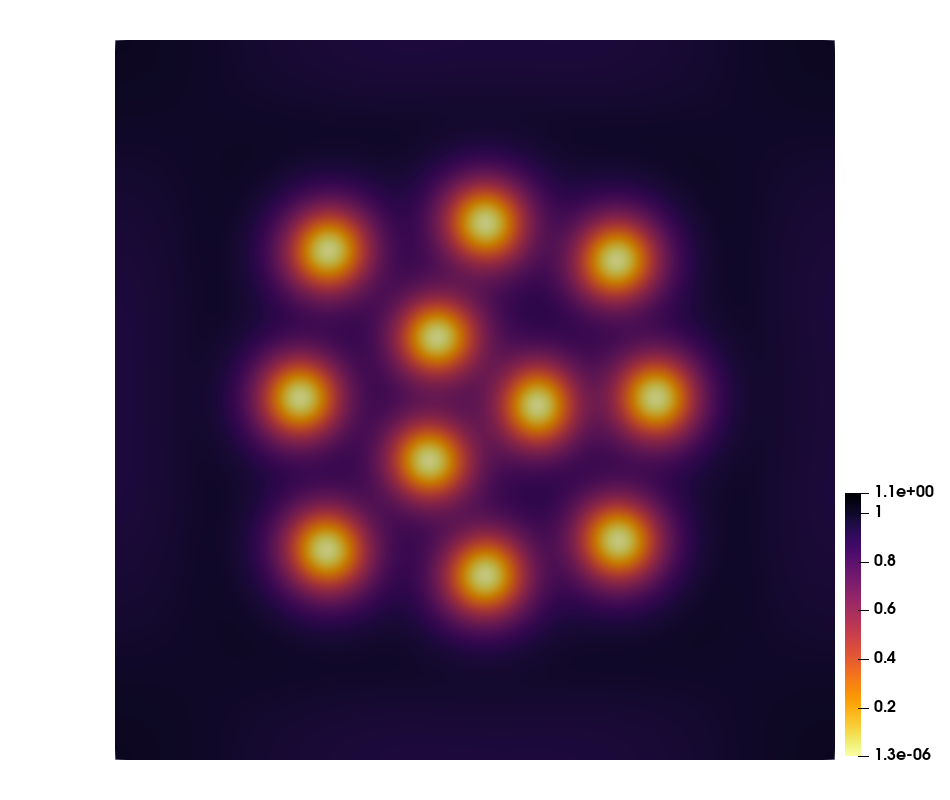}
    \end{minipage}
    \begin{minipage}{0.15\textwidth}
    \includegraphics[scale=0.07]{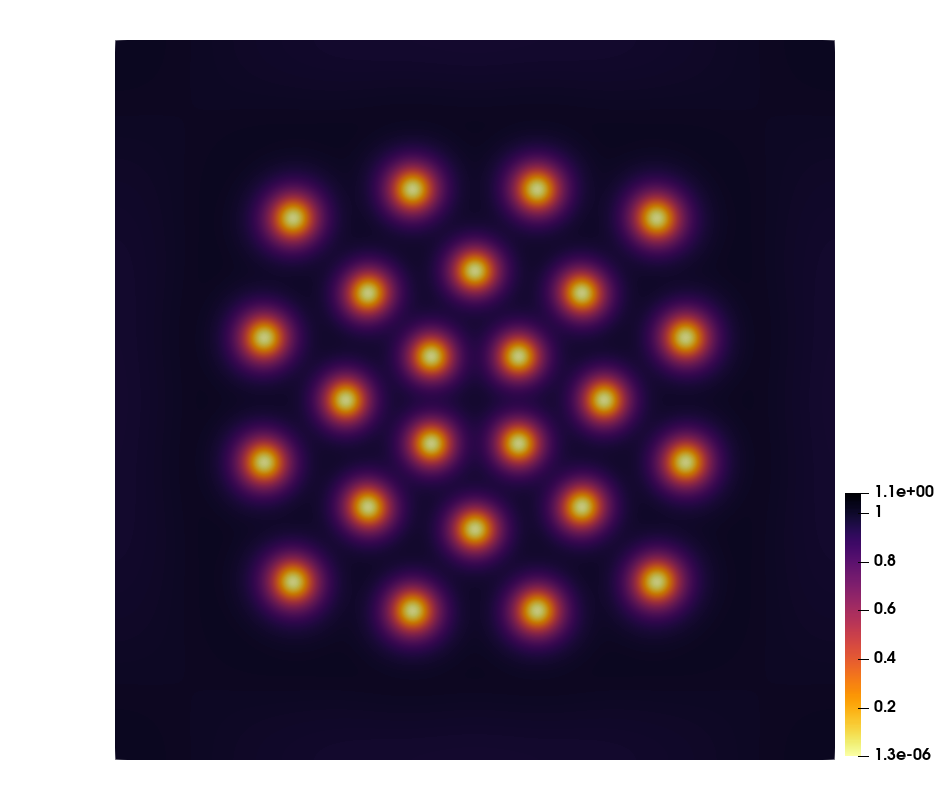}
    \end{minipage}
    \begin{minipage}{0.15\textwidth}
    \includegraphics[scale=0.07]{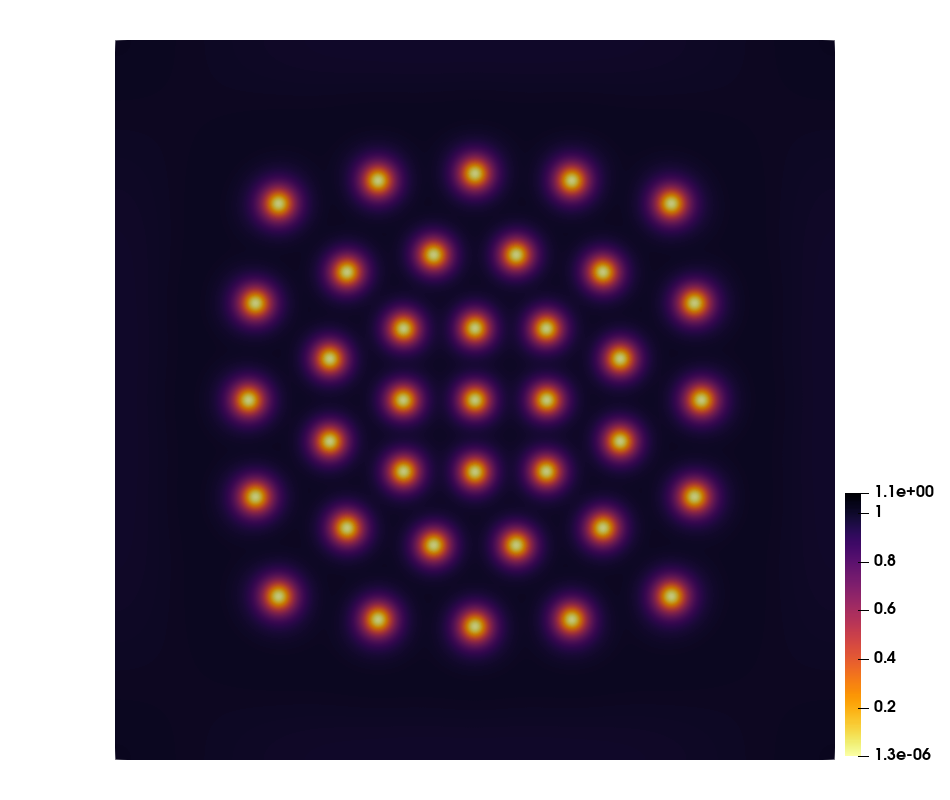}
    \end{minipage}
    \begin{minipage}{0.15\textwidth}
    \includegraphics[scale=0.07]{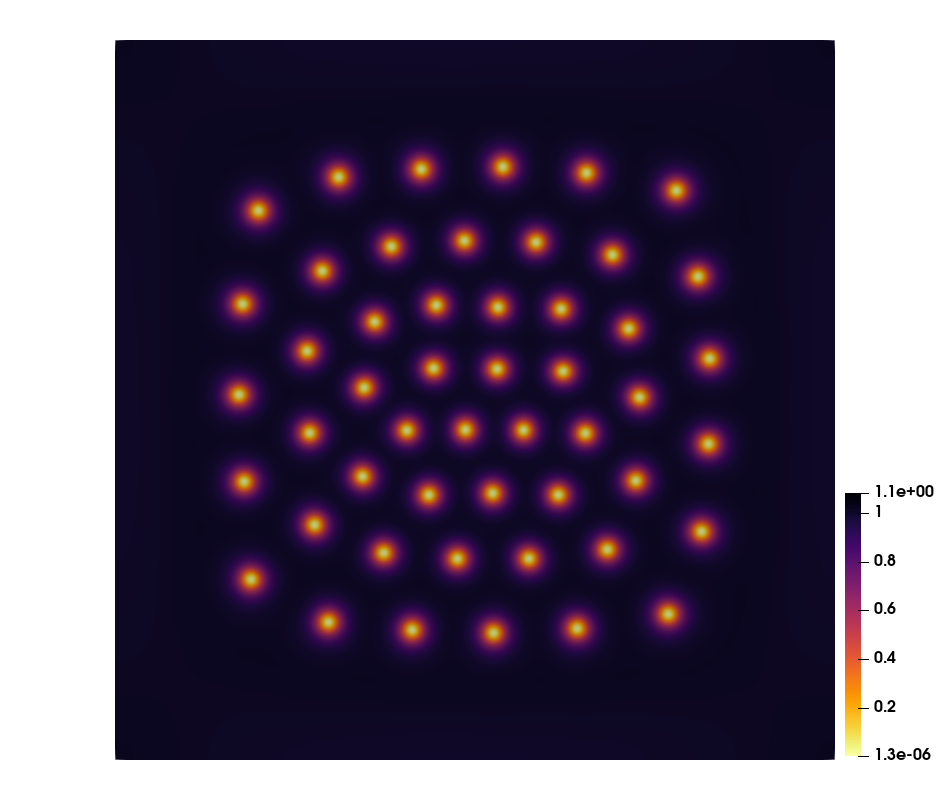}
    \end{minipage}
    \caption{Densities $|u|^2$ of computed minimizers for $\kappa = 10 ,25, 50, 75, 100$ (left to right), corresponding to Table~\ref{tab:energy_full}. Dark indicates values close to 1; light indicates values close to 0. First row: reference minimizer computed with the classical FE iterative solver and best heuristic initial values $\varphi_j$ (best). Second and third row: pure NN models \texttt{GL-F1} (NN) and \texttt{GL-F2} (NN). Fourth and fifth row: hybrid approach \texttt{GL-F1} (hyb.) and \texttt{GL-F2} (hyb.) using initial values from the NN.}
    \label{fig:minimizers_full_square}
\end{figure}

Similar to the previous experiment, the first block of Table~\ref{tab:energy_full} shows the computed minimizers obtained with the classical finite element iterative solver using heuristic initial values $\varphi_j$. The reference minimizers with the lowest energy levels are again highlighted in bold.

For the full model, we trained two NNs, denoted by \texttt{GL-F1} and \texttt{GL-F2}, and identified architectures and training strategies that yield to minimizers of low energy. The corresponding parameters are listed in Table \ref{tab:parameters}, and the resulting energy levels of the computed GL minimizers are reported in the second block of Table \ref{tab:energy_full}. The training times of \texttt{GL-F1} and \texttt{GL-F2} were~7 and~4.5 hours, respectively, on the HoreKa machine at KIT. As in the reduced-model setting, the pure NN approaches produce energy levels close to, but usually slightly above, those of the reference minimizers. The only exception is the case $\kappa = 75$, where the networks find a minimizer with significantly lower energy, suggesting already a candidate for the global GL minimizer.

We then apply the hybrid approach to the full problem to refine the NN predictions. The corresponding energy levels for the hybrid approach using \texttt{GL-F1} and \texttt{GL-F2} are shown in the third block of Table~\ref{tab:energy_full}. For $\kappa = 10$ and $50$, the hybrid approach recovers the same energy states as the best classical methods, while for $\kappa = 25, 75,$ and $100$ it yields states with significantly lower energy. Hence, using NN predictions as initial values can improve the iterative solver and lead to lower-energy states.

Figure \ref{fig:minimizers_full_square} shows plots of the densities $|u|^2$ of the computed minimizers. The observations are very similar to the case of the reduced model. The NNs yield rough approximations of the vortex pattern in the minimizers that, however, need to be refined by the iterative solver, i.e., the hybrid approach. In particular, considering the cases $\kappa = 25, 75, 100$ the hybrid approach compute a significant different vortex pattern of lower energy. 

\subsection{Fast-training neural networks} 

Finally, we consider a scenario in which, for a given GL model, the goal is to train a NN on demand, with a limited training time of about~1.5 hours. As a basis for the NNs, we use the same architecture and hyperparameters as for \texttt{GL-F1} and train five different NNs by varying the initial seeds in the training process. We then exploit the hybrid approach to the given GL test model. This experiment evaluates the reliability and quality of the approximations achievable with the \texttt{GL-F1} architecture and hyperparameters. Note that we do not use the trained model \texttt{GL-F1}, but instead train the same architecture from scratch with only 4 epochs. The obtained networks are termed \texttt{GL-U*}, $\texttt{*} \in \{1,2,3,4,5\}$.

\begin{table}[h!]
\centering
\small
\begin{tabular}{c c c c c c}
\toprule
$\kappa$ & 10 & 25 & 50 & 75 & 100 \\
\midrule
$\varphi_j$ (best) & 0.10440628 & 0.06491247 & 0.04148069 & 0.03446474 & 0.02615202 \\
\midrule
\makebox[18mm][l]{\texttt{GL-U1} (hyb.) }
& 0.10440628 &  0.06445609 & 0.04148069 & 0.03176793 & 0.02615239 \\
\makebox[18mm][l]{\texttt{GL-U2} (hyb.) }
& 0.10440628 &  \textbf{0.06445608} & 0.04160989 & \textbf{0.03171263} & 0.02628949 \\
\makebox[18mm][l]{\texttt{GL-U3} (hyb.) }
& 0.10440628 &  0.06445617 & 0.04160988 & 0.03171263 & \textbf{0.02614767} \\
\makebox[18mm][l]{\texttt{GL-U4} (hyb.) }
& \textbf{0.10440628} &  0.06445608 & \textbf{0.04148069} & 0.03171263 & 0.02628949 \\
\makebox[18mm][l]{\texttt{GL-U5} (hyb.) }
& 0.10440628 &  0.06445615 & 0.04148069 & 0.03171263 & 0.02660915 \\
\bottomrule
\end{tabular}
\caption{Energy levels $E(u,\MagPot)$ of computed energy minimizers for the full GL model on the unit square domain. First block: classical FE iterative solver with heuristic initial values $\varphi_j$ for $u$ and $0$ for $\MagPot$. Second block: hybrid approach using the fast-trained 
\texttt{GL-U*}
(hyb.), $\texttt{*} \in \{ 1, 2, 3, 4, 5\}$.}
\label{tab:unit_square}
\end{table}

\begin{figure}[h!]
    \flushleft
    \begin{minipage}{0.1\textwidth}
    \centering
    $\varphi_j$ (best)
    \end{minipage}
    \begin{minipage}{0.15\textwidth}
    \includegraphics[scale=0.07]{full_with_iv1_kappa10__plot_u_abs_init_11_1_1_kappa_10_h_0.00390625_tol_1e-12.png}
    \end{minipage}
    \begin{minipage}{0.15\textwidth}
    \includegraphics[scale=0.07]{full_with_iv4_kappa25__plot_u_abs_init_14_1_1_kappa_25_h_0.00390625_tol_1e-12.png}
    \end{minipage}
    \begin{minipage}{0.15\textwidth}
    \includegraphics[scale=0.07]{full_with_iv1_kappa50__plot_u_abs_init_11_1_1_kappa_50_h_0.00390625_tol_1e-12.png}
    \end{minipage}
    \begin{minipage}{0.15\textwidth}
    \includegraphics[scale=0.07]{full_with_iv4_kappa75__plot_u_abs_init_14_1_1_kappa_75_h_0.00390625_tol_1e-12.png}
    \end{minipage}
    \begin{minipage}{0.15\textwidth}
    \includegraphics[scale=0.07]{full_with_iv1_kappa100__plot_u_abs_init_11_1_1_kappa_100_h_0.00390625_tol_1e-12.png}
    \end{minipage} 
    \\
    \begin{minipage}{0.1\textwidth}
    \centering
    best \texttt{GL-U*} (hyb.)
    \end{minipage}
    \begin{minipage}{0.15\textwidth}
    \includegraphics[scale=0.07]{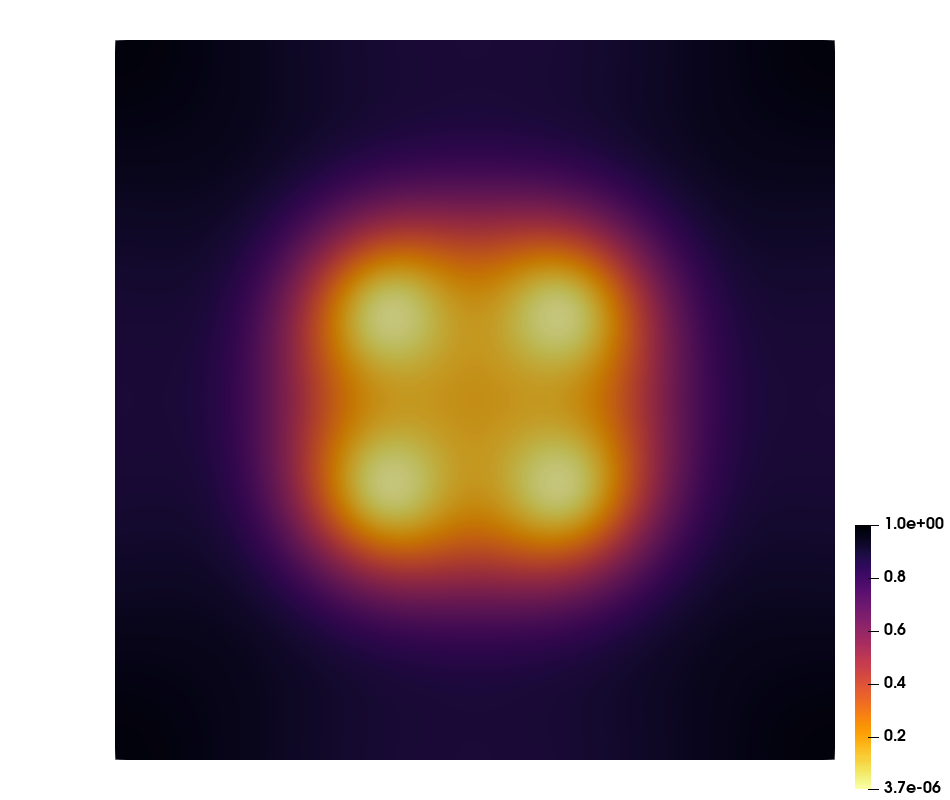}
    \end{minipage}
    \begin{minipage}{0.15\textwidth}
    \includegraphics[scale=0.07]{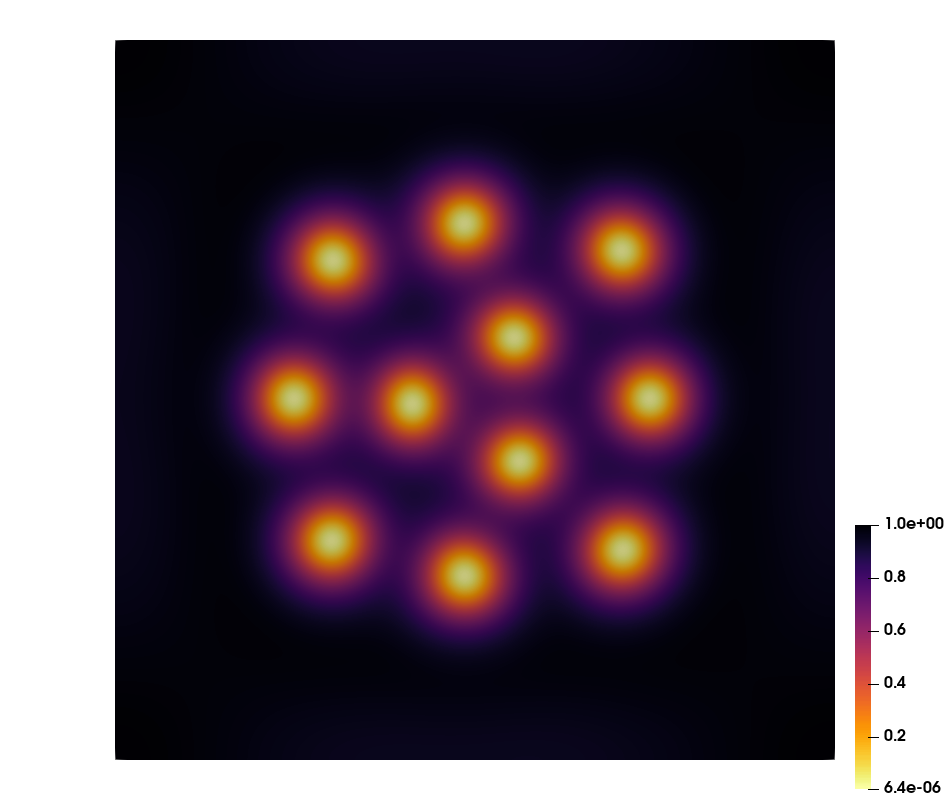}
    \end{minipage}
    \begin{minipage}{0.15\textwidth}
    \includegraphics[scale=0.07]{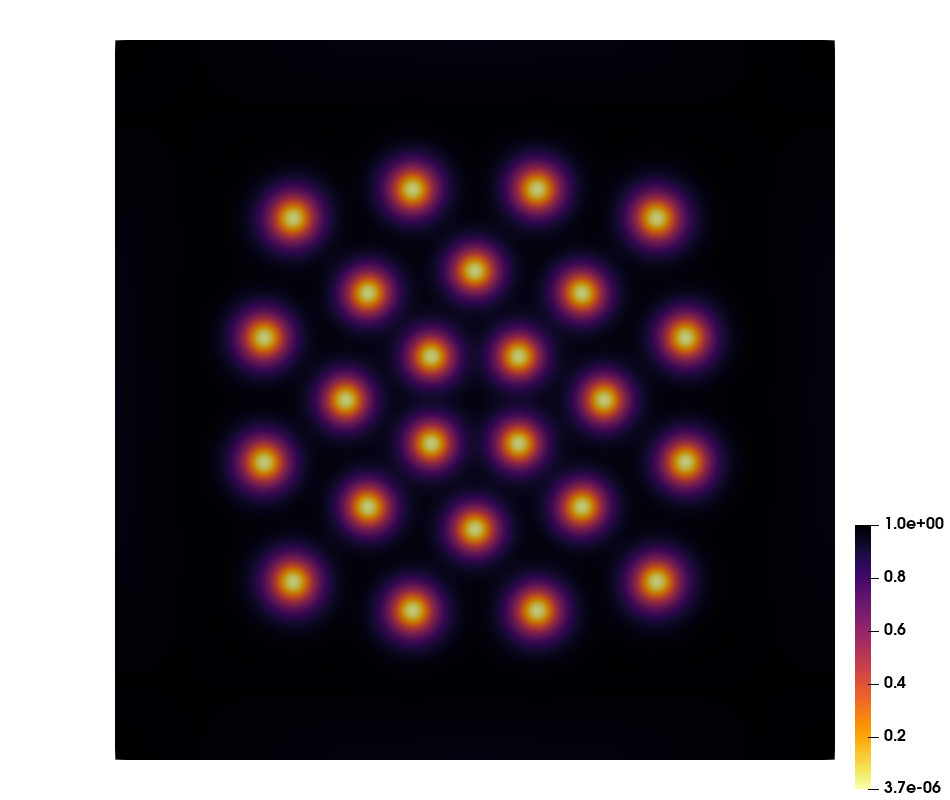}
    \end{minipage}
    \begin{minipage}{0.15\textwidth}
    \includegraphics[scale=0.07]{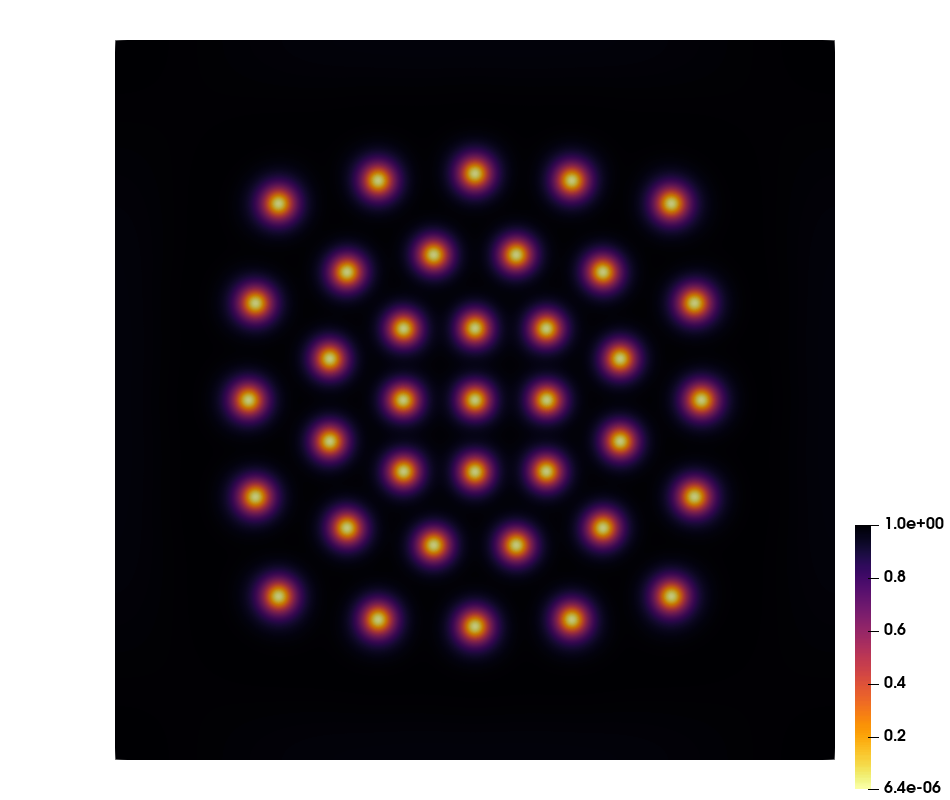}
    \end{minipage}
    \begin{minipage}{0.15\textwidth}
    \includegraphics[scale=0.07]{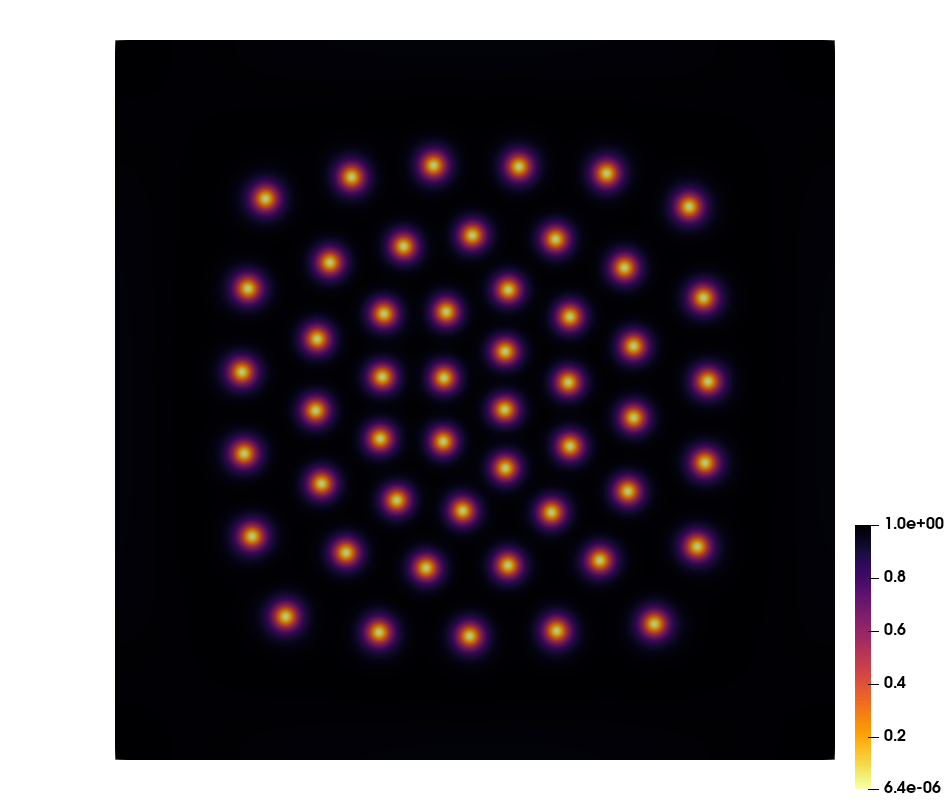}
    \end{minipage} 
    \caption{Densities $|u|^2$ of computed minimizers for $\kappa = 10 ,25, 50, 75, 100$ (left to right), corresponding to Table~\ref{tab:unit_square}. Dark indicates values close to 1; light indicates values close to 0. First row: reference minimizer computed with the classical FE iterative solver and best heuristic initial values $\varphi_j$ (best). Second row: best hybrid approach using the fast-trained  \texttt{GL-U*} (hyb.), $\texttt{*} \in \{ 1, 2, 3, 4, 5\}$.}
    \label{fig:minimizers_unit_square}
\end{figure}

For the first test problem, we consider again the unit square and the external magnetic field $\MagF$ from \eqref{eq:hext}. The energy levels of the resulting GL minimizers computed with the hybrid approach with the five different NNs are shown in Table~\ref{tab:unit_square}. Again, the minimizers of the lowest energy are highlighted in bold. If there are multiple minimizers with the same energy, we highlight the one that requires the fewest iterations. We observe that there is only a minor variation in the energy levels using the hybrid method with \texttt{GL-U*}, trained with different seeds. Compared to the classical method with the heuristic initial values, the obtained energy levels are in favor of the NNs (see, e.g., the lower energy levels for $\kappa = 25$ or $\kappa = 75$). This can be recognized also by considering the densities of the obtained minimizers that are shown in Figure~\ref{fig:minimizers_unit_square}. The different vortex patterns confirm that indeed a different minimizing state is obtained by the hybrid method.

\begin{table}[h!]
\centering
\small
\begin{tabular}{c c c c c c}
\toprule
$\kappa$ & 10 & 25 & 50 & 75 & 100 \\
\midrule
$\varphi_j$ (best) & 0.06958856 & 0.04353387 & 0.02864807 & 0.02261634 & 0.01901139 \\
\midrule
\makebox[18mm][l]{\texttt{GL-L1} (hyb.) }
& 0.06958856 & 0.04310817 & 0.02856779 & 0.02204774 & 0.01830038 \\
\makebox[18mm][l]{\texttt{GL-L2} (hyb.) }
& \textbf{0.06958856} & 0.04310817 & 0.02854874 & \textbf{0.02204774} & 0.01830132 \\
\makebox[18mm][l]{\texttt{GL-L3} (hyb.) }
& 0.06958856 & \textbf{0.04310817} & 0.02862300 & 0.02204774 & 0.01830037 \\
\makebox[18mm][l]{\texttt{GL-L4} (hyb.) }
& 0.06958856 & 0.04310817 & 0.02862300 & 0.02204485 & \textbf{0.01829811} \\
\makebox[18mm][l]{\texttt{GL-L5} (hyb.) }
& 0.06958856 & 0.04310817 & \textbf{0.02854874} & 0.02204774 & 0.01830038 \\
\bottomrule
\end{tabular}
\caption{Energy levels $E(u,\MagPot)$ of computed energy minimizers for the full GL model on the L-shaped domain. First block: classical FE iterative solver with heuristic initial values $\varphi_j$ for $u$ and $0$ for $\MagPot$. Second block: hybrid approach using the fast-trained \texttt{GL-L*} (hyb.), $\texttt{*} \in \{ 1, 2, 3, 4, 5\}$.}
\label{tab:Lshape}
\end{table}

\begin{figure}[h!]
    \flushleft
    \begin{minipage}{0.1\textwidth}
    \centering
    $\varphi_j$ (best)
    \end{minipage}
    \begin{minipage}{0.15\textwidth}
    \includegraphics[scale=0.07]{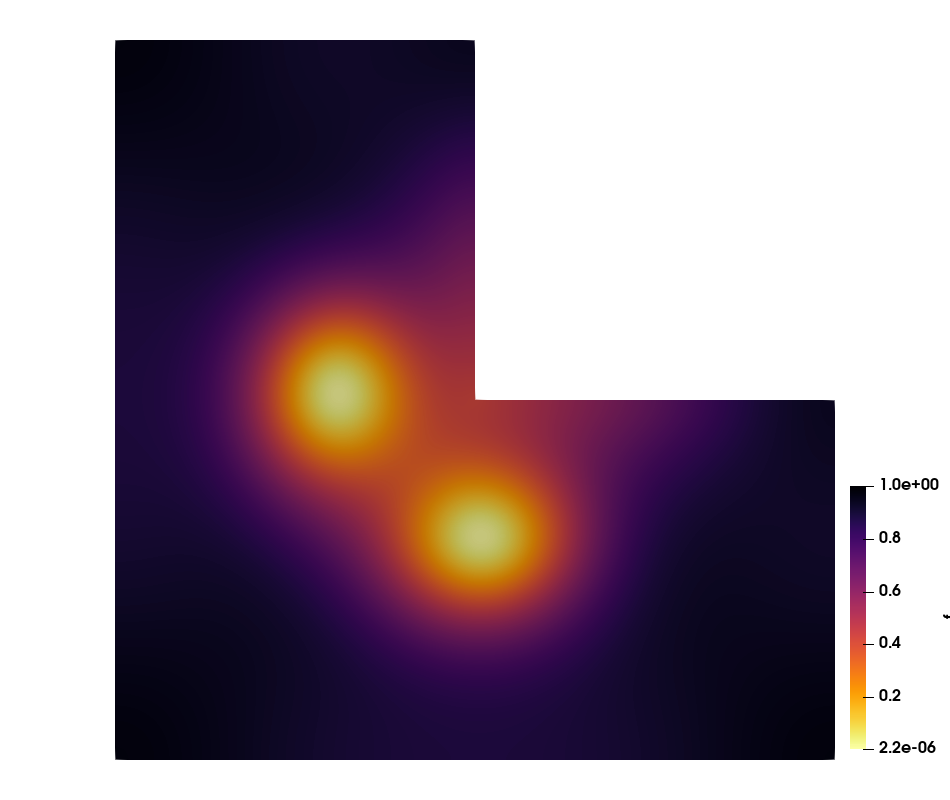}
    \end{minipage}
    \begin{minipage}{0.15\textwidth}
    \includegraphics[scale=0.07]{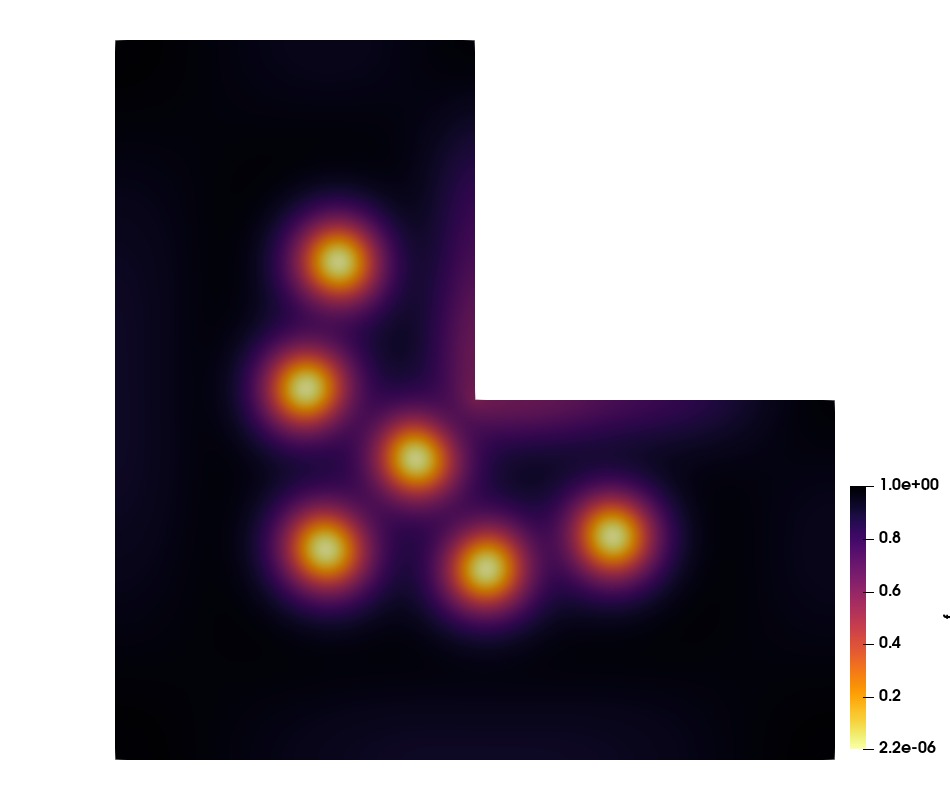}
    \end{minipage}
    \begin{minipage}{0.15\textwidth}
    \includegraphics[scale=0.07]{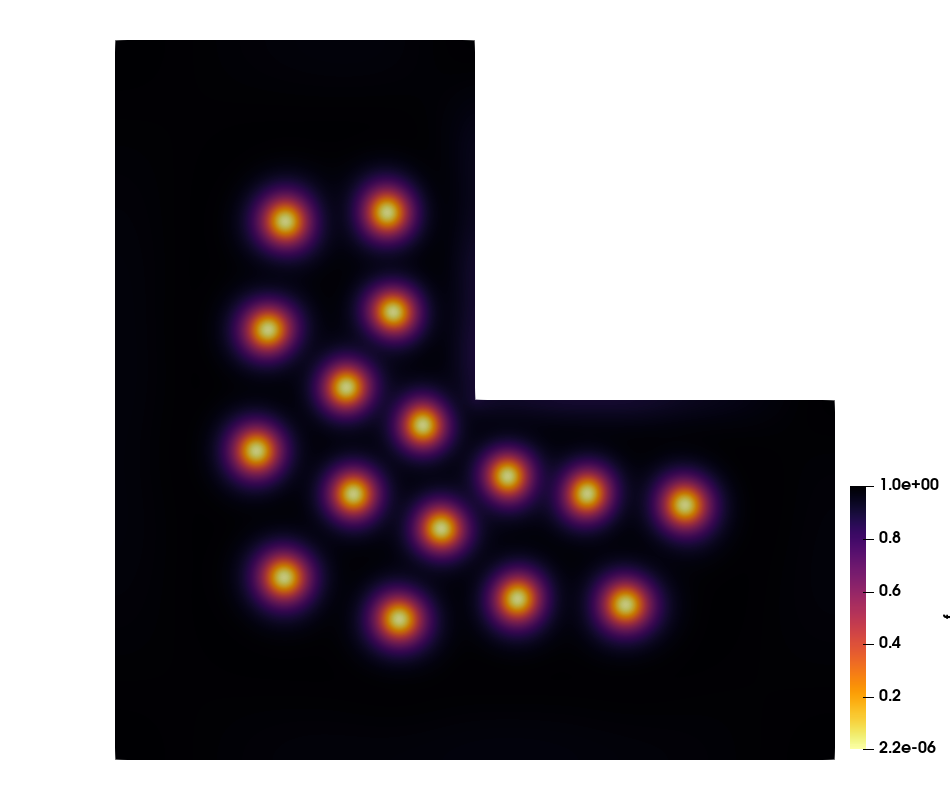}
    \end{minipage}
    \begin{minipage}{0.15\textwidth}
    \includegraphics[scale=0.07]{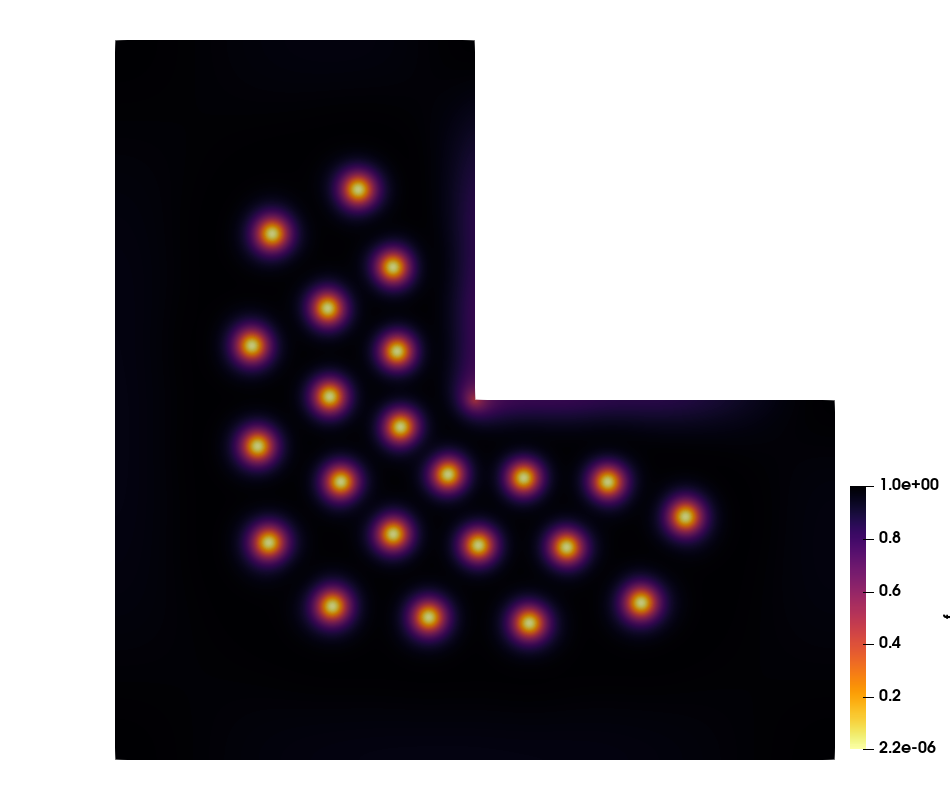}
    \end{minipage}
    \begin{minipage}{0.15\textwidth}
    \includegraphics[scale=0.07]{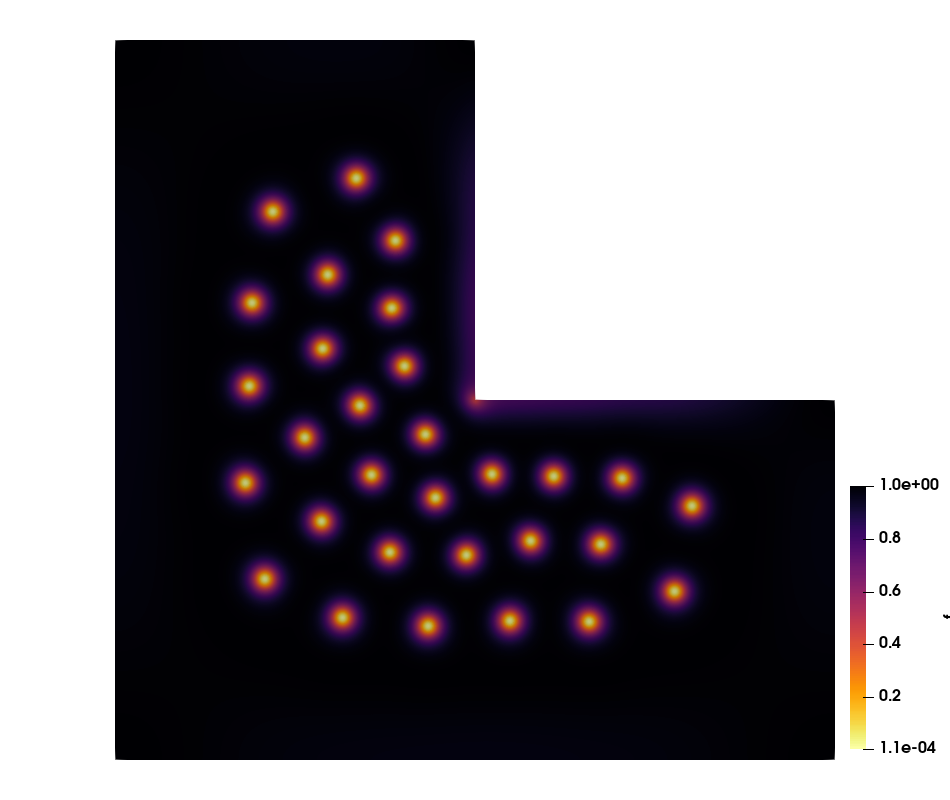}
    \end{minipage} 
    \\
    \begin{minipage}{0.1\textwidth}
    \centering
    best \texttt{GL-L*} (hyb.)
    \end{minipage}
    \begin{minipage}{0.15\textwidth}
    \includegraphics[scale=0.07]{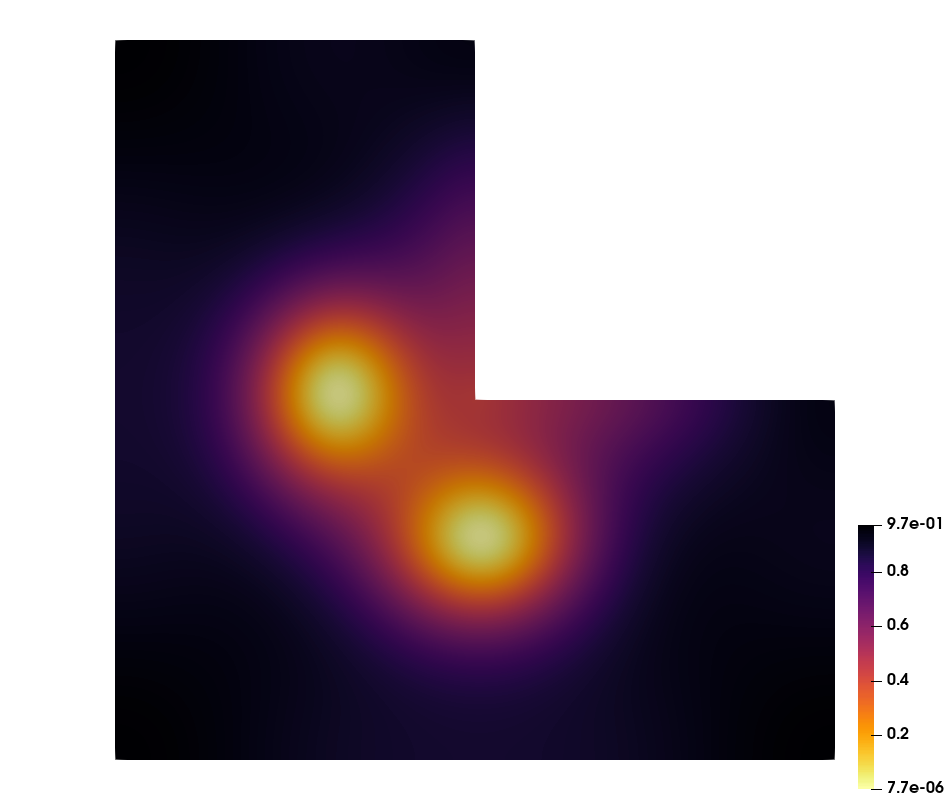}
    \end{minipage}
    \begin{minipage}{0.15\textwidth}
    \includegraphics[scale=0.07]{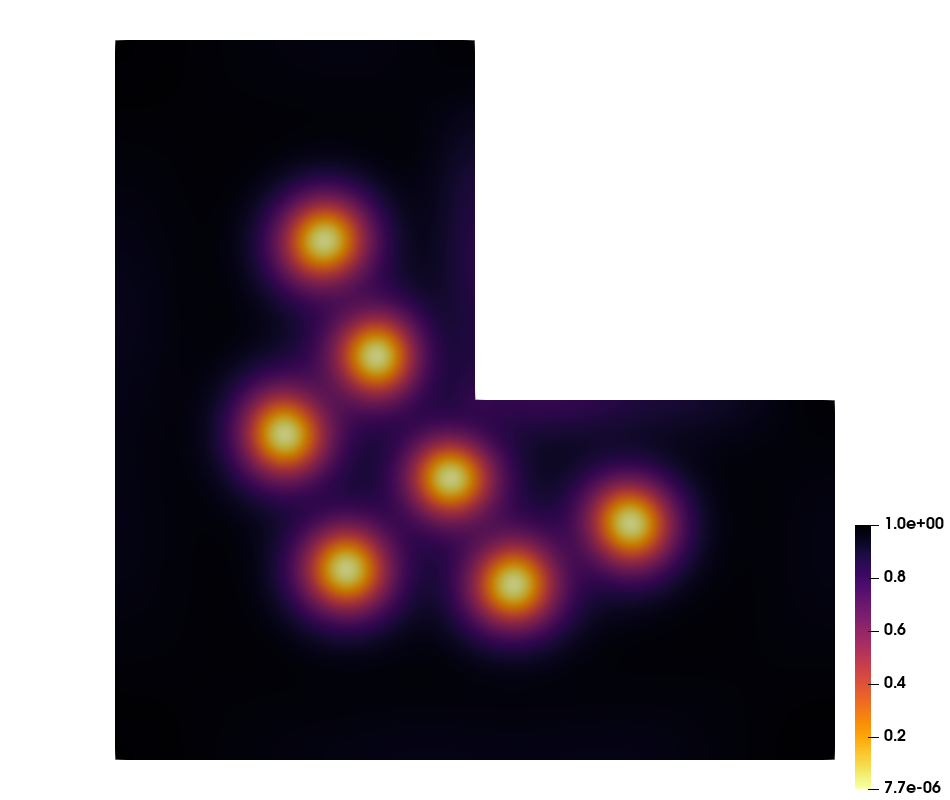}
    \end{minipage}
    \begin{minipage}{0.15\textwidth}
    \includegraphics[scale=0.07]{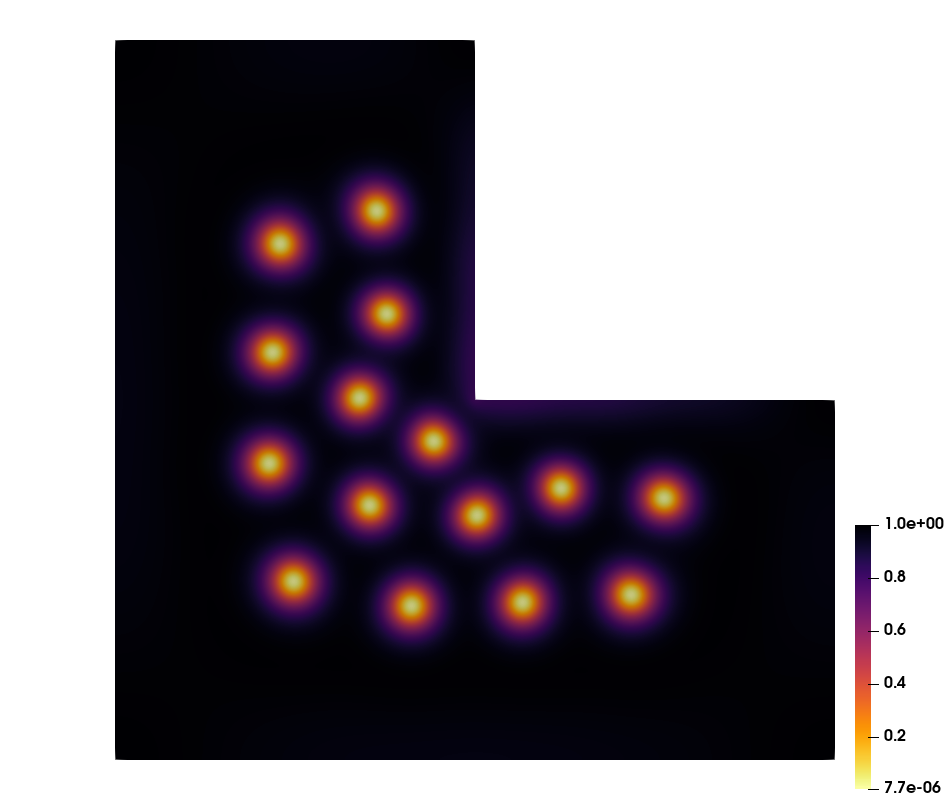}
    \end{minipage}
    \begin{minipage}{0.15\textwidth}
    \includegraphics[scale=0.07]{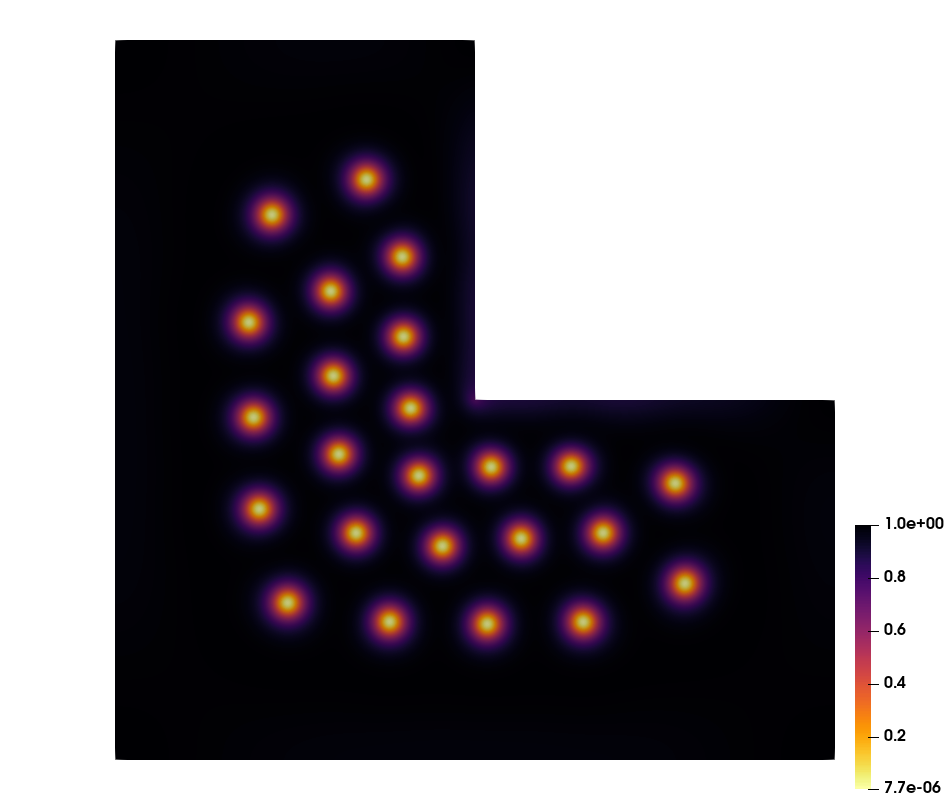}
    \end{minipage}
    \begin{minipage}{0.15\textwidth}
    \includegraphics[scale=0.07]{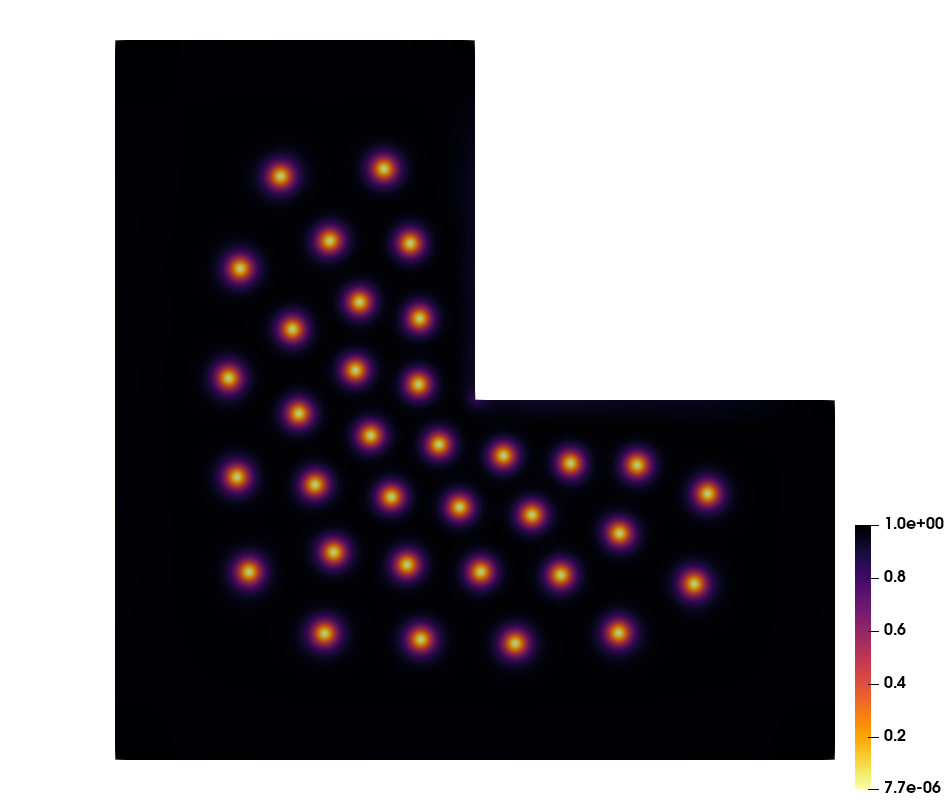}
    \end{minipage} 
    \caption{Densities $|u|^2$ of computed minimizers for $\kappa = 10 ,25, 50, 75, 100$ (left to right), corresponding to Table~\ref{tab:Lshape}. Dark indicates values close to 1; light indicates values close to 0. First row: reference minimizer computed with the classical FE iterative solver and best heuristic initial values $\varphi_j$ (best). Second row: best hybrid approach using the fast-trained  \texttt{GL-L*} (hyb.), $\texttt{*} \in \{ 1, 2, 3, 4, 5\}$.}
    \label{fig:minimizers_Lshape}
\end{figure}

For the second test problem, we replace the computational domain by an L-shaped domain, while keeping the external magnetic field $\MagF$ from \eqref{eq:hext}. We employ the architecture and hyperparameters of \texttt{GL-F1}, but as before train the model from scratch, only adjusting the learning rate and its scheduling. We restrict the modifications of the model to these hyperparameters, since they are both essential for successful training and can be identified through a small number of test runs. In practice, tuning the learning rate and its schedule is typically necessary for successful and effective NN training. The obtained networks are termed \texttt{GL-L*}, $* \in \{1,2,3,4,5\}$ and its architecture and hyperparameters can be found in Table~\ref{tab:parameters}.

The energy levels of the computed minimizers by the hybrid approach for the five different models are shown in Table \ref{tab:Lshape}. For comparison, we again employ the classical FE iterative solver with the heuristic initial values $\varphi_j$ and include the best obtained energy level in Table~\ref{tab:Lshape} as well. Again, we conclude that for all cases except for $\kappa = 10$ the NN-enhanced approach yields a minimizing state of lower energy than what is possible with the classical FE iterative solver and heuristic initial values. For $\kappa = 10$ the same state is obtained. The difference in the energy states computed by the classical FE method and the hybrid method is also confirmed by comparing the densities, $|u|^2$, as shown in Figure \ref{fig:minimizers_Lshape}.

Summarizing the two experiments, we conclude that the results achieved with the hybrid-approach are robust with respect to the random seeds used in the NN stage. Furthermore, the \texttt{GL-F1} seems to be reliable within the limitation of the experiments conducted. Therefore, the architecture of \texttt{GL-F1} can be considered as a good candidate for the NN stage in the hybrid approach when applying the method to other GL problems.

\section{Conclusion}

In this work, we proposed a neural network approach as well as a neural network-enhanced (hybrid) finite element method to compute minimizers of the Ginzburg--Landau energy \eqref{eq:energy_functional}.
One of the main challenges in finding global minimizers is to choose good initial values for the minimization process as the results are highly sensitivity with respect to this choice.
We have shown that neural networks can be used to approximate Ginzburg--Landau minimizers, often leading to already satisfactory energy levels. 
However, using these results as initial values or predictors for the finite element solver 
we can achieve even better results. 
Leveraging the computational advantage of neural networks exploiting massive parallelization on GPUs, we have further demonstrated that even for a short offline phase,
compared to the computational time of the finite element solver, our neural network-enhanced approach can lead to minimizers of lower energy compared to a pure finite element approach with empirically chosen initial values. 
The proposed hybrid method does not necessarily imply faster convergence in terms of iteration steps, but we are able to find minimizers of lower energy, which we were not able to find with the pure finite element approach.

While the current neural network architecture achieves promising results and provides a good initial validation of the approach, we recognize that additional refinements and optimizations could further enhance its performance. However, these are beyond the scope of this study.

\section*{Acknowledgments}

B.~D\"orich and R.~Maier acknowledge funding from the Deut\-sche Forschungsgemeinschaft (DFG, German Research Foundation) -- Project-ID 258734477 -- SFB 1173. C.~D\"oding acknowledge funding from the Deutsche Forschungsgemeinschaft (DFG, German Research Foundation) under Germany's Excellence Strategy – EXC-2047/1 – 390685813. M.~Crocoll acknowledges funding by the Deutsche Forschungsgemeinschaft (DFG, German Research
Foundation) -- Project-ID 545165789. 
Parts of this work were conducted during the authors' stay at the
Hausdorff Research Institute for Mathematics funded by the Deutsche Forschungsgemeinschaft (DFG, German Research Foundation) under Germany's Excellence Strategy – EXC-2047/2 – 390685813.
\\
Further, we would like to thank Tim Buchholz for his valuable input on the MPI implementation of the FE solver and the research group IANM3 at KIT for including us in their computing project.

The authors gratefully acknowledge the computing time provided on the high-performance computer HoreKa by the National High-Performance Computing Center at KIT (NHR@KIT). This center is jointly supported by the Federal Ministry of Education and Research and the Ministry of Science, Research and the Arts of Baden-Württemberg, as part of the National High-Performance Computing (NHR) joint funding program (https://www.nhr-verein.de/en/our-partners). HoreKa is partly funded by the German Research Foundation (DFG).

This work is supported by the Helmholtz Association Initiative and Networking Fund on the HAICORE@KIT partition.

\end{document}